\documentclass[review]{elsarticle}

\usepackage{hyperref}
\usepackage{amsfonts}
\usepackage{verbatim}
\usepackage{amsbsy}

\newcommand{\be}{\begin{equation}}
\newcommand{\ee}{\end{equation}}
\newcommand{\bea}{\begin{eqnarray}}
\newcommand{\eea}{\end{eqnarray}}
\newcommand{\mbfx}{{\mathbf{x}}}

\newcommand{\mbfv}{{\mathbf{v}}}

\newcommand{\mbft}{{\mathbf{t}}}
\newcommand{\mbfn}{{\mathbf{n}}}

\newcommand{\hole}{{hole}}
\newcommand{\corner}{{c}}

\journal{Thin-walled structures}

\begin{document}

\begin{frontmatter}

\title{A weighted extended B-spline solver for bending and buckling of
stiffened plates}



\author{Joris C.G. Verschaeve}

\address{Universitet i Oslo, Oslo, Norway}
\ead{joris@math.uio.no}

\begin{abstract}
  The weighted extended B-spline method [H\"ollig (2003)] is applied to
  bending and buckling problems of plates with different shapes
  and stiffener arrangements. The discrete equations
  are obtained from the energy contributions of the
  different components constituting the system by means of the Rayleigh-Ritz approach. 
  The pre-buckling or plane stress is computed by means 
  of Airy's stress function. A boundary data extension
  algorithm for the weighted extended B-spline method
  is derived in order to solve for inhomogeneous 
  Dirichlet boundary conditions. A series of benchmark
  tests is performed touching various aspects 
  influencing the accuracy of the method.
\end{abstract}


\end{frontmatter}


\section{Introduction}

The growing need for larger container ships \cite{Economist}
led to renewed interest in computational methods for plate bending and
plate buckling in the maritime industry. One of the main
challenges in the construction of 
modern container vessels is to provide a sufficient
ultimate strength of the structure while keeping the material
usage minimal. The development of more accurate and faster computational
methods is one aspect in helping the
industry obtaining these goals.\\[1ex]

The physics and methods for plate bending and buckling problems
with stiffener arrangements are treated and reviewed in 
\cite{Timoshenko1961,VentselKrauthammer2001,Chajes1974,Bedair2013}. 
There is a relatively long history 
on the numerical treatment of the plate bending and buckling equations,
for which the finite element method
\cite{Ciarlet1978,PowellSabin1977,BarikMukhopadhyay1998,LiApplegarthBullBettessBondThompson1997}, 
the boundary element method 
\cite{CostaBrebbia1985,ManolisBeskosPineros1986,TanakaBercin1998},
the finite strip method \cite{HoTham1994} and the 
quadrature element method \cite{ZhongPanYu2011}
probably account for the most widespread methods. Another class
of methods are the sometimes called semi-analytic methods
based on Navier's sine solution. Among these
we find the method by \cite{BrubakHellesland2007b} 
for rectangular plates and stiffener arrangements
and the one by \cite{DjelosevicTepicTanackovKostelac2013} for buckling
of plates of more complex geometries. 
As a method, which is more related to the present work, 
B-splines have been in use for plate 
deformation, vibration and buckling problems,
for example for rectangular domains
\cite{MizusawaKajitaNaruoka1980,WuWuHuang2008,Sherar2004Thesis,YangChenZhangHe2013}
or in connection with the isoparametric stripe approach
for more complicated
geometries \cite{AuCheung1993,EccherRasmussenZandonini2006}.
A general trend in computational solid mechanics 
is the 
integration between CAD and structural analysis, which has
led to the usage of B-splines in connection with
NURBS for the isogeometric approach. 
The isogeometric approach with NURBS has also been used
for plate bending problems \cite{Raknes2011Master,Kim2013Thesis,BeiraoDaVeigaBuffaLovadinaMartinelliSangalli2012,ShojaeeaIzadpanahaValizadeha2012,LiZhangZheng2013,RealiGomez2015,LeeKim2013}.\\[1ex]

On the other hand, the weighted extended B-spline method
\cite{Hoellig2003} follows a similar aim as 
the isogeometric approach with NURBS,
namely to facilitate the
integration between CAD and structural analysis. Whereas
the isogeometric approach is based on using
the same discretization for the structural analysis as for the CAD,
the weighted extended B-spline method aims to facilitate
integration of CAD and structural solver by describing the
boundary of the domain in an embedded fashion,
allowing for a flexible treatment of complicated geometries,
while leading to sparse matrices and being higher order accurate.
\\[1ex]

As such, a single example of bending of a clamped plate
has been solved in \cite{Hoellig2003} by means
of the weighted extended B-spline method.
In the present treatise this method
is applied to the bending and buckling problem of
Kirchhoff plates of various shapes with and without stiffeners. 
Using the energy formulation of the
system, the Rayleigh-Ritz approach is used
to obtain the discrete equations. 
The pre-buckling stress is computed by means of Airy's stress
function. A scheme to solve for inhomogeneous Dirichlet boundary
conditions in the framework of the weighted extended B-spline
method is derived in order to handle the traction boundary conditions
for Airy's stress function. The method is applied to a number of
benchmark cases. Thereby different issues affecting the accuracy,
such as discontinuities or
singularities, are discussed. 
As the method displays higher order accuracy,
it is particularly well adapted for eigenvalue problems \cite{Boyd2001}
as arising in plate buckling. \\[1ex]

The present work is organized as follows. In section
\ref{sec:physics}, the physical problem is presented.
The weighted extended B-spline method by \cite{Hoellig2003} is
briefly summarized in section \ref{sec:numerics}. In this
section, we shall also explain the adaptation of the
weighted extended B-spline method to the present
computation of the plate
bending and buckling problem and the pre-buckling stress.
Results for a number of benchmark cases are 
presented in section \ref{sec:results}. The present
treatise is concluded in section \ref{sec:conclusions}. 

\section{Physical problem} \label{sec:physics} 

\subsection{Bending and buckling of plate and stiffener}

\begin{figure}
\includegraphics[width=\linewidth]{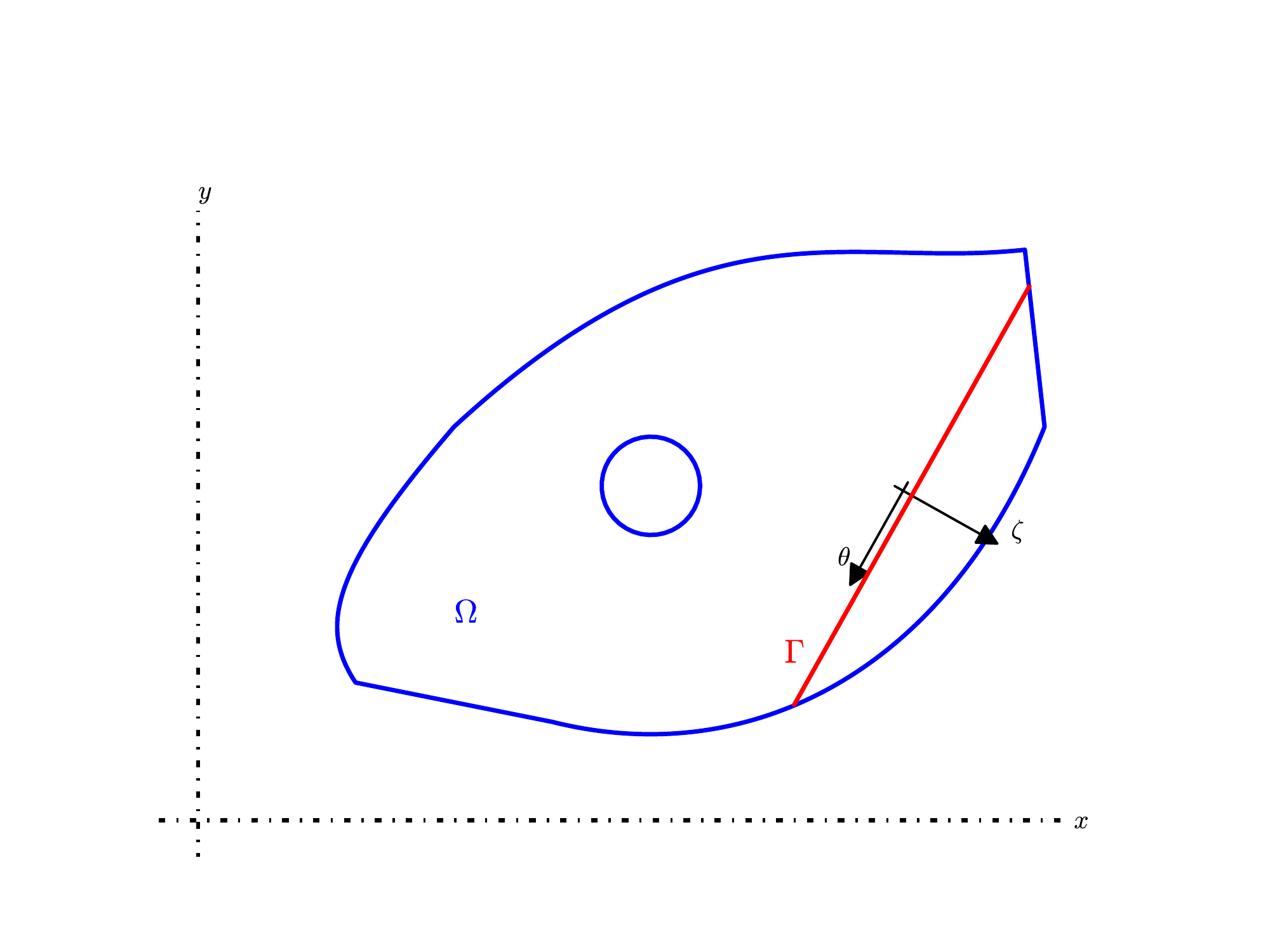}
\caption{A plate $ \Omega $ with a stiffener connected
to the plate at $ \Gamma $.}
\label{fig:sketchPlateBeam}
\end{figure}

In figure \ref{fig:sketchPlateBeam}, the geometry of a plate is sketched. 
It is described by
a domain $ \Omega \subset \mathbb{R}^2 $ not necessarily simply connected
(we allow for holes). At the straight line $ \Gamma $,
a stiffener is attached to the plate. 
For the domain $ \Omega $,
we use the coordinate system given by $ (x,y,z)$. The orientation of
the stiffener gives rise to a coordinate system 
defined by $ (\zeta, \eta, \theta) $, where $ \theta $
is the arclength. The $ \theta $-axis is parallel to the stiffener,
whereas $ (\zeta, \eta ) $ are in the plane perpendicular to the
axis of the stiffener with the $ \zeta $-axis lying in the $(x,y)$
plane of the plate. 
For a Kirchhoff plate, the displacement of the middle surface of
the plate is entirely described by the vertical displacement $ w $. 
The bending energy $ U_{plate} $  of a Kirchhoff plate is then given
by \cite{VentselKrauthammer2001}:
\be
U_{plate} = \frac{D}{2} 
\int \limits_{\Omega} 
\left( \nabla^2 w \right)^2 - 2 \left( 1 - \nu \right)
\left\{ 
\frac{\partial^2 w}{\partial x^2} 
\frac{\partial^2 w}{\partial y^2} 
- \left( \frac{\partial^2 w}{\partial x \partial y} \right)^2 \right\}
d \,  \Omega \label{eq:Energy1}
\ee
If a distribution $ p(x,y) $ of lateral forces per area 
is applied to the plate,
the work $ W_{lateral} $ done by these lateral forces has to be added to
the total energy,
\be
W_{lateral} = - \int \limits_{\Omega} p w \, d \Omega.
\ee
On the other hand, when considering buckling, in-plane forces $ \mathbf{T} $
are applied at the boundary $ \partial \Omega $ of the plate. 
These forces $ \mathbf{T} = ( T_x , T_y ) $ give rise to a stress field
$ \sigma $ in the plate 
which shall be called the pre-buckling stress
field and whose computation shall be described in section \ref{sec:planeStress}. 
Once the plate buckles, the work $ W_{buckling} $ done by the in-plane forces
can be computed by \cite{VentselKrauthammer2001}:
\be
W_{buckling} = 
\frac{1}{2} \int \limits_{\Omega} \sigma_{xx} \left( \frac{\partial w}{\partial x} \right)^2 + \sigma_{yy} \left( \frac{\partial w}{\partial y} \right)^2
+ 2 \sigma_{xy} \frac{\partial w}{\partial x} \frac{\partial w}{\partial y} \, d\Omega. \label{eq:Energy3}
\ee
Equations (\ref{eq:Energy1}-\ref{eq:Energy3}) describe all the relevant
energy contributions for plate bending and buckling without stiffeners.
Since at $ \Gamma $, cf. figure
\ref{fig:sketchPlateBeam}, a stiffener is attached to the plate, its
energy contribution needs to be added to the total energy of the system. 
In the present discussion,
the cross-section of the stiffener is assumed to stay constant
during deformation. For a class of stiffeners the torsional and
warping energies are negligible compared to the bending energy. 
The present discussion is restricted to these type of stiffeners. 
The displacement of the stiffener is thus in vertical direction
of an equal amount as the plate, since no detachment of the stiffener
is allowed. The bending energy of the stiffener can then be written as:
\be
U_{stiffener} = \frac{EI}{2} \int \limits_{\Gamma} 
\left( \frac{ d^2 w}{d \theta^2} \right)^2 \, d \theta, \label{eq:energyBeam}
\ee
where $ EI $ is the bending stiffness in vertical direction of the stiffener.
If opposing
forces of equal magnitude $ T_{s} $ 
are applied at both ends in axial direction of the stiffener,
the work $ W_{stiffener} $ due to axial shorting of the stiffener needs
to be accounted for. According to \cite{Chajes1974} (formula 5.42),
$ W_{stiffener} $ can in our case be written as:
\be
W_{stiffener} =  \frac{T_s}{2} \int \limits_{\Gamma}
\left\{ \left( \frac{d w}{d \theta} \right)^2 
+ r_0^2 \left( \frac{ d }{ d \theta} \frac{\partial w}{ \partial \zeta} \right)^2 - 2 \zeta_0 \frac{d w}{d \theta} 
\left( \frac{ d }{ d \theta} \frac{\partial w}{ \partial \zeta} \right) \right\}
\, d \theta, \label{eq:axialForcing}
\ee 
where $ r_0 $ is the radius of gyration 
and $ \zeta_0 $ the $ \zeta $-coordinate of the centroid of the cross section:
\be
A r_0 = \int \limits_A (\zeta^2 + \eta^2) \, d\zeta d\eta,
\quad A \zeta_0 = \int \limits_A \zeta \, d\zeta d\eta ,
\ee
where $ A $ is the area of the cross section. \\

Having defined all the relevant energy contributions, the problem
of plate bending by lateral loads and of plate buckling by in-plane
loads can be defined as follows. 

\paragraph{Plate bending} When lateral loads are applied, the plate
and the stiffener undergo bending and the total energy $ E_{bending} $ 
of the system can be written as:
\be
E_{bending} = U_{plate}  + U_{stiffener} + W_{lateral}.  \label{eq:energyBending}
\ee
The governing equations of the system can be obtained by 
variational minimization of (\ref{eq:energyBending}). For the
present system, we obtain:
\bea
D \nabla^4 w &=& p \quad \mbox{in} \quad \Omega \setminus \Gamma, \\
D \left[ \lim_{\epsilon \rightarrow 0^+} \frac{\partial^3 w}{\partial \zeta^3}
- \lim_{\epsilon \rightarrow 0^-} \frac{\partial^3 w}{\partial \zeta^3} \right]
+ EI \frac{\partial^4 w}{\partial \theta^4} &=& 0 \quad \mbox{on} \quad \Gamma,
\label{eq:jumpCondition}
\eea
in addition to boundary conditions on the boundary of $ \Omega $.  

\paragraph{Plate buckling} When in-plane loads are applied, the
plate and the stiffener undergo sudden buckling once a critical value
of the loading has been reached. The energy $ E_{buckling} $ of the problem
is given by:
\be
E_{buckling} = U_{plate} + U_{stiffener} 
+ \lambda \left( W_{buckling} + W_{stiffener} \right),
\label{eq:masterBuckling}
\ee
where $ \lambda $ is a parameter controlling the intensity of the
in-plane loading. In the framework of small displacement theory,
it plays the role of an eigenvalue allowing for non-trivial solutions
of $ w $ minimizing equation (\ref{eq:masterBuckling}). \\

\subsection{Pre-buckling stress} \label{sec:planeStress}

As mentioned above, the pre-buckling stress enters equation
(\ref{eq:Energy3}) and is unknown a priory. A
weighted extended B-spline
formulation for the plane stress problem in terms 
of the horizontal displacements can be found in \cite{Hoellig2003}.
However, in the present case, the boundary conditions
are given in terms of the in-plane forces
$ \mathbf{T} $ at the boundaries of $ \Omega $ which result
into boundary conditions for the stress tensor $ \sigma $. The
horizontal displacements at the boundaries are unknown a priori
and the formulation in \cite{Hoellig2003} cannot be applied straightforwardly. 
As such iterative solvers might be applied accounting for the undetermined
solid body motions. However, as we are employing a direct solver,
a more practical approach is to use the Airy stress function formulation. 
The present plane stress 
problem is therefore expressed in terms of Airy's stress function $ \Phi$. 
When dealing
with multiple connected domains, a necessary condition
for the Airy stress function to exist and to be smooth
is that the total
force and torque at the outer boundary and at the boundary of each hole vanishes
separately \cite{FosdickSchuler2003}.
In the present context, free holes, i.e. absence of traction,
are of main interest, such that the Airy stress function formulation
is applicable. \\

The stress components are expressed via Airy's stress function as: 
\be
\sigma_{xx} = \frac{\partial^2 \Phi}{\partial y^2} \quad
\sigma_{yy} = \frac{\partial^2 \Phi}{\partial x^2} \quad
\sigma_{xy} = -\frac{\partial^2 \Phi}{\partial y \partial x},
\label{eq:Airy}
\ee
and the traction boundaries can be written as 
\bea
T_x &=& \frac{d}{ds} \frac{\partial \Phi}{\partial y}  \label{eq:bcPhiA}\\
T_y &=& - \frac{d}{ds} \frac{\partial \Phi}{\partial x} \label{eq:bcPhiB},
\eea
where $ s $ is the arclength of the boundary. When the domain is 
multiply connected by $ n +1 $ boundaries, i.e. one outer boundary
and $ n $ inner boundaries (holes), the boundary conditions
(\ref{eq:bcPhiA}) and (\ref{eq:bcPhiB}) translate to 
the following expressions for each boundary $ \partial \Omega_i $, 
$ i= 0,\ldots, n $ \cite{JaswonSymm1977}:
\be
\left. \begin{array}{lcl} \Phi & = & \gamma_i + \alpha_i x +
\beta_i y + F_i \\
\frac{\partial \Phi}{\partial n} & = & \alpha_i n_x + \beta_i n_y + N_i
\end{array} \right\} \quad \mbox{for} \quad (x,y) \in \partial \Omega_i,
\ee
where $ \alpha_i $, $ \beta_i $ and $ \gamma_i $ are constants of
integration, $ n_x $ and $ n_y $ are the components of the
outward pointing normal on the boundary $ \partial \Omega_i $ and
the functions $ F_i $ and $ N_i $ are defined as follows:
\bea
F_i(s) &=& \int \limits_{0}^s \frac{d x}{d s}(s') f_x(s') + \frac{d y}{ds} (s')
f_y(s') \, ds' \\
N_i(s) &=& \frac{dy}{ds}(s) f_x(s) - \frac{dx}{ds}(s) f_y(s) \\
f_x(s) &=& - \int \limits_0^s T_y(s') \, ds' \label{eq:f1}\\
f_y(s) &=& \int \limits_0^s T_x(s') \, ds'\label{eq:f2}
\eea
In terms of Airy's stress function, the membrane energy $ E_m $ of the plate
becomes:
\bea
 E_m[\Phi] & = & \frac{h}{2E} \int \limits_{\Omega} \left( \Delta \Phi \right)^2
+ 2 ( 1 + \nu ) \left\{ \left( \frac{ \partial^2 \Phi}{\partial x \partial y} \right)^2  - \frac{ \partial^2 \Phi}{\partial x^2}
\frac{ \partial^2 \Phi }{ \partial y^2} \right\} \, d \Omega \\
&=& \frac{h}{2E} \int \limits_{\Omega} \left( \Delta \Phi \right)^2 \, d \Omega
- \frac{h}{2E} 
\sum \limits_{i=0}^n \int \limits_{\partial \Omega_i} f_x T_x + f_y T_y \,
d \sigma 
\label{eq:energyStress}
\eea
where $ f_x $ and $ f_y $ are defined in equations (\ref{eq:f1}) and (\ref{eq:f2}), respectively. The second term in equation (\ref{eq:energyStress})
is just a constant. An equivalent energy $ E^*_m $ for $ \Phi $ 
can therefore be written as: 
\be
E_m^*[\Phi] = \int \limits_{\Omega} \left( \Delta \Phi \right)^2 \, d \Omega,
\ee
which after variational minimization leads to 
\be
\Delta^2 \Phi = 0 \quad \mbox{in} \quad \Omega. \label{eq:govPhi}
\ee
Equation (\ref{eq:govPhi}) tells us that the method of choice
for the plane stress problem would be a boundary integral
or boundary element solver \cite{JaswonSymm1977}, since it would reduce
the two-dimensional problem into a one-dimensional one. The
resulting stress field could then be fed to the solver presented
in section \ref{sec:numerics} in order to solve the buckling
problem (\ref{eq:masterBuckling}). However, for illustration
purposes we shall use the present solver in order to solve
the plane stress problem for Airy's stress function, since
it shows how inhomogeneous Dirichlet boundary conditions can
be handled in the framework of the weighted extended B-spline method. \\

When given a multiple connected domain, with $ n+1 $
boundaries, we shall solve $ 4(n+1) $ times the following
decoupled problems:
\bea
\Delta^2 \Phi_{ik} &=& 0 \quad x \in \Omega \quad i = 0,\ldots, n \quad
k = 0,\ldots,3 \label{eq:preBuckling1}
\eea
with the following boundary conditions:
\bea
\Phi_{i0} &=& F_i, \quad \frac{\partial \Phi_{i0}}{\partial n} = N_i \quad x \in \partial \Omega_i
\label{eq:preBuckling2}\\
\Phi_{i1} &=& 1, \quad \frac{\partial \Phi_{i1}}{\partial n} = 0 \quad x \in \partial \Omega_i \\
\Phi_{i2} &=& x, \quad \frac{\partial \Phi_{i2}}{\partial n} = n_x \quad x \in \partial \Omega_i \\
\Phi_{i3} &=& y, \quad \frac{\partial \Phi_{i3}}{\partial n} = n_y \quad x \in \partial \Omega_i \\
\Phi_{ik} &=& 0, \quad \frac{\partial \Phi_{ik}}{\partial n} = 0 \quad x \in \partial \Omega_j \quad \forall j \neq i \quad k = 0,\ldots, 3. \label{eq:preBuckling6}
\eea
The linear combination
\be
\Phi = \sum \limits_{i=0}^n \Phi_{i0} + \gamma_i \Phi_{i1}
+ \alpha_i \Phi_{i2} + \beta_i \Phi_{i3} ,
\ee
satisfies the governing equation (\ref{eq:govPhi}) and the
boundary conditions (\ref{eq:bcPhiA}) and (\ref{eq:bcPhiB}). \\

In the framework of the boundary integral or boundary element method
\cite{JaswonSymm1977} constraints for 
the undetermined constants $ \alpha_i $, $ \beta_i $, $ \gamma_i $,
$ i = 0,\ldots, n $ can be formulated
exploiting compatibility between displacements and stresses. However,
for a volume based method (or area based in two dimensions)
as the present method, such an approach
is not straightforward. We shall instead use the minimum energy 
principle by writing:
\bea
E^*_m[\Phi] &=& E^*_m\left[ \sum \limits_{i=0}^n \Phi_{i0} + \gamma_i \Phi_{i1}
+ \alpha_i \Phi_{i2} + \beta_i \Phi_{i3} \right]\\
&=& E^*_m(\dots,\alpha_i,\beta_i,\gamma_i,\ldots)\\
&=& \sum \limits_{i=0}^n \sum \limits_{j=0}^n 
\int \limits_{\Omega} 
 \left( \Delta \Phi_{i0} + \gamma_i \Delta\Phi_{i1} + \alpha_i \Delta\Phi_{i2} 
+ \beta_i \Delta\Phi_{i3} \right) \nonumber \\
& & \quad \quad \left( \Delta \Phi_{j0} + \gamma_j \Delta\Phi_{j1} + \alpha_j \Delta\Phi_{j2} 
 + \beta_j \Delta\Phi_{j3} \right) \, d\Omega
\eea
The coefficients $ \alpha_i $, $ \beta_i $, $ \gamma_i $, $ i = 0,\ldots, i $
are then found by minimizing the function 
$ E^*_m(\dots,\alpha_i,\beta_i,\gamma_i,\ldots) $ under the constraint that
\bea
\int \limits_{\Omega} \Phi \, d\Omega & = & 0, \\
\int \limits_{\Omega} \frac{\partial}{\partial x}\Phi \, d\Omega & = & 0, \\
\int \limits_{\Omega} \frac{\partial}{\partial y}\Phi \, d\Omega & = & 0,
\eea
since any function of the form:
\be
\Phi + a + bx + cy,
\ee
where $ a $, $ b $ and $ c $ are arbitrary constants, satisfies
equations (\ref{eq:bcPhiA}), (\ref{eq:bcPhiB}) and (\ref{eq:govPhi}). \\[1ex]

\section{Numerical scheme} \label{sec:numerics}

This section presents briefly the main results of
the weighted extended B-spline method, presented in detail 
in \cite{HoelligReifWipper2001,Hoellig2003,HoelligHoernerHoffacker2010,HoelligHoernerPfeil2012}. The weighted extended B-spline method
is a finite element method based on B-splines combined with
an embedded description of
the boundary. The method can be considered a Cartesian grid method,
which avoids grid generation, since no body fitted mesh needs to be generated.
The vertical displacement $ w $ of the plate is expanded on the
weighted extended B-spline basis functions $ B_\mathbf{k} $ \cite{Hoellig2003}:
\be
w = \sum \limits_{\mathbf{i} \in \mathcal{I} } w_\mathbf{i} B_\mathbf{i}(x,y), \label{eq:basicExpansion}
\ee
where $ \mathbf{i} = ( i_x, i_y) $ is a two-dimensional index
which can be mapped onto a global index $ i$. The
set $ \mathcal{I} $ denotes the set of inner indices, which shall
be explained below.
The construction of the
weighted extended B-spline basis functions shall now be presented
along general lines, for details we refer to \cite{Hoellig2003}. \\

\begin{figure}
\includegraphics[width=\linewidth]{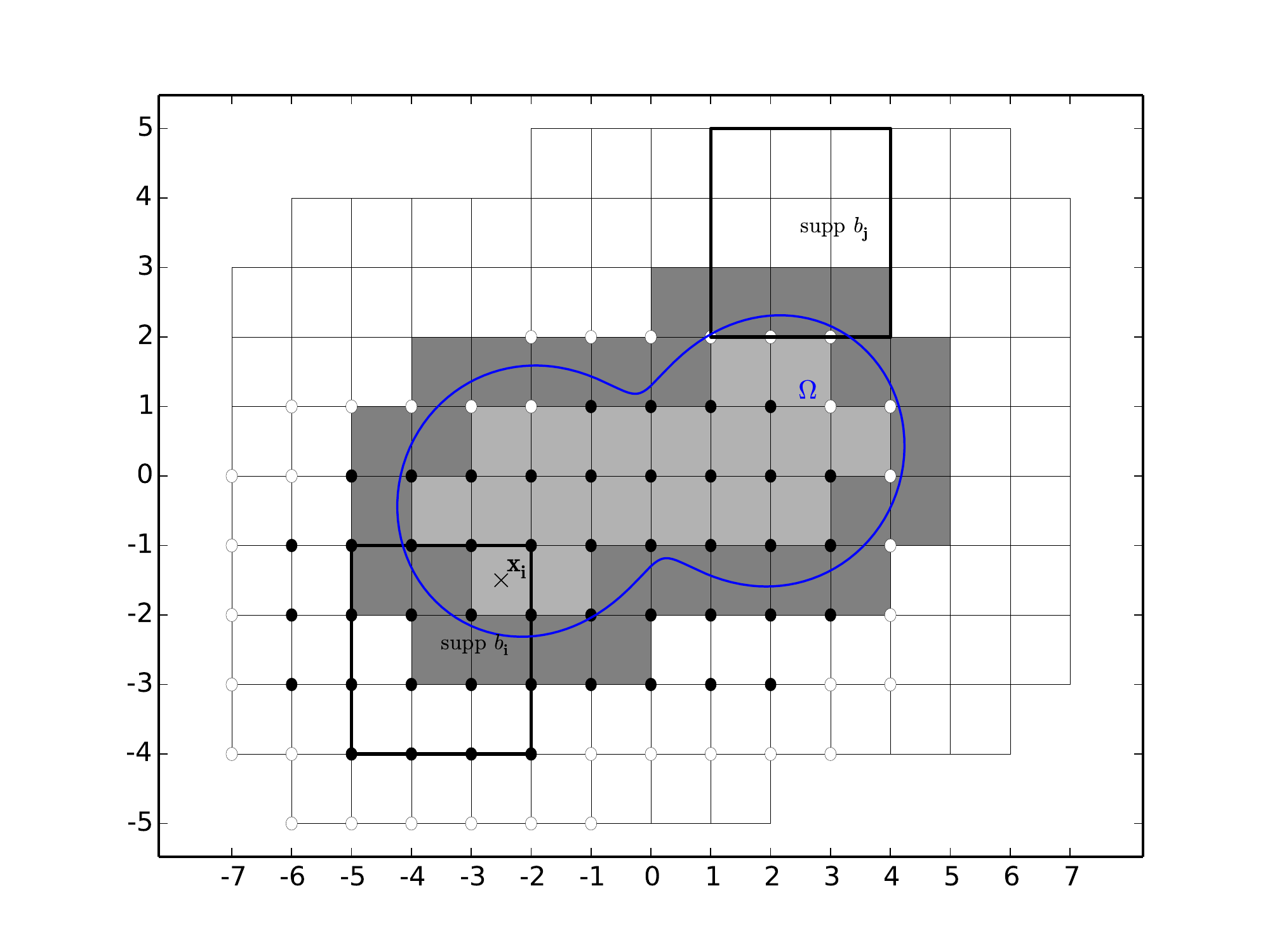}
\caption{Discretization of a domain $ \Omega $. The boundary of $ \Omega $ (blue line)
is embedded in a grid consisting of uniform rectangular cells. Cells
can either be cut by the boundary (boundary cells, dark gray color),
lie completely inside the domain (interior cells, light gray color)
or lie outside of the domain (exterior cells, white color).}
\label{fig:cellAndDomain}
\end{figure}

The discretization of a two dimensional domain $ \Omega $
is given by square cells of side length $ h$ as sketched in figure \ref{fig:cellAndDomain}.
Figure \ref{fig:cellAndDomain} has been drawn using figure 4.6 in \cite{Hoellig2003} as example.
The cells can be categorized into interior, boundary or exterior cells, depending on
whether they are entirely inside $\Omega $, they are cut by the boundary of $ \Omega $
or they are entirely outside of $ \Omega$.
The non-weighted basis function $ b_{\mathbf{k}} $ is a tensor
product of the one-dimensional B-spline basis function $ b^p_k $
in $ x $ and $ y $ direction of degree $ p $, cf. \cite{Hoellig2003}:
\be
b_{\mathbf{k}} (x,y) = b^p_{k_x}(x) b^p_{k_y}(y).
\ee
The indexing of the basis functions is
by means of the index of the lower left cell of their support. In figure
\ref{fig:cellAndDomain}, we plotted the supports of $ b_{\mathbf{i}} = b_{(-5,-4)} $
and $ b_{\mathbf{j}} = b_{(1,2)} $ in the case of $ p = 2 $. When the support
of $ b_{\mathbf{i}} $ contains at least one interior cell, it is labeled
as inner B-spline,
cf. \cite{Hoellig2003}. On the other hand, when the support of $ b_{\mathbf{j}} $
contains no interior cell but the intersection between the support and
$ \Omega $ is nonempty, the B-spline is labeled as outer B-spline. The
support of all other B-splines does not intersect $ \Omega $ and they 
can be discarded from the computations. The set of indices of
all inner B-splines is denoted $ \mathcal{I} $. For the set
of indices of all outer B-splines we use $ \mathcal{J} $. The notation
follows the one used in \cite{Hoellig2003}. 

In order to conform the basis functions to the boundary, the basis
functions $ b_{\mathbf{k}} $ are weighted by a weight function $ \omega(x,y) $:
\be
B^{\omega}_{\mathbf{k}} (x,y) = \omega(x,y) b_{\mathbf{k}}(x,y).
\ee
The role of the weight function is to modify the behavior of
the basis functions at the boundary such 
that the basis functions satisfy the boundary conditions at the
boundary. In particular, if $ d $ is the distance to a point
on the boundary,
in a region sufficiently close to this point, a Taylor expansion
of the weight function might be written as
\be
w \approx w_0 + w_1 d + w_2 d^2 + \ldots
\ee
If $ w_0 $ is chosen to be zero
on some part of the boundary, all basis functions
will be zero at this part of the boundary. This way simply supported
boundary conditions can be imposed for the plate. When in addition
we require that $ w_1 = 0 $ on some part of the boundary,
clamped boundary conditions can be modeled for this part of the
boundary. When $ w_0 $ is nonzero,
the plate is free at the boundary. This way the different boundary conditions
for a plate can be introduced in a relatively straightforward way. \\

A major difficulty of the embedded boundary description 
is the appearance of basis functions $ B^{\omega}_{\mathbf{k}} $
with arbitrarily small support, where the support of 
$ B^{\omega}_{\mathbf{k}} $ is given by the intersection of the support of
$ b_{\mathbf{k}} $ and $ \Omega $. 
For example the support of the 
basis function $ B^{\omega}_{(4,1)} $ in figure \ref{fig:cellAndDomain}
is only a small fraction of the area of a cell. 
This leads to extremely ill conditioned stiffness matrices. 
This problem has been solved by \cite{Hoellig2003} 
by means of the extension algorithm. 
As the inner basis functions $ B^{\omega}_{\mathbf{i}} $
dispose of a support at least the size of a cell, the
aim of the extension algorithm is to extend the
inner basis functions $ B^{\omega}_{\mathbf{i}} $ by
the outer basis functions $ B^{\omega}_{\mathbf{j}} $
and thus creating a new set of basis functions $ B_{\mathbf{i}} $. 
The formula for the weighted extended B-spline of degree $ p $ is
given by (see box 4.9 on page 48 in \cite{Hoellig2003}):
\be
B_{\mathbf{i}} = \frac{\omega}{\omega(\mathbf{x}_\mathbf{i}) } 
\left[ b_{{\mathbf{i}}} + \sum \limits_{\mathbf{j} \in \mathcal{J}(\mathbf{i}) } 
e_{\mathbf{i},\mathbf{j}} b_{\mathbf{j}} \right] \label{eq:webSpline}
\ee
where $ \mathcal{J}(\mathbf{i}) $ is the set of all outer
indices adjacent to the inner index $ \mathbf{i} $. 
In chapters 4 and 8 in \cite{Hoellig2003}, 
the procedure of computing $ \mathcal{J}(\mathbf{i}) $
given an inner index $ \mathbf{i} $ is explained in detail.
The 
precise definition of the coefficients
$ e_{\mathbf{i}\mathbf{j}} $ in front of the outer basis functions
is given in box (4.9) on page 48 in \cite{Hoellig2003}. 
The weight function in equation (\ref{eq:webSpline}), 
which depends on the geometry and the boundary
conditions of the problem is central to the definition 
of the weighted extended B-spline basis functions. 
The weight functions used in the
present treatise will be defined in section \ref{sec:results}
for each case considered.\\

\subsection{Bending and buckling of plate and stiffener}

Once we have defined
a set of weighted extended B-spline functions for a given
geometry and discretization, we can apply the Rayleigh-Ritz approach
to find a discretization of problems (\ref{eq:energyBending}) and (\ref{eq:masterBuckling}). 
If $ \mathbf{w} = (\ldots, w_{\mathbf{i}}, \ldots) $ is the column vector
containing all expansion coefficients of (\ref{eq:basicExpansion}), the discretization
of problem (\ref{eq:energyBending}) might be written as
\be
\mathbf{A} \mathbf{w} = \mathbf{b}, \label{eq:discrete1}
\ee
whereas for problem (\ref{eq:masterBuckling}), we obtain
\be
\mathbf{A} \mathbf{w} = \lambda \mathbf{B} \mathbf{w}, \label{eq:discrete2}
\ee
The elements of the stiffness matrix $ \mathbf{A} $ are given by:
\bea
\lefteqn{A_{\mathbf{k}\mathbf{l}} } \nonumber \\
&=& D \int \limits_{\Omega} \left(\frac{\partial^2 B_{\mathbf{k}}}{\partial x^2} 
+ \nu \frac{\partial^2 B_{\mathbf{k}}}{\partial y^2} \right) 
\frac{\partial^2  B_{\mathbf{l}} }{\partial x^2}
+ 2(1 - \nu) \frac{\partial^2 B_{\mathbf{k}} }{\partial x \partial y} 
\frac{\partial^2 B_{\mathbf{l}}}{\partial x \partial y} 
+ \left(\frac{\partial^2 B_{\mathbf{k}} }{\partial y^2} 
+ \nu \frac{\partial^2 B_{\mathbf{k}} }{\partial x^2} \right) 
\frac{\partial^2 B_{\mathbf{l}}}{\partial y^2} \, d \Omega \nonumber \\
& & \quad + EI \int \limits_{\Gamma} \frac{\partial^2 B_{\mathbf{k}}}{\partial \theta^2}
\frac{\partial^2 B_{\mathbf{l}}}{\partial \theta^2}\, d s
\eea
The elements of the second member $ \mathbf{b} $ due to the lateral loading $ p $
are defined by:
\be
b_{\mathbf{k}} = \int \limits_{\Omega} B_{\mathbf{k}} p \, d \Omega
\ee
The matrix $ \mathbf{B} $ accounting for the shortening of 
plate and stiffener has the following elements:
\bea
B_{\mathbf{k}\mathbf{l}} &=& 
\int \limits_{\Omega} \sigma_{xx} \frac{\partial B_{\mathbf{k}}}{\partial x} 
\frac{\partial B_{\mathbf{l}}}{\partial x} 
+ \sigma_{yy} \frac{\partial B_{\mathbf{k}}}{\partial y} 
\frac{\partial B_{\mathbf{l}}}{\partial y} 
+ \sigma_{xy} \left( \frac{\partial B_{\mathbf{k}}}{\partial x} 
\frac{\partial B_{\mathbf{l}}}{\partial y} +
\frac{\partial B_{\mathbf{k}}}{\partial y} 
\frac{\partial B_{\mathbf{l}}}{\partial x} \right) \, d\Omega \nonumber \\
&& + T_s \int \limits_{\Gamma} \frac{\partial B_{\mathbf{k}}}{\partial \theta} 
\frac{\partial B_{\mathbf{l}}}{\partial \theta}  + r_0^2 
\frac{\partial^2 B_{\mathbf{k}}}{\partial \theta \partial \eta} 
\frac{\partial^2 B_{\mathbf{l}}}{\partial \theta \partial \eta}
- \eta_0 \left( \frac{\partial B_{\mathbf{k}}}{\partial \theta } 
\frac{\partial^2 B_{\mathbf{l}}}{\partial \theta \partial \eta}
+ 
\frac{\partial^2 B_{\mathbf{k}}}{\partial \theta \partial \eta}
 \frac{\partial B_{\mathbf{l}}}{\partial \theta } \right) \, d s \nonumber
\eea
In general, the solution of problems (\ref{eq:energyBending}) and
(\ref{eq:masterBuckling}) is only $ \mathcal{C}^2 $ continuous across the
stiffener location. This can be seen from the jump condition 
(\ref{eq:jumpCondition}) accounting for an additional term in 
the force balance of the plate due to the stiffener.
However the basis functions $ B_{\mathbf{k}} $ will in general
be of higher smoothness than the solution at the stiffener location. 
As we shall see in section \ref{sec:results}, this leads to reduced
convergence rates. \\

Due to the finite support of the basis functions, the matrices
in (\ref{eq:discrete1}) and (\ref{eq:discrete2}) are sparse.
In \cite{Hoellig2003,HoelligHoernerPfeil2012} special
iterative schemes exploiting the sparseness are devised. 
Since the problems in the present work are of moderate size,
we shall use a sparse LU solver to invert (\ref{eq:discrete1}).
The sparse LU solver has been written by Tim Davis and can be
downloaded freely at \cite{DavisWebsite}. 
For problem (\ref{eq:masterBuckling}) the generalized
eigenvalue solver from Lapack \cite{Lapack} has been used. \\

\subsection{Pre-buckling stress}  \label{sec:numericsStress}

\begin{figure}
\includegraphics[width=\linewidth]{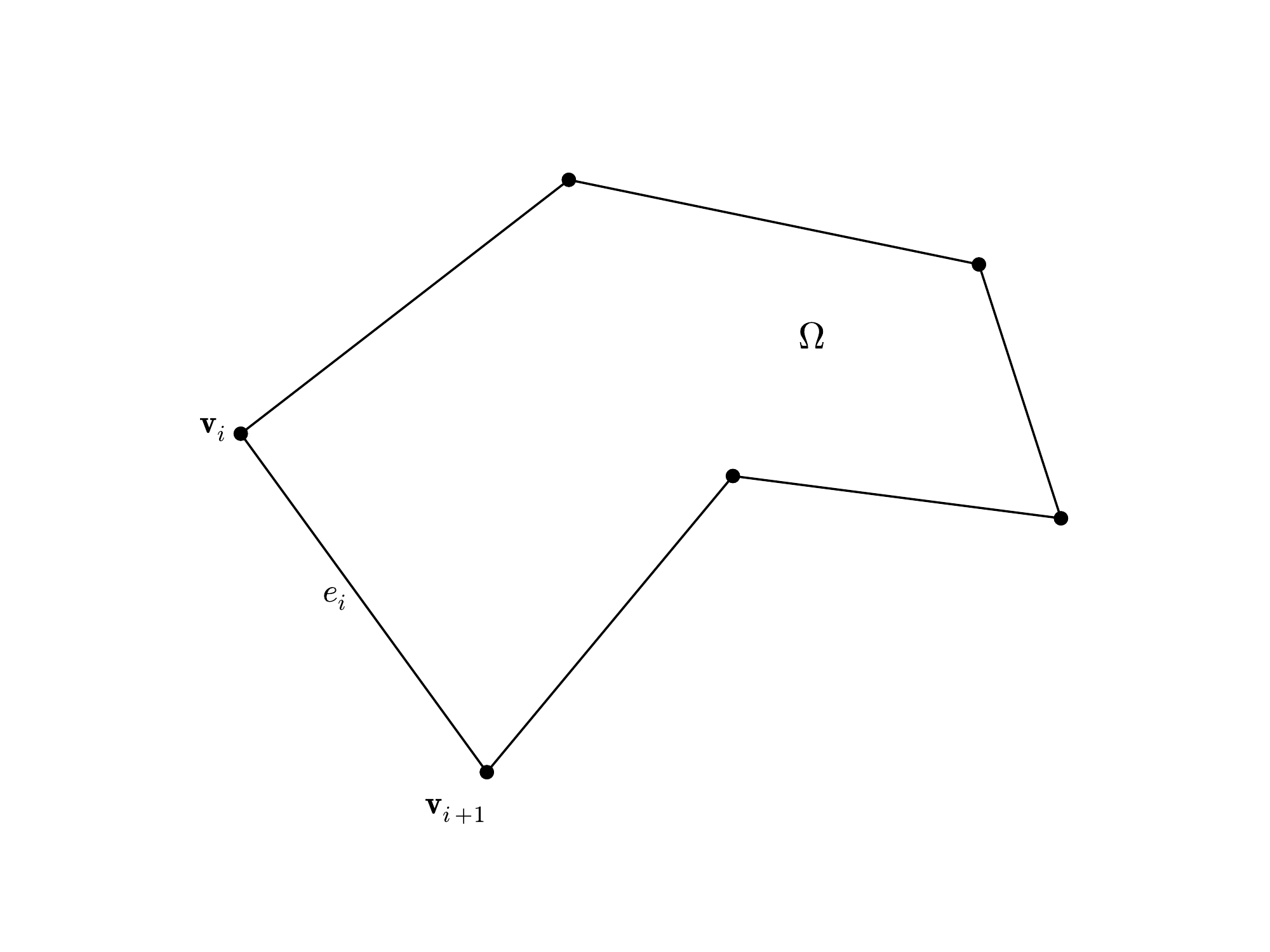}
\caption{A simple (non-convex) polygon with vertices $ \mbfv_{i} $ and
edges $ e_{i} $.}
\label{fig:polygon}
\end{figure}

Concerning the pre-buckling stress problem, equations (\ref{eq:preBuckling1}
-\ref{eq:preBuckling6}),
each sub-problem can be written as:
\bea
\Delta^2 \hat{u} = 0 \quad \mbox{in} \quad \Omega \label{eq:extension1}\\
\hat{u} = f \quad \mbox{on} \quad \partial \Omega \\
\frac{\partial \hat{u}}{\partial n}  = g \quad \mbox{on} \quad \partial \Omega
\label{eq:extension2} \\
\eea
This will be solved by formulating an extension $ \tilde{u} $ of the boundary
data, such that  $ \tilde{u} $ satisfies the boundary conditions: 
\bea
\tilde{u} = f \quad \mbox{on} \quad \partial \Omega \\
\frac{\partial \tilde{u}}
{\partial n}  = g \quad \mbox{on} \quad \partial \Omega. 
\eea
Problem (\ref{eq:extension1}-\ref{eq:extension2}) can then be cast into an equivalent problem:
\bea
\Delta^2 u = -\Delta^2 \tilde{u} \quad \mbox{in} \quad \Omega \label{eq:inhom1}\\
u = 0 \quad \mbox{on} \quad \partial \Omega \\
\frac{\partial u}{\partial n}  = 0 \quad \mbox{on} \quad \partial \Omega. \label{eq:inhom3}
\eea
The solution $ \hat{u} $ is then simply:
\be
\hat{u} = u + \tilde{u}. 
\ee
For each $ \Phi_{ik} $ in equation (\ref{eq:preBuckling1}),
we need thus to 
solve a problem of type (\ref{eq:inhom1}-\ref{eq:inhom3}), where
the boundary conditions (\ref{eq:preBuckling2}-\ref{eq:preBuckling6}) 
enter $ \tilde{u} $ in the
second member of (\ref{eq:inhom1}). Since we are using a LU solver
for the inversion of the stiffness matrix a repeated
solution of equations (\ref{eq:inhom1}-\ref{eq:inhom3})
for different second members can be performed efficiently. The difficulty
lies in formulating a function extension formulation for $ \tilde{u} $.
A general transfinite interpolation algorithm as in \cite{RvachevSheikoShapiroTsukanov2001}
introduces (numerical) singularities at the vertices of the domain,
cf. figure \ref{fig:polygon},
even in case of convex polygons,
which might not only reduce the order of convergence
of the scheme but, in addition, the integral of $ \Delta \tilde{u} $ might
not exist in $ \Omega $ rendering the scheme unusable. For convex polygons, a transfinite
interpolation leading to smooth interpolants as long as the boundary data
is smooth is given in \cite{VaradyRockwoodSalvi2011} 
and \cite{SalviVaradyRockwood2014}.
As, we are not primarily interested in transfinite interpolation, but
only want to find an extension function to the boundary data, we propose
an extension algorithm in \ref{sec:extension}
which works also for non-convex polygons, c.f. figure \ref{fig:polygon}, as long
as they are simple. The algorithm might also be extended to cases,
where more general segments replace the edges of the polygon. 

\clearpage
\section{Results} \label{sec:results}

A series of test cases for plate bending and buckling is 
computed for different geometries. In section \ref{sec:annular}
bending and buckling of an annular plate is considered
and compared to reference cases in the literature. The same
is done for the geometry of a rectangular plate with and
without holes in section \ref{sec:rectangular}. 
As an example of a more complicated shape, we consider bending and
buckling of a polygonal plate with holes in section \ref{sec:polygonal}. 

\subsection{Annular plate} \label{sec:annular}

\begin{figure}
\includegraphics[width=\linewidth]{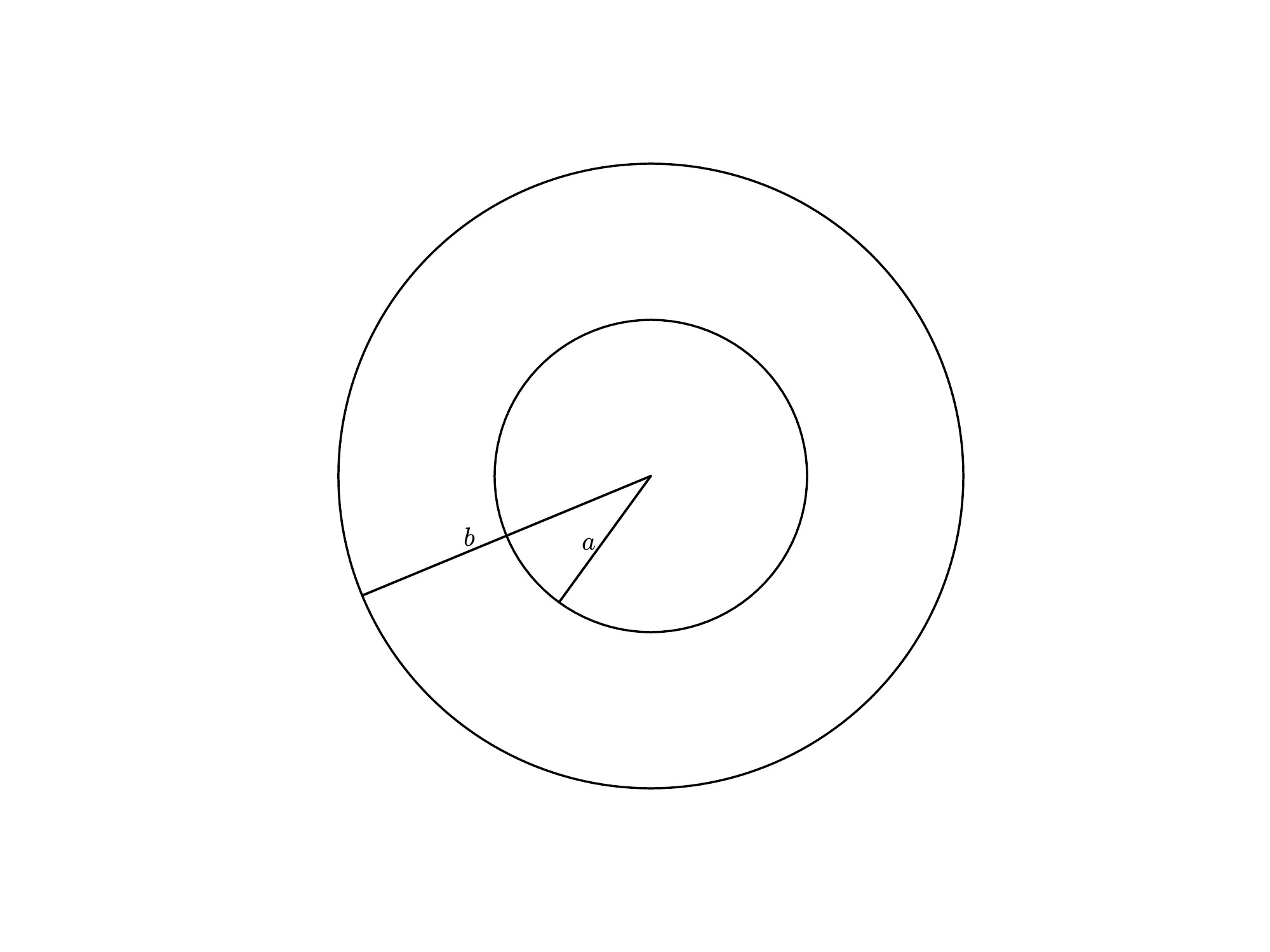}
\caption{Geometry of the annular 
plate with clamped outer boundary and free inner boundary.
A constant lateral forcing is applied onto the plate.}
\label{fig:annularPlateGeometryDeformation}
\end{figure}

The geometry of the annular plate is sketched in figure \ref{fig:annularPlateGeometryDeformation}.
At the outer boundary with radius $ b $ clamped boundary conditions
are applied, whereas the inner boundary with radius $ a$ is free. 
The weight function $ \omega $, entering the definition of
the basis functions, equation (\ref{eq:webSpline}), is for this
case defined by:
\be
\omega(x,y) = \left( x^2 + y^2 - b^2 \right)^2. 
\ee

\subsubsection{Bending of an annular plate by a lateral load}

\begin{figure}
\includegraphics[width=\linewidth]{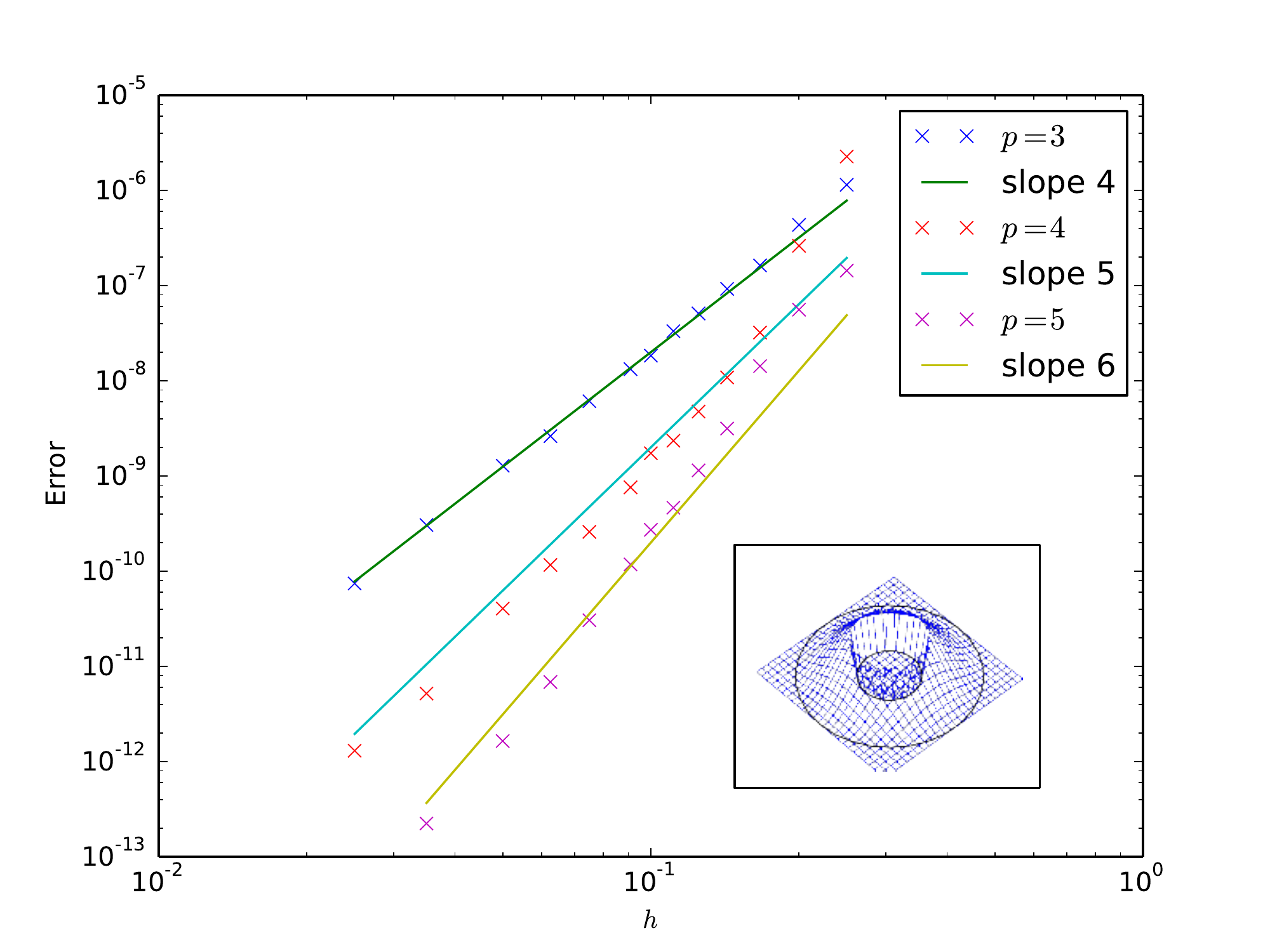}
\caption{Convergence of the error of the numerical solution with respect to
the cell size $ h $ for the bending problem of an annular plate
under constant lateral loading using splines of different degree $ p $.}
\label{fig:convergence:BendingAnnularPlate}
\end{figure}

When considering plate bending for a 
constant lateral loading $ p(x,y) = p_0 $, an analytic solution 
can be found \cite{VentselKrauthammer2001}:
\be
w(r) = c_1 \log(r) + c_2 r^2 \log(r) + c_3 r^2 + \frac{p_0}{64D} r^4,
\label{eq:annularPlateAnalyticalSolution}
\ee
where $ c_1 $, $ c_2 $ and $ c_3 $ are functions of $ a,b,\nu $ and $ D $.
Their exact definition is rather laborious but can easily be obtained
by symbolic computation using the boundary conditions at the inner and
outer boundary. 
Choosing $ a = 0.5345 $, $ b = 1.5432 $,
$ p_0 = 1.74586 $, $ \nu = 0.3 $ and $ D = 1.234 $, we apply the present solver
for different degrees $p $ and resolutions $ h $ to the above problem. A surface plot
of the solution is embedded in the convergence plot
in figure \ref{fig:convergence:BendingAnnularPlate}, showing
that the maximum deformation is obtained at the inner boundary and that
the solution is axisymmetric. The present solver does, however, not employ
the axisymmetry of the problem. It solves it on a rectangular Cartesian grid. 
In figure \ref{fig:convergence:BendingAnnularPlate}, we observe that the embedded boundary 
treatment of the weighted extended B-spline solver is indeed higher order accurate. 
For sufficiently small $ h$,the error converges with the theoretical order of convergence for
smooth problems, which is $ p+1$. 

\subsubsection{Buckling of an annular plate by in-plane loads}

\begin{figure}
\centering
\includegraphics[width=\linewidth]{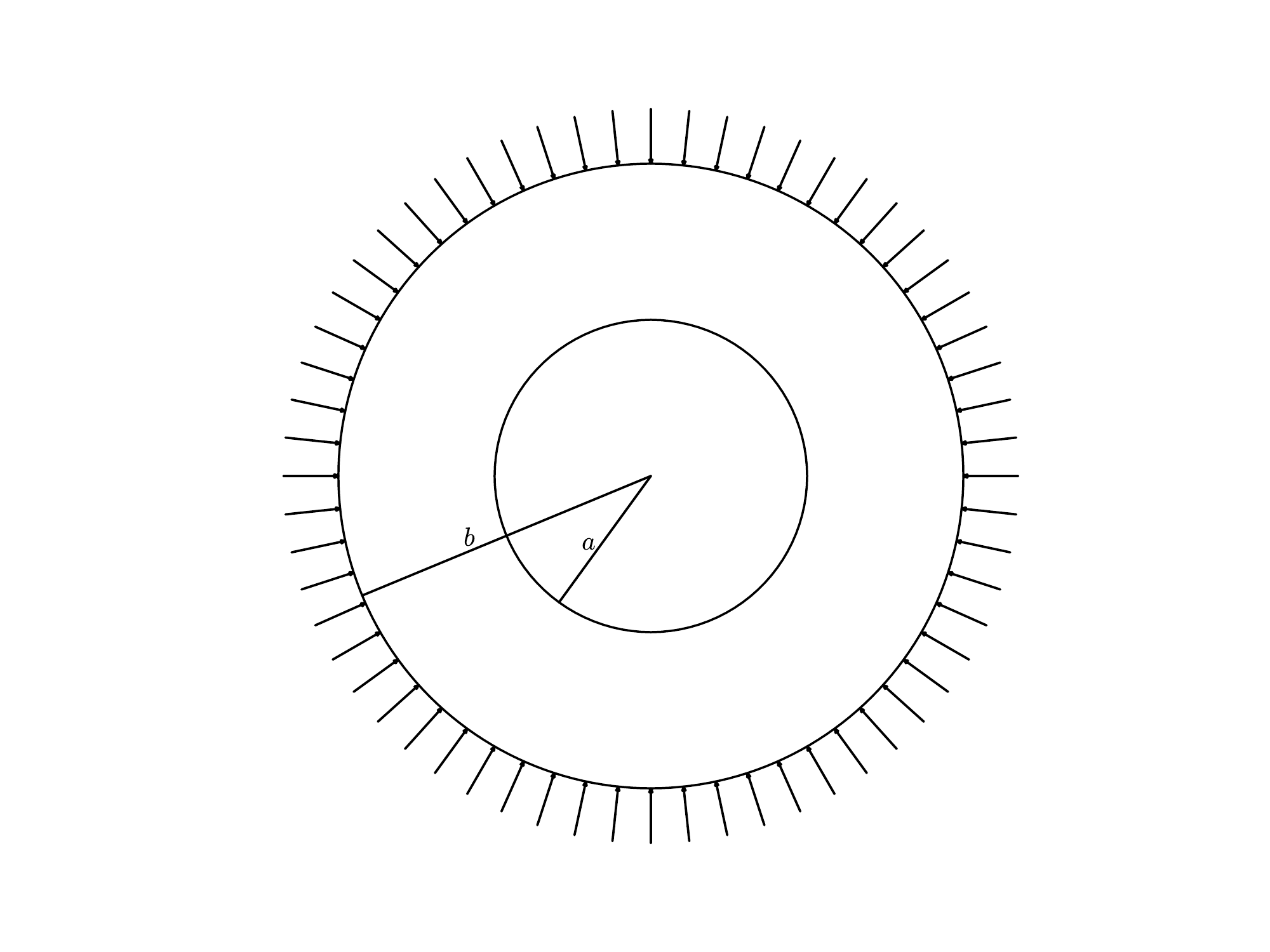}
\caption{Geometry of the annular 
plate with free inner boundary and clamped
outer boundary.
A constant in plane forcing in normal direction is applied onto the
outer boundary.}
\label{fig:annularPlateGeometry}
\end{figure}

\begin{figure}
\centering
\includegraphics[width=\linewidth]{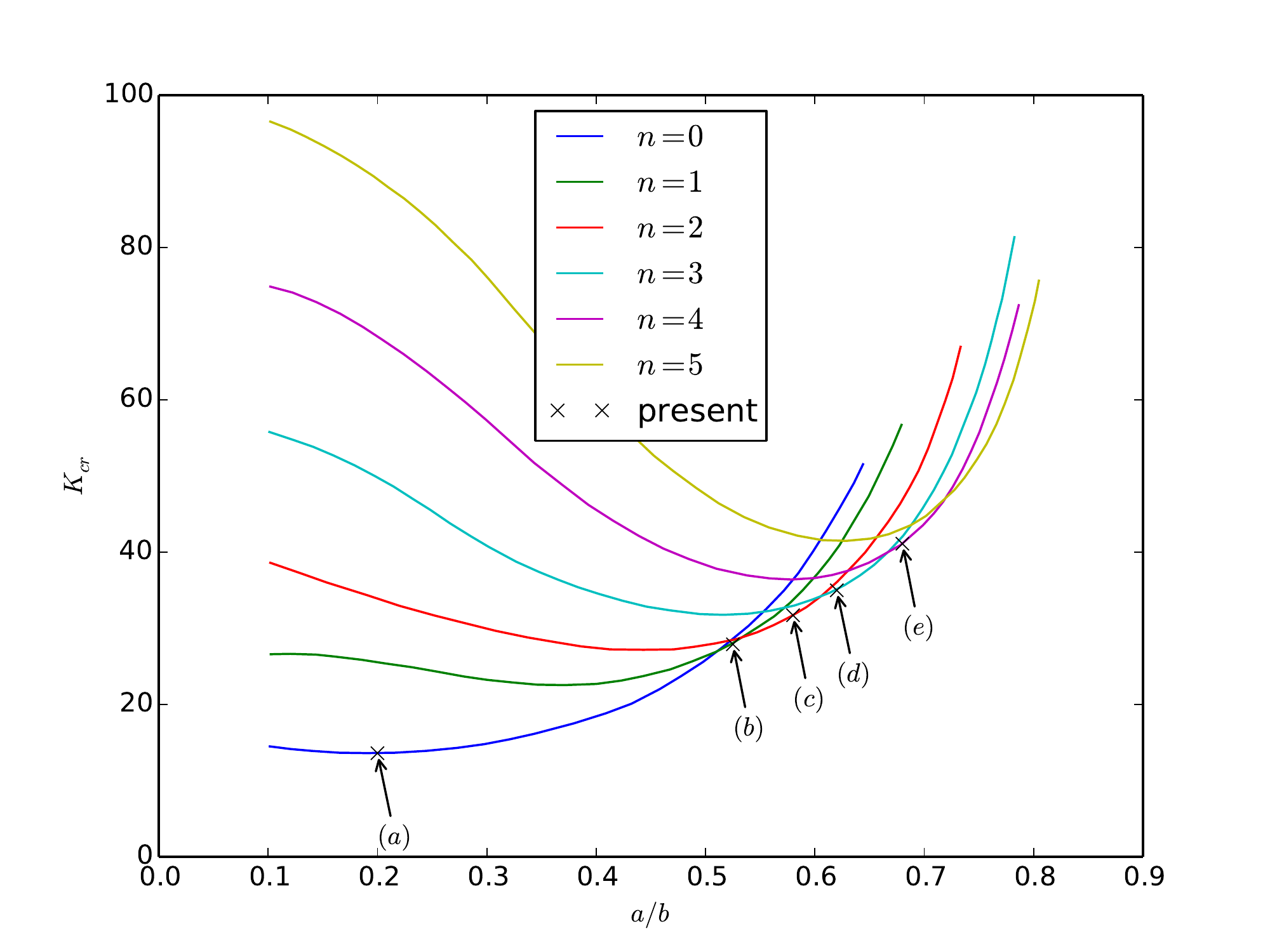}
\caption{Critical buckling stresses for an annular plate with
free inner and clamped outer boundary condition, as sketched in figure
\ref{fig:annularPlateGeometry}. The lines for the modes $ n = 0,\ldots 5 $ 
are scanned form figure 3 in \cite{JillellaPeddieson2013}. The buckling modes corresponding
to the examples $ (a),(b),(c),(d) $ and $ (e) $ computed by means of 
the present method are plotted in figure \ref{fig:bucklingModesAnnularPlate}.}
\label{fig:bucklingAnnularPlate}
\end{figure}

\begin{figure}
\centering
\includegraphics[width=0.3\linewidth]{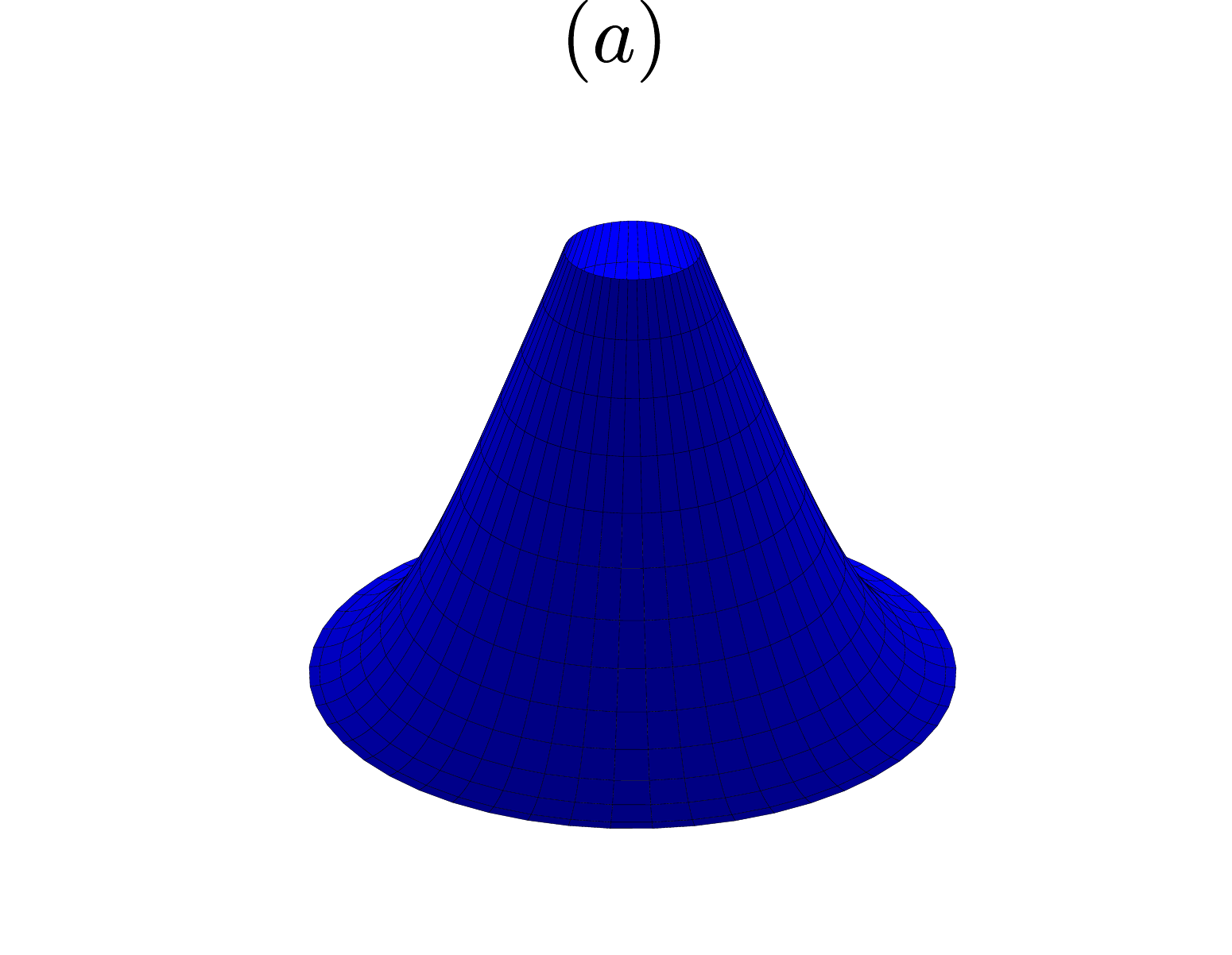}
\includegraphics[width=0.3\linewidth]{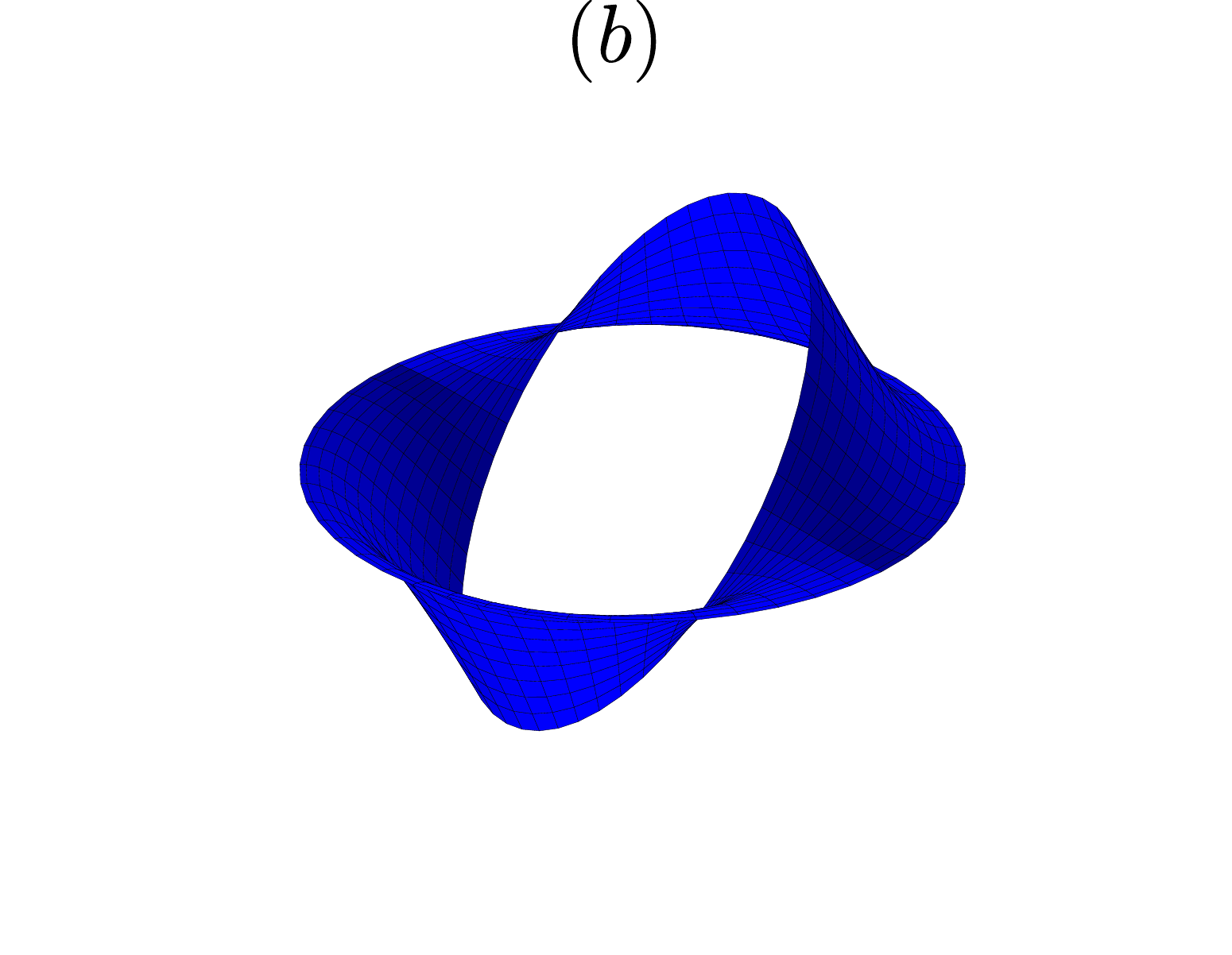}
\includegraphics[width=0.3\linewidth]{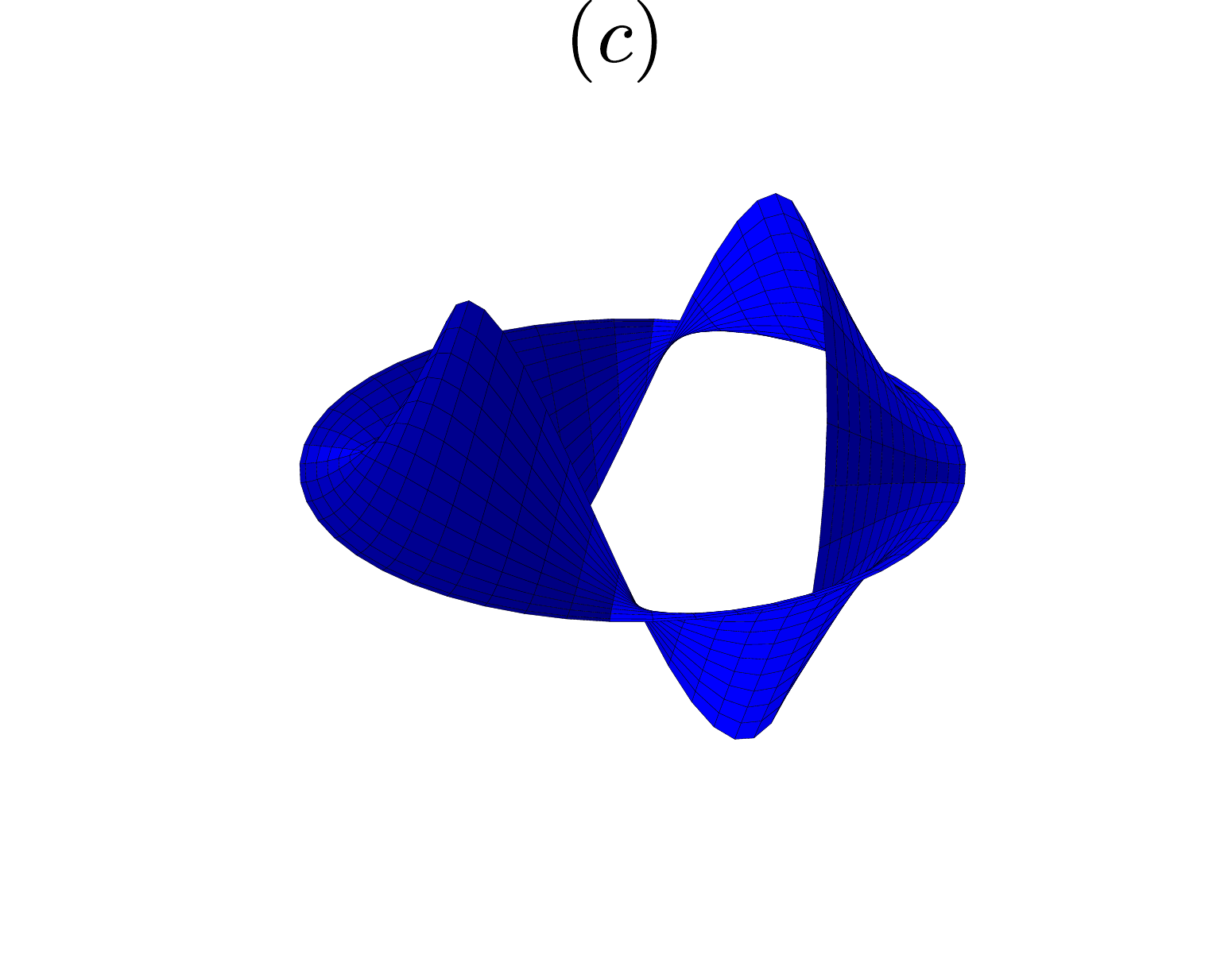}\\
\includegraphics[width=0.3\linewidth]{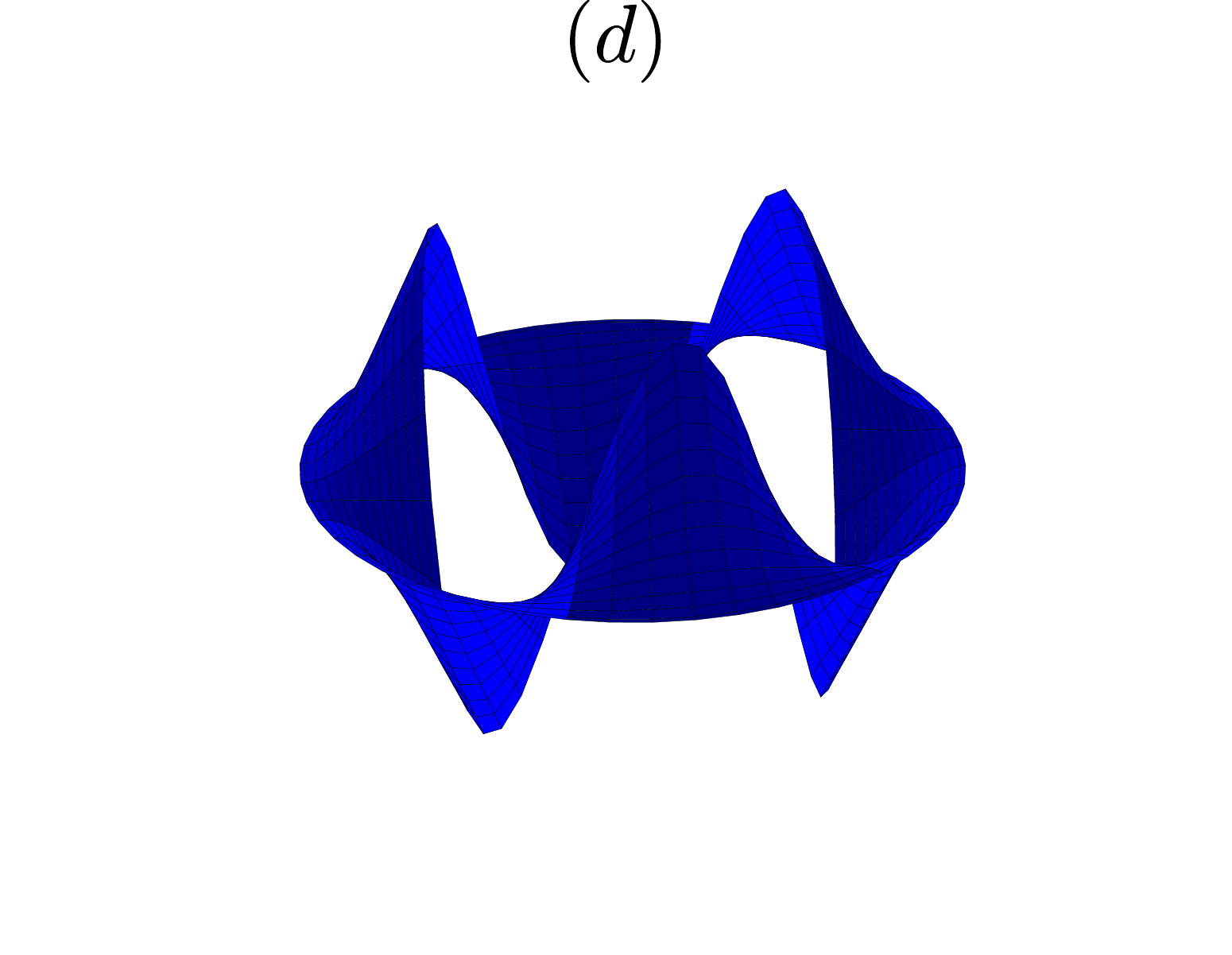}
\includegraphics[width=0.3\linewidth]{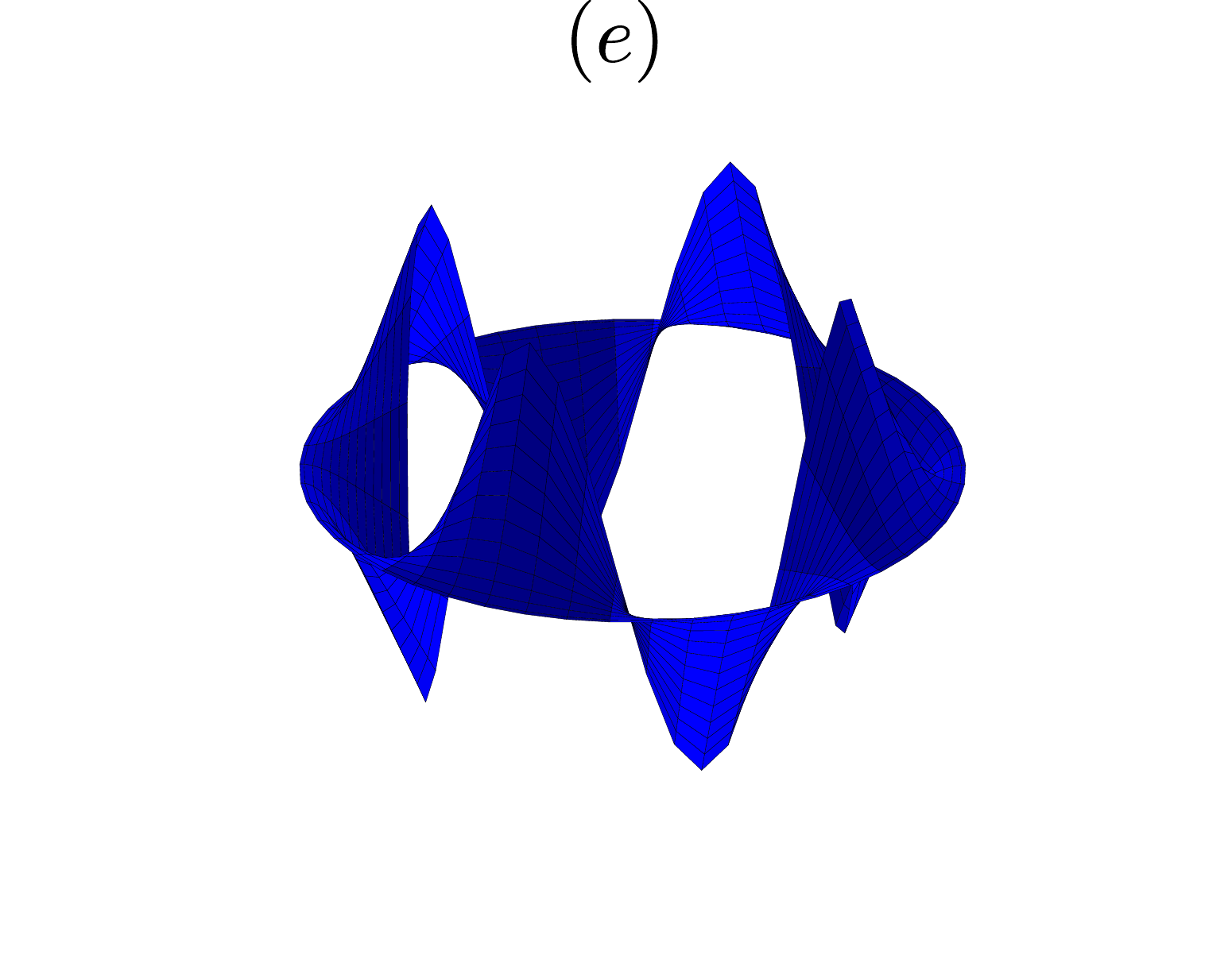}
\caption{Buckling modes of an annular plate with
free inner and clamped outer boundary condition for the cases
denoted by $ (a),(b),(c),(d) $ and $ (e) $ in figure
\ref{fig:bucklingAnnularPlate}.}
\label{fig:bucklingModesAnnularPlate}
\end{figure}

In this section, instead of a lateral loading, a radial compression force is applied at
the outer boundary of an annular plate,
cf. figure \ref{fig:annularPlateGeometry}.
The case of buckling of an annular plate has found some attention
in the literature, cf. \cite{Coman2007,ComanBassom2009,JillellaPeddieson2013} 
and references therein.
Concerning the choice of boundary conditions
at the outer boundary, the more interesting case is for clamped 
boundary conditions. For this relatively simple geometry, the 
stress components of $ \sigma $ can be found analytically
and are given by \cite{Coman2007}:
\bea
\sigma_{rr} &=& - \frac{ 1 - a^2/r^2}{1 - a^2/b^2}\\
\sigma_{\phi\phi} & = & - \frac{ 1 + a^2/r^2}{1 - a^2/b^2} \\
\sigma_{r\phi} & = & 0
\eea
in polar coordinates $ (r,\phi)$. For the case $ \nu = 0.3 $, the buckling
stresses for the different modes have been plotted in figure 3
in \cite{JillellaPeddieson2013}. These graphs have been digitized 
from figure 3
in \cite{JillellaPeddieson2013} and are plotted together with
the present results in
figure \ref{fig:bucklingAnnularPlate}. We solve equation (\ref{eq:masterBuckling})
for five aspect ratios $ a/b $, cases $ (a) $ to
$ (e) $ in table \ref{tab:bucklingAnnularPlate},
using the present weighted extended B-spline solver. The outer radius
has been set to $ b = 2.28 $. The critical stresses $ \lambda_c $ from
equation (\ref{eq:masterBuckling}), are nondimensionalized by
means of the $ K $-value:
\be
K = \frac{\lambda_c}{D/b^2},
\ee
and are reported in table \ref{tab:bucklingAnnularPlate} for 
two resolutions $ h = 0.2 $ and $ h = 0.1 $. As can be seen from
the values found, approximately five significant digits are
identical for both resolutions. The reference values have been
read from figure 3 in \cite{JillellaPeddieson2013}. 
The buckling modes corresponding to the cases in table \ref{fig:bucklingAnnularPlate}
are plotted in figure \ref{fig:bucklingModesAnnularPlate}.
For increasing ratio $ a/b $, the solution displays an increasing
number of buckles in angular direction. From this relatively simple test
case, we can conclude that the present solver based
on an embedded description of the boundary provides accurate solutions
for buckling problems. 


\begin{table} 
\centering
\caption{K-values for the buckling of an annular plate for different ratios $ a/b $ and resolutions $ h$.}
\label{tab:bucklingAnnularPlate}

 \begin{tabular}{c|c|r|r|r|c} 
case & $a/b$ & $K (h = 0.2) $ & $K (h = 0.1) $ & $K^{ref} $ \\
 \hline \hline 
(a) & 0.2 & 13.60464726310300 & 13.60389138752100 & 13.6  \\
(b) & 0.525 & 27.90167587045000 & 27.90151625370600 & 28.0  \\
(c) & 0.58 & 31.71513802673000 & 31.71489313775400 & 31.7  \\
(d) & 0.62 & 34.99334267376700 & 34.99266753385800 & 35.1  \\
(e) & 0.68 & 41.10888319580700 & 41.10806291507800 & 41.1 
\end{tabular}
\end{table}

\clearpage
\subsection{Rectangular plate} \label{sec:rectangular}

For a rectangular plate defined by $ [-a/2,a/2] \times [-b/2,b/2] $,
the weight function in equation (\ref{eq:webSpline}) is
given for simply supported boundary conditions by
\be
\omega(x,y) = \left( x + \frac{a}{2} \right) \left( \frac{a}{2} - x \right)
\left( y + \frac{b}{2} \right) \left( \frac{b}{2} - y \right).
\label{eq:weightRectangle}
\ee
When dealing with clamped boundary conditions, we have to square the 
right hand side of (\ref{eq:weightRectangle}). In section
\ref{sec:bendingRectangularPlateAlignedStiffener}, where we rotate
the plate with respect to the grid, the weight function
(\ref{eq:weightRectangle}) is rotated with the plate. In section 
\ref{sec:stressRectangularPlateHole}, the pre-buckling 
stress in a plate with central hole is computed. For this case,
we need in addition a simple weight for the boundary of
the round hole:
\be
\omega_{hole}(x,y) = x^2+y^2 - \frac{d^2}{4},
\ee
where $ d $ is the diameter of the hole. 

\subsubsection{Bending of a rectangular plate with aligned stiffener}
\label{sec:bendingRectangularPlateAlignedStiffener}

\begin{figure}
\includegraphics[width=\linewidth]{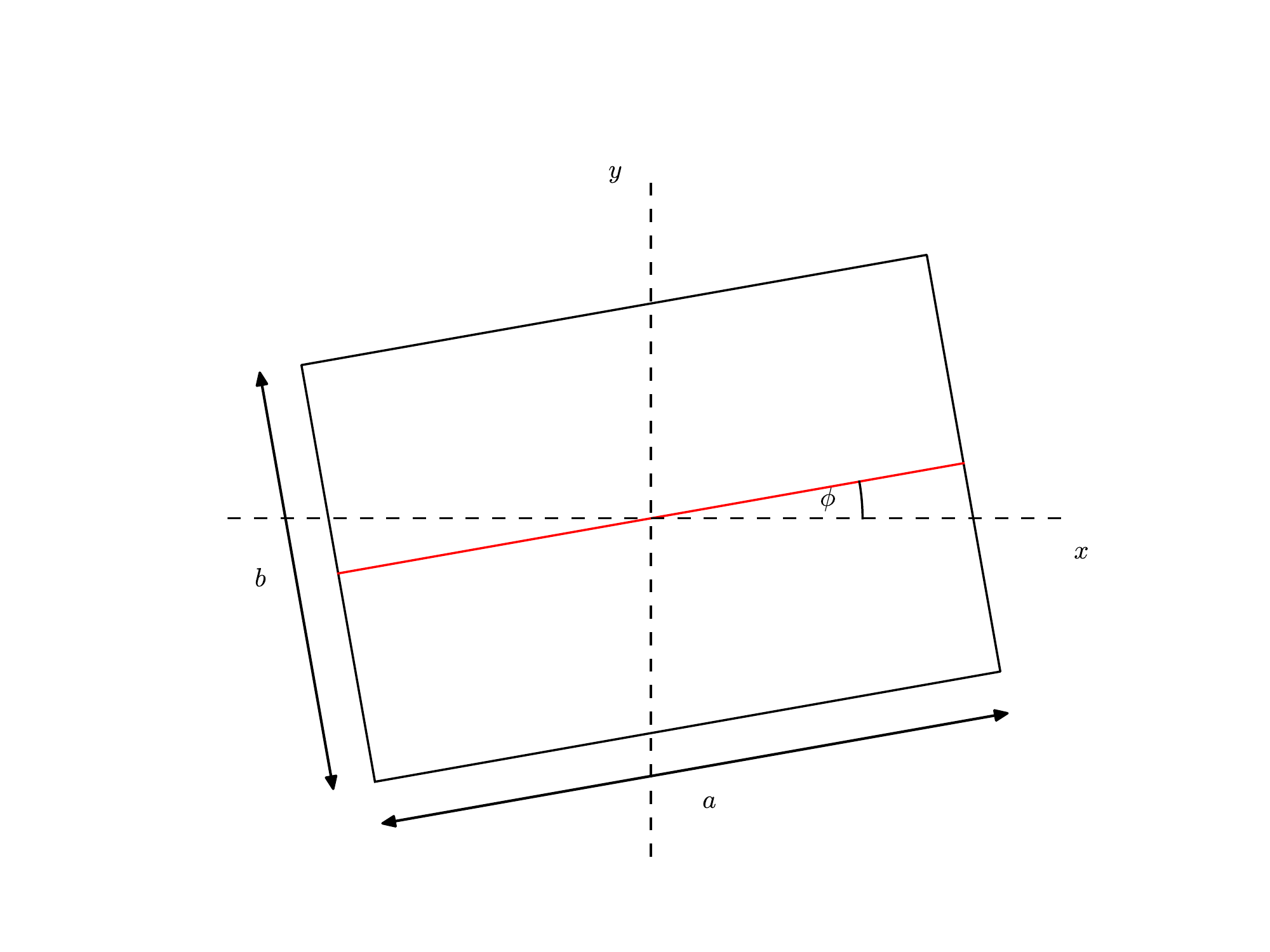}
\caption{Geometry of a rectangular plate reinforced by a single stiffer (red color) parallel to the lower and upper boundary of the plate.
The plate is rotated by an angle $ \phi $
with respect to the coordinate axes $ (x,y) $ defining the orientation of 
the computational cells. A uniform lateral loading is applied onto the plate. }
\label{fig:sketchRotatedStiffenedPlate}
\end{figure}

\begin{figure}
\includegraphics[width=\linewidth]{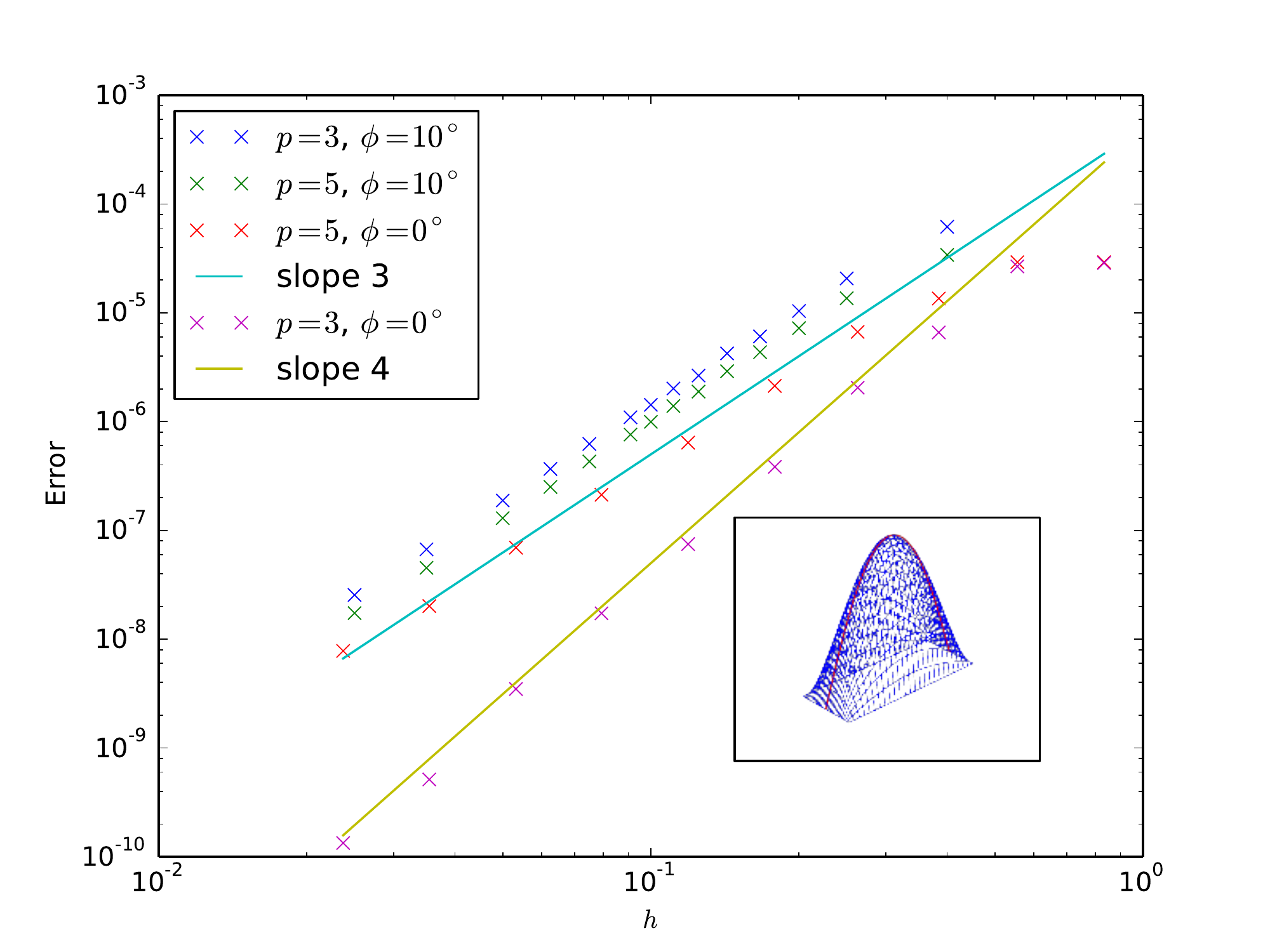}
\caption{Convergence of the error of the numerical solution 
with respect to the cell size $ h $ for 
the bending problem sketched in figure \ref{fig:sketchRotatedStiffenedPlate}.
The degree of the B-splines is given by $ p = 3,5 $ and the inclination
is either $ \phi = 0^{\circ} $ or $ \phi = 10^{\circ} $.} 
\label{fig:convergencePlateRotated}
\end{figure}

A rectangular plate with side lengths $ a $ and $ b$ is supported by a
stiffener located at the middle line parallel to the lower 
and upper boundary of the plate,
cf. figure \ref{fig:sketchRotatedStiffenedPlate}. 
In order to test the present embedded boundary description, we rotate the
plate by an angle $ \phi $ such that the geometry of the plate is no longer
aligned with the grid. The plate is loaded by a constant loading $ p_0 $
and we use simply supported boundary conditions at the edges of the domain.
The parameters of the simulation are given by:
\be
a = 5 \quad b = 1 \quad D = 1.234 \quad EI = 150 \quad p_0 = 1.543. 
\ee
The angle $ \phi $ is set to $ 10^\circ $. For comparison,
we perform a simulation with aligned geometry, meaning that we take $ \phi =  0^\circ $. 
A reference solution, computed on a fine grid, is used in order to compute
the error of the numerical solution. In figure \ref{fig:convergencePlateRotated},
the error of the solution is plotted in function of the cell size $ h$ for
splines of degree $ p = 3,5$. For the aligned case $ \phi = 0^\circ $, we
observe that the convergence rate for the splines of degree $ p = 3 $ is
approximately four. However, for splines of degree $ p = 5 $, we observe a smaller rate
at approximately three. This can be explained by the jump condition (\ref{eq:jumpCondition})
at the beam. As the third normal derivative of the solution is discontinuous,
we expect a reduction of the convergence rate, since the third
derivatives of the basis functions with degree $ p = 5 $ are continuous.
The third derivatives of the basis functions with $ p = 3 $ are however
not continuous at the grid lines and we therefore see the theoretical
order of convergence, namely four. When rotating the domain by $ 10^\circ $,
the stiffener is not aligned with the grid lines anymore and we observe
a reduced convergence rate of approximately three for both degrees. 

\subsubsection{Buckling of a rectangular plate with aligned stiffener}

\begin{figure}
\includegraphics[width=\linewidth]{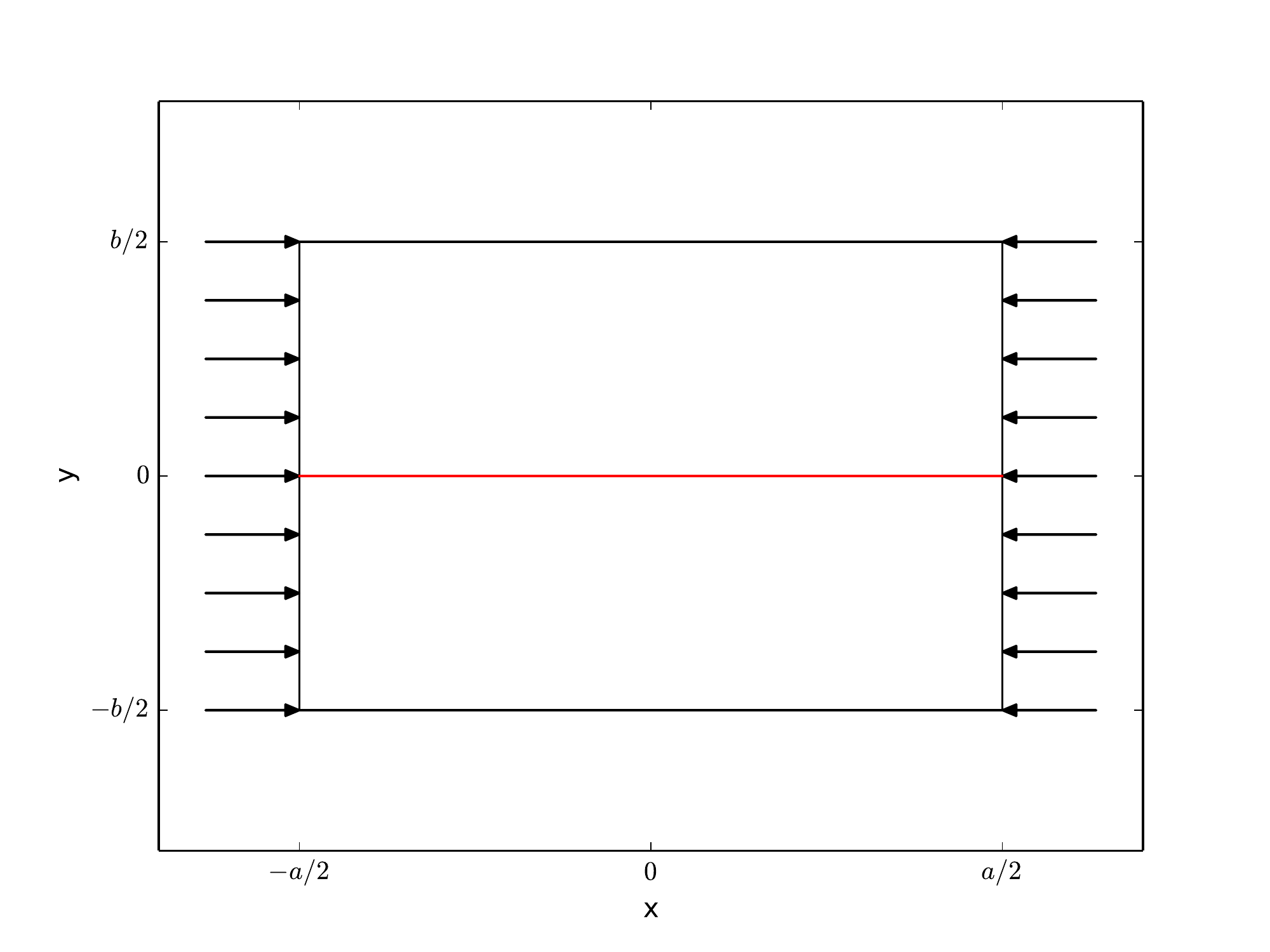}
\caption{Geometry of a rectangular plate reinforced by a single stiffener
(red color) parallel to the lower and upper boundary.
Uniform in-plane loading is applied at the left and right boundaries.}
\label{fig:sketchBucklingStiffenedPlate}
\end{figure}

\begin{figure}
\includegraphics[width=\linewidth]{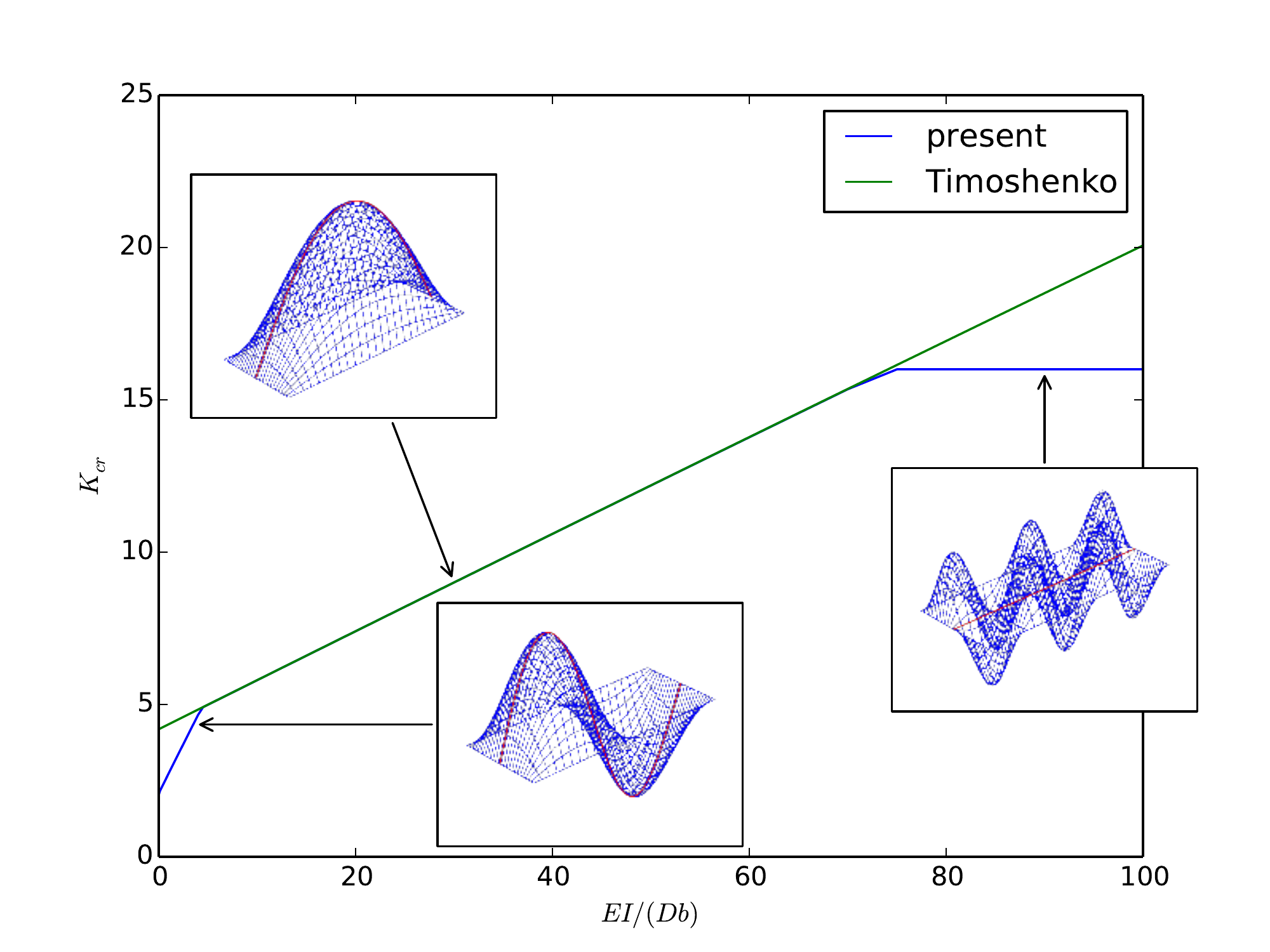}
\caption{Buckling of a rectangular plate reinforced by a single
stiffener parallel to the lower and upper boundary of the plate 
with dimensions $ a/b = 5/2 $ and 
aspect ratio between the cross section areas of the rip and
the plate given by $ A/(bh) = 1/2 $.}
\label{fig:bucklingRectangularPlateRip}
\end{figure}

This case corresponds to example 9.11 in \cite{Timoshenko1961}. As before,
a simply supported plate is supported by a stiffener along its middle line. However,
this time an in-plane force normal to the left and right boundary is
applied, cf. figure \ref{fig:sketchBucklingStiffenedPlate}.
The stress field for this case is simply given by:
\be
\sigma_{xx} = -1, \quad \sigma_{yy} = \sigma_{xy} = 0. 
\ee
The magnitude $ T_s $ of the axial forcing for the stiffener,
cf. equation (\ref{eq:axialForcing}), is given as in \cite{Timoshenko1961} by:
\be
T_s = b \delta, \quad \delta = \frac{A}{bh}, 
\ee
where $ A $ is the cross section area of the stiffener. Since the solution
is symmetric around the stiffener, we have $ \frac{ \partial w}{\partial \zeta} = 0 $ 
on $ \Gamma $, meaning that in the formula for the axial shortening (\ref{eq:axialForcing})
only the first term is nonzero. As in \cite{Timoshenko1961}, we introduce
the other nondimensional parameters controlling the system by:
\be
\beta = \frac{a}{b} \quad \gamma = \frac{EI}{bD}. 
\ee
The $K$-value for this system is defined as in \cite{Timoshenko1961}:
\be
K = \frac{\lambda_{cr} b^2 h}{\pi^2 D}.
\ee
An analytic value for $ K $ is given by the smaller root of
equation (k) in example 9.11 in \cite{Timoshenko1961}. 
This value is plotted as a green line
in figure \ref{fig:bucklingRectangularPlateRip} for the present choice of
\be
\beta = \frac{5}{2} \quad \delta = \frac{1}{2}. 
\ee
Formula (k) in \cite{Timoshenko1961},
is, however, only valid in the case of the plate buckling into a single buckle. As
can be observed from figure \ref{fig:bucklingRectangularPlateRip}, for
weak stiffeners, the buckled plate displays two buckles, whereas for strong
stiffeners, only the plate displays buckling, whereas the stiffener 
stays straight, similar to what \cite{BrubakHellesland2007} (figure 6) 
have found before. 
For a single buckle the present method gives buckling stresses
remarkably close to the approximate formula by \cite{Timoshenko1961}.

\subsubsection{Bending of a square plate with non-aligned stiffener}

\begin{figure}
\includegraphics[width=\linewidth]{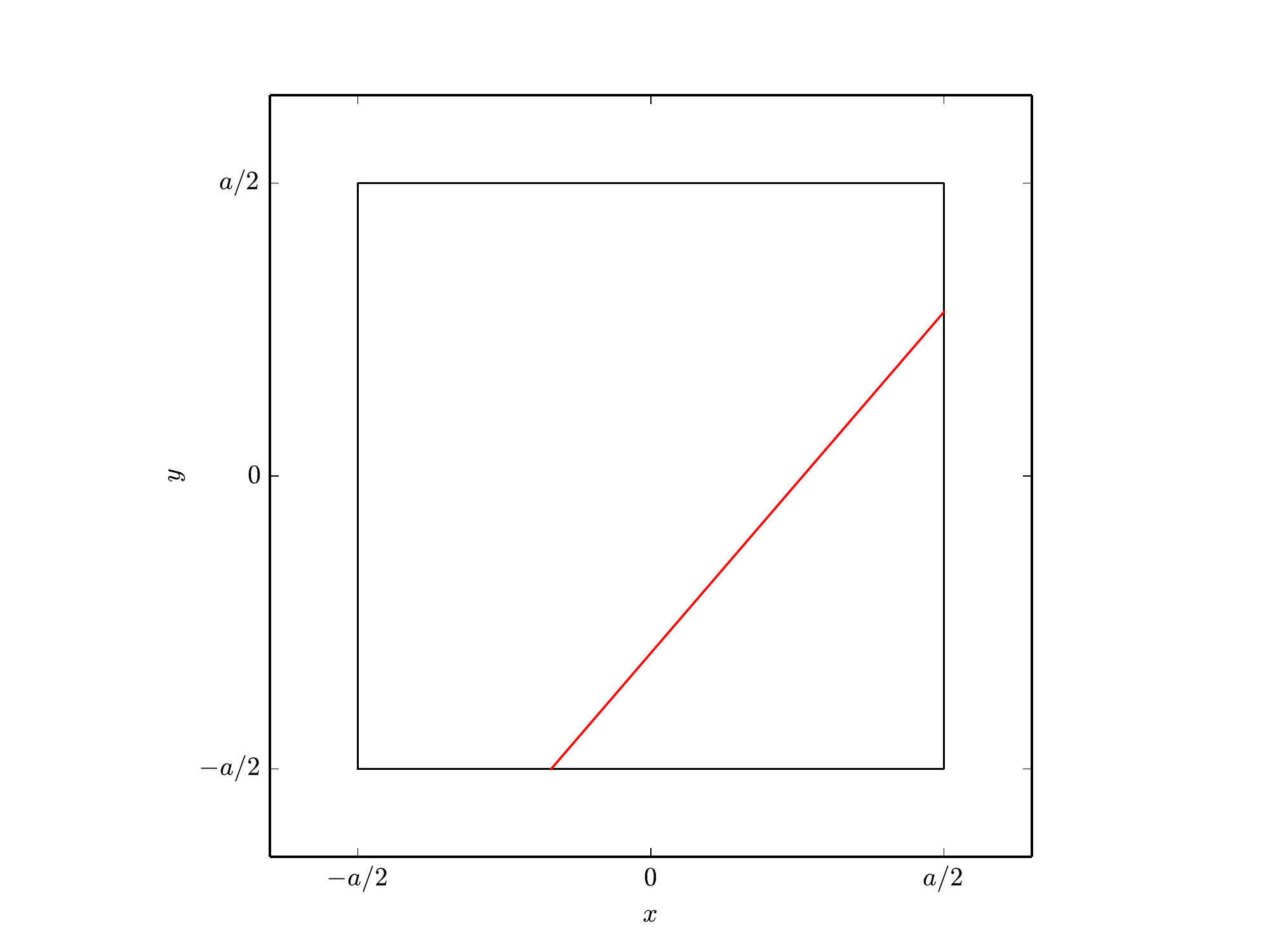}
\caption{Geometry of a square plate reinforced by a single stiffer (red color) going from $ ( -0.17a, -a/2) $ to $ (a/2,0.28a) $.
A uniform lateral loading is applied onto the plate. }
\label{fig:sketchStiffenedPlateDiagonalBeam}
\end{figure}

\begin{figure}
\includegraphics[width=\linewidth]{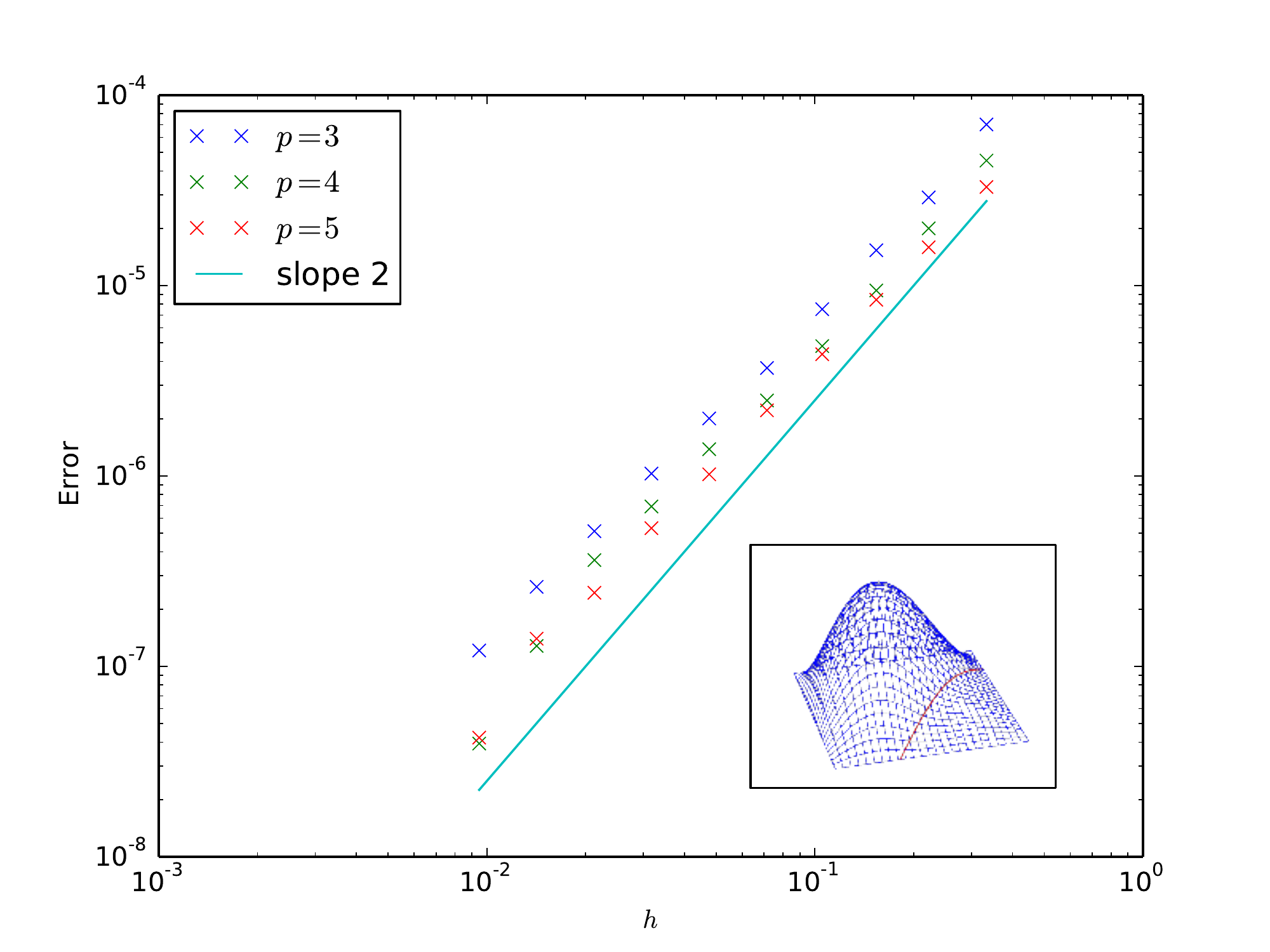}
\caption{Convergence of the error of the numerical solution 
with respect to the cell size $ h $ for
the bending problem sketched in figure \ref{fig:sketchStiffenedPlateDiagonalBeam}.
The degree of the B-splines is given by $ p = 3,4,5 $.} 
\label{fig:convergenceRectangularPlateAndBeam}
\end{figure}

In section \ref{sec:bendingRectangularPlateAlignedStiffener}, we observed that 
when the stiffener is no longer aligned with the grid lines, the rate of
convergence is reduced due to the jump in the third derivative of the solution
at the location of the beam. When faced with a stiffener crossing a simply
supported square plate
at a specific angle, cf. figure \ref{fig:sketchStiffenedPlateDiagonalBeam}, 
we expect a reduction of the order of convergence. This can also 
be observed in the convergence plot, figure \ref{fig:convergenceRectangularPlateAndBeam},
for the present choice of parameters:
\be
a = 2 \quad D = 1.234 \quad EI = 150 \quad p_0 = 1.543. 
\ee
The order of convergence in figure \ref{fig:sketchStiffenedPlateDiagonalBeam} is,
however, lower than anticipated from equation (\ref{eq:jumpCondition}), i.e. two instead of three.
The reason for this might lie in the development of a singularity at the intersections
between stiffener and boundary of $ \Omega $. A plot of the spatial
distribution of the error for a specific resolution $ h =0.02 $, cf. figure \ref{fig:sketchStiffenedPlateDiagonalBeam}, supports this possibility. 
Away from the beam the numerical solution seems to be close to the reference solution. 
Close to the beam we observe some smaller wiggles, whereas at the obtuse angles
between boundary and stiffener, the error of the solution is largest. The appearance of
singularities between stiffener and boundary 
reduces the accuracy of the method even if we had rotated the domain such that 
the stiffener would be aligned with one grid axis. 

\begin{figure}
\includegraphics[width=\linewidth]{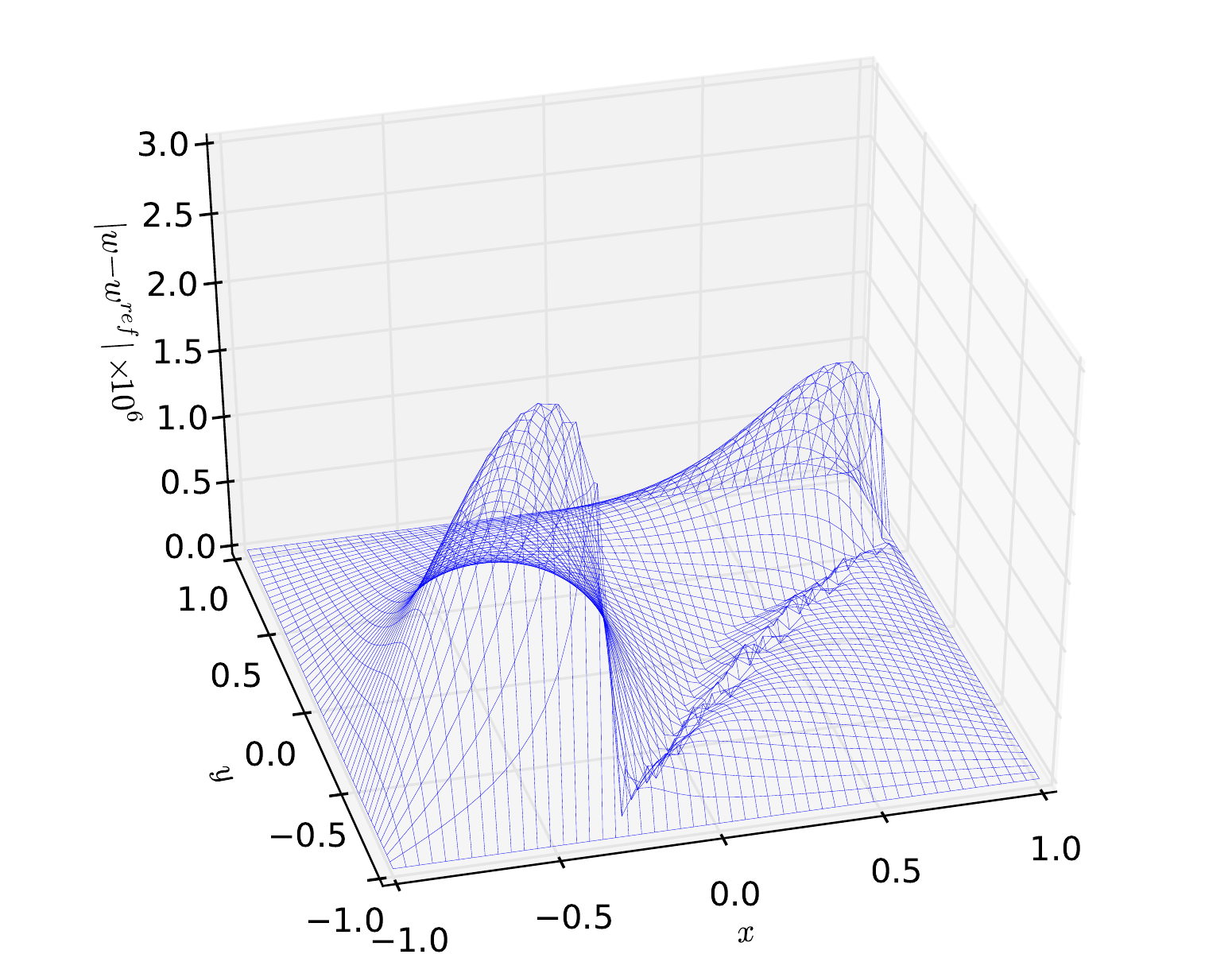}
\caption{Absolute difference between the solution $ w $ for $ p = 3 $ and $ h = 0.02 $ and the reference solution $ w^{ref}$ ($ p = 3 $ and $ h =0.005 $) for
the problem sketched in figure \ref{fig:sketchStiffenedPlateDiagonalBeam}.}
\label{fig:plateDiagonalBeamError}
\end{figure}

\subsubsection{Buckling of a square plate with non-aligned stiffener}

\begin{figure}
\includegraphics[width=\linewidth]{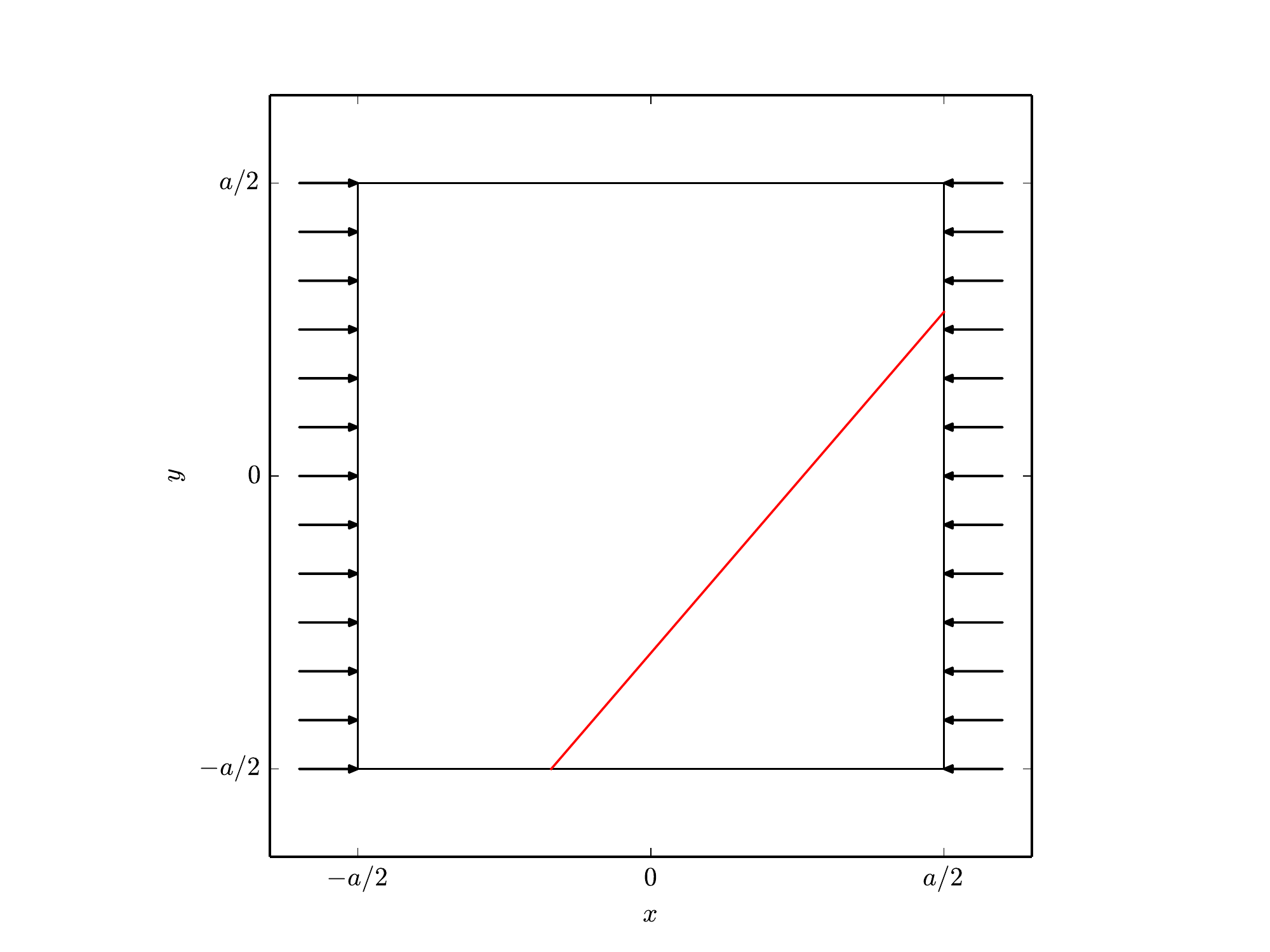}
\caption{Geometry of a square plate reinforced by a single stiffer (red color) going from $ ( -0.17a, -a/2) $ to $ (a/2,0.28a) $.
A uniform in plane lateral loading is applied at the right and left boundary. }
\label{fig:sketchStiffenedPlateDiagonalBeamBuckling}
\end{figure}

\begin{figure}
\includegraphics[width=\linewidth]{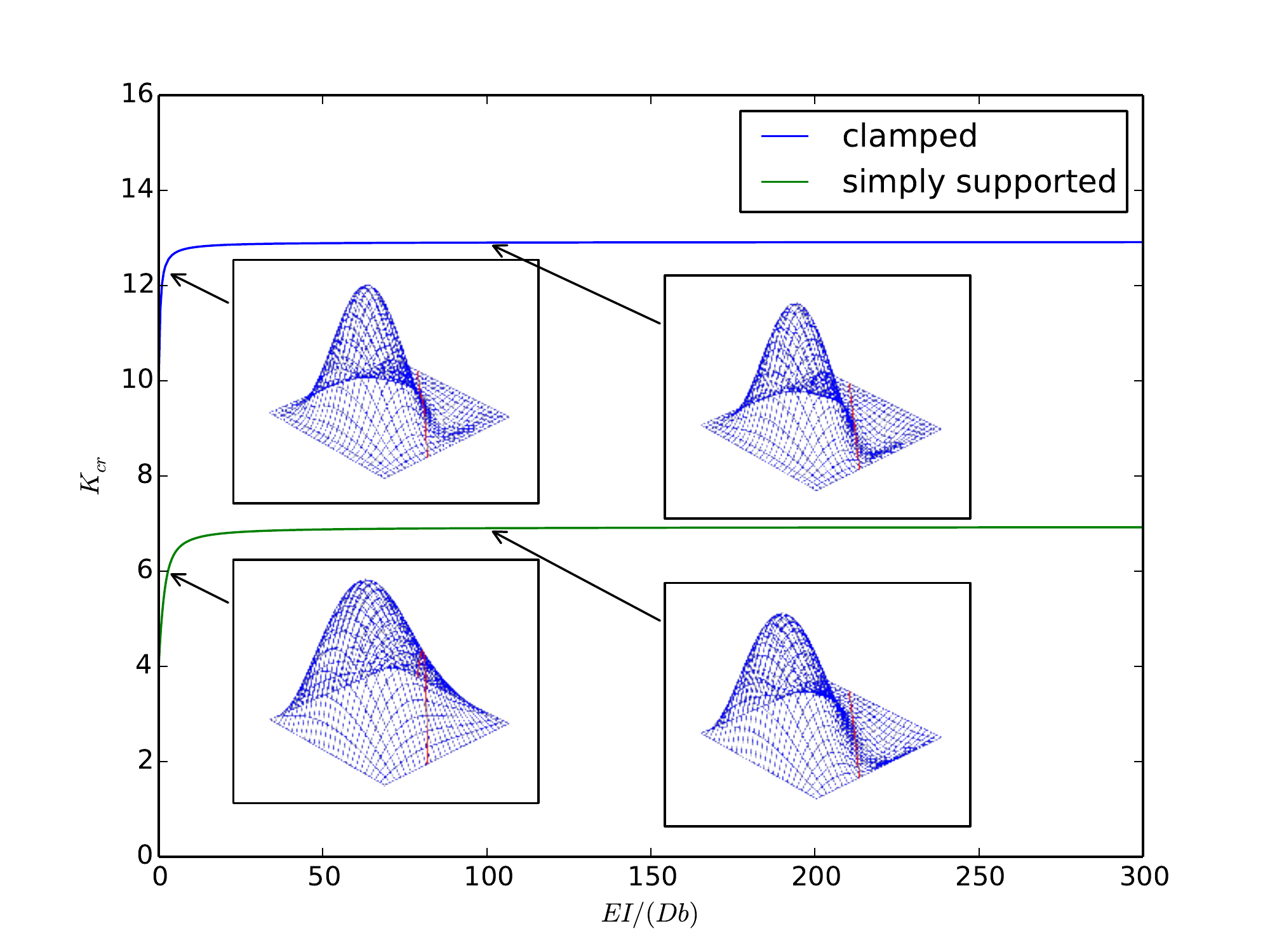}
\caption{Buckling of a square plate reinforced by a single
stiffener (red color) going from $ ( -0.17a, -a/2) $ to $ (a/2,0.28a) $ 
for simply supported and clamped boundary conditions.}
\label{fig:bucklingRectangularPlateAndBeam}
\end{figure}

Instead of a lateral loading, an in-plane loading can be applied to the left and
right boundaries of the domain, cf. figure \ref{fig:sketchStiffenedPlateDiagonalBeamBuckling}.
Assuming that no axial forcing is applied onto the stiffener, $ T_s = 0 $, we
solve (\ref{eq:masterBuckling}) for a simply supported and clamped square plate for
different values of $ EI/(Db)$. 
As can be observed from figure \ref{fig:bucklingRectangularPlateAndBeam}, the
buckling stress is a smooth function of the stiffness ratio. For large ratios,
the stiffener is almost flat, whereas for $ EI = 0 $, we obtain the reference $ K $
values for simply supported and clamped plates, i.e. $ 4 $ and $ 9.4 $, respectively. 

\subsubsection{Pre-buckling stress of a rectangular plate with central hole}
\label{sec:stressRectangularPlateHole}

\begin{figure}
\includegraphics[width=\linewidth]{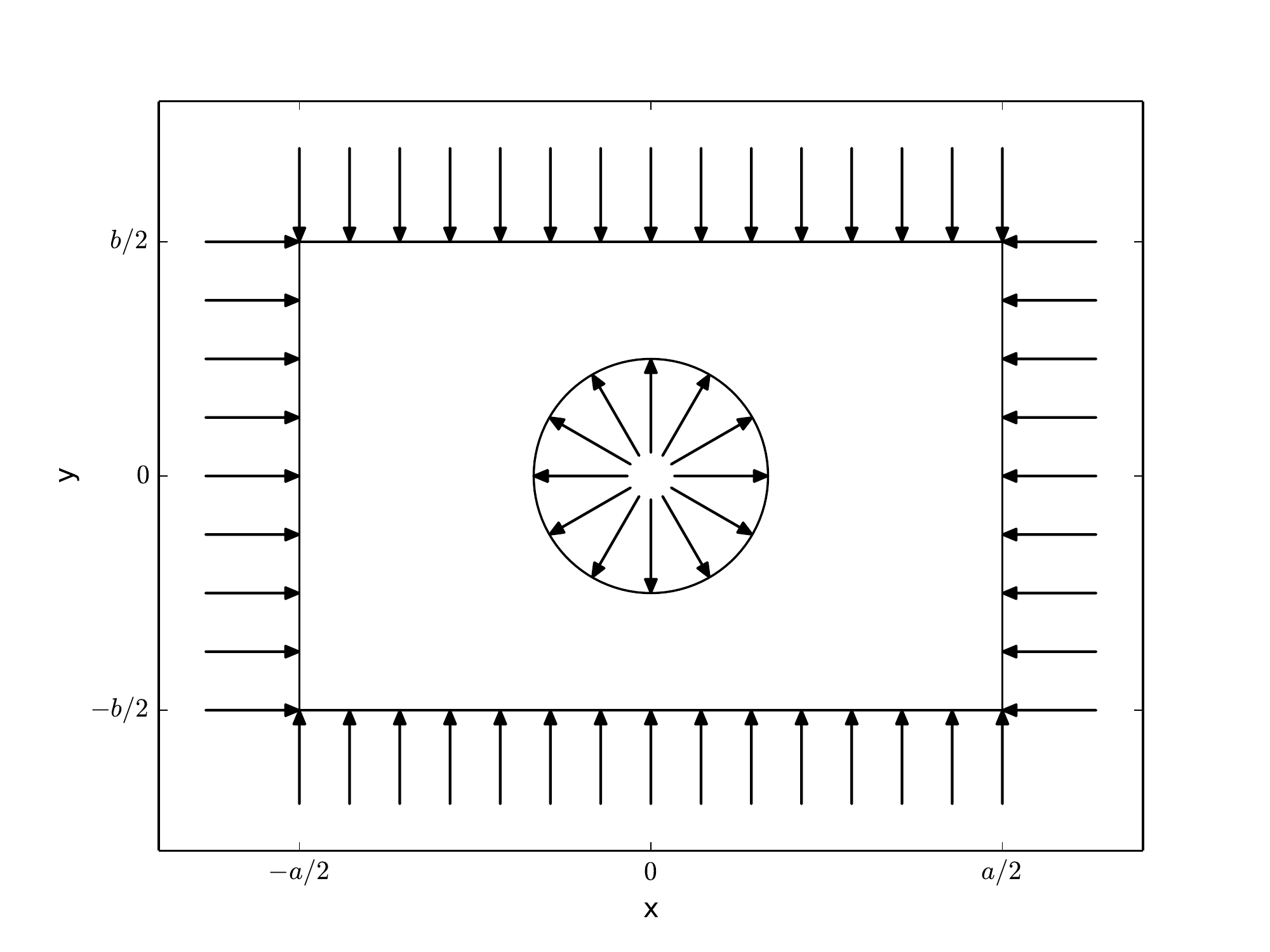}
\caption{Geometry of a rectangular plate $ [-a/2,a/2] \times [-b/2,b/2] $
with central hole with diameter $ d$. A uniform in-plane compression force is
applied in normal direction to the boundaries, such that the 
resulting stress tensor is given by $ \sigma_{xx} = \sigma_{yy} = -1 $
and $ \sigma_{xy} = 0 $.}
\label{fig:constantStressCentralHole}
\end{figure}

\begin{figure}
\includegraphics[width=\linewidth]{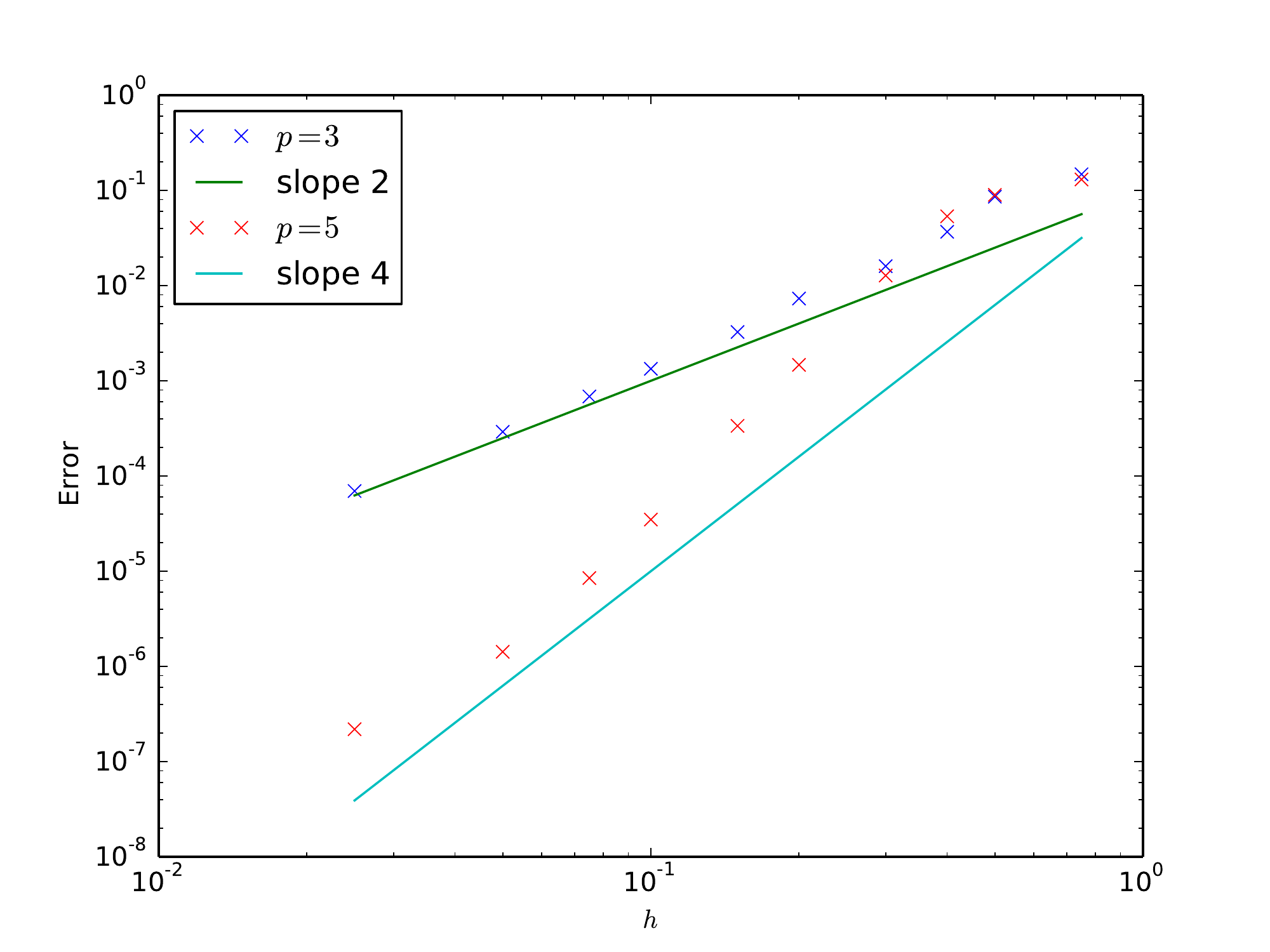}
\caption{Error convergence of the pre-buckling stress problem sketched in
figure \ref{fig:constantStressCentralHole}. Since the error is measured on
the stress components which are given by the second derivatives of Airy's stress
function, the order of convergence is reduced by two.}
\label{fig:convergenceStressCentralHole}
\end{figure}

In the following, we add a hole in the center of the rectangular plate, 
cf. figure \ref{fig:constantStressCentralHole}. When considering buckling, 
we are faced with the difficulty that no
analytic formula is known for the pre-buckling stress for
a rectangular plate with a central hole with free boundary conditions.
Therefore, the pre-buckling stress
will be computed numerically by the approach presented
in section \ref{sec:numericsStress}. In a first attempt, we shall 
verify the present method by applying it to the case when uniform
normal compression is applied both at the outer boundary and at the hole,
cf. figure \ref{fig:constantStressCentralHole}. For this configuration,
the reference solution is given by
\be
\sigma_{xx} = \sigma_{yy} = -1 \quad \sigma_{xy} = 0. 
\ee
For a plate with parameters 
\be
a = 4 \quad b = 2 \quad d = 1,
\ee
where $ d $ is the diameter of the hole, the error of the numerical solution
in function of $ h $ for splines of degree $ 3 $ and $ 5 $ is plotted in
figure \ref{fig:convergenceStressCentralHole}. Since the error is measured on
the stress components which are given by the second derivatives of Airy's stress
function, the order of convergence is reduced by approximately two. \\

\begin{figure}
\includegraphics[width=\linewidth]{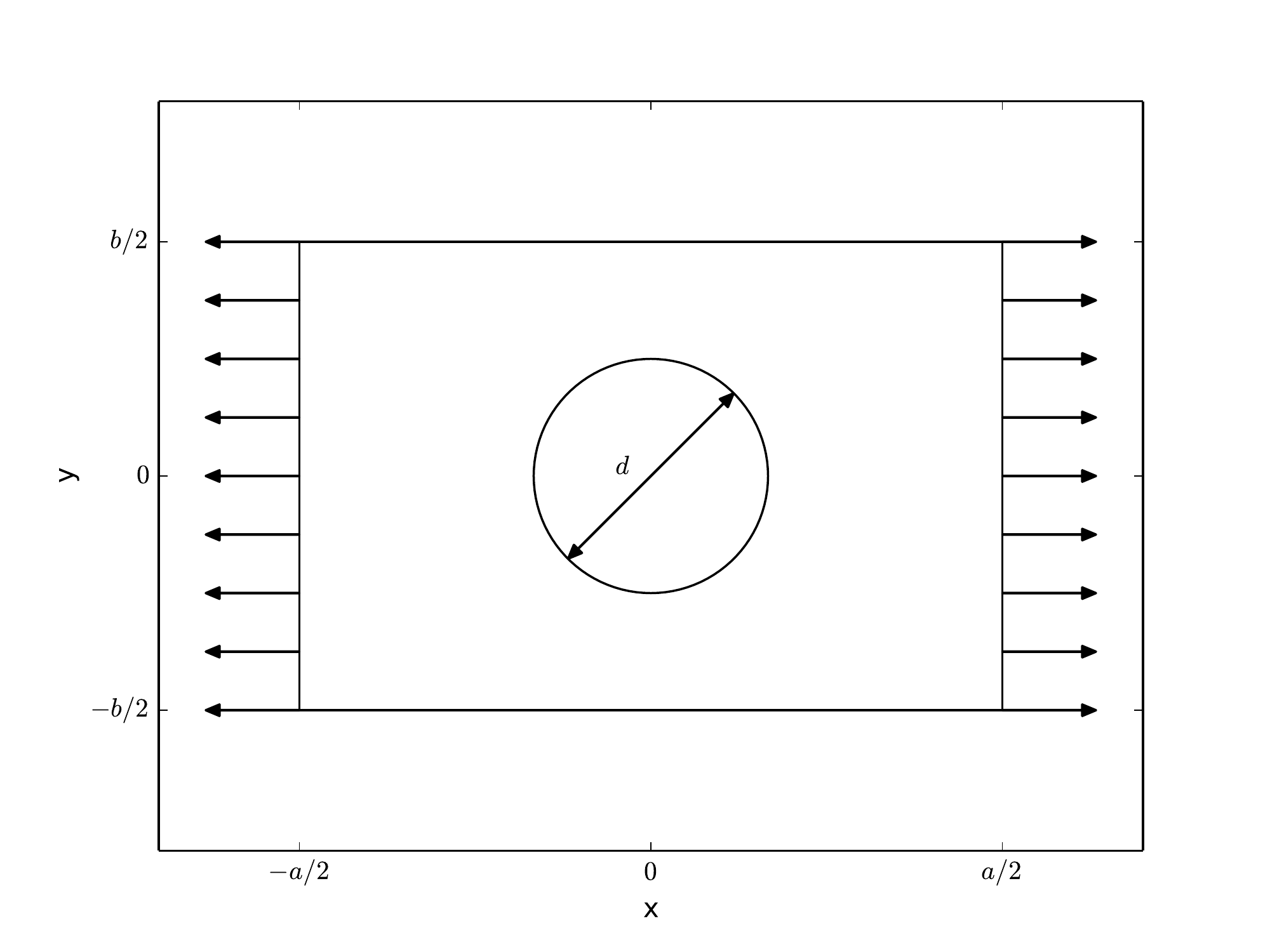}
\caption{Geometry of a rectangular plate $ [-a/2,a/2] \times [-b/2,b/2] $
with central hole with diameter $ d$. A uniform in-plane tension force is
applied on the right and left boundaries. The set-up corresponds to
the plane stress problem solved in \cite{Hayashi1989}.}
\label{fig:StressCentralHole}
\end{figure}

\begin{figure}
\includegraphics[width=\linewidth]{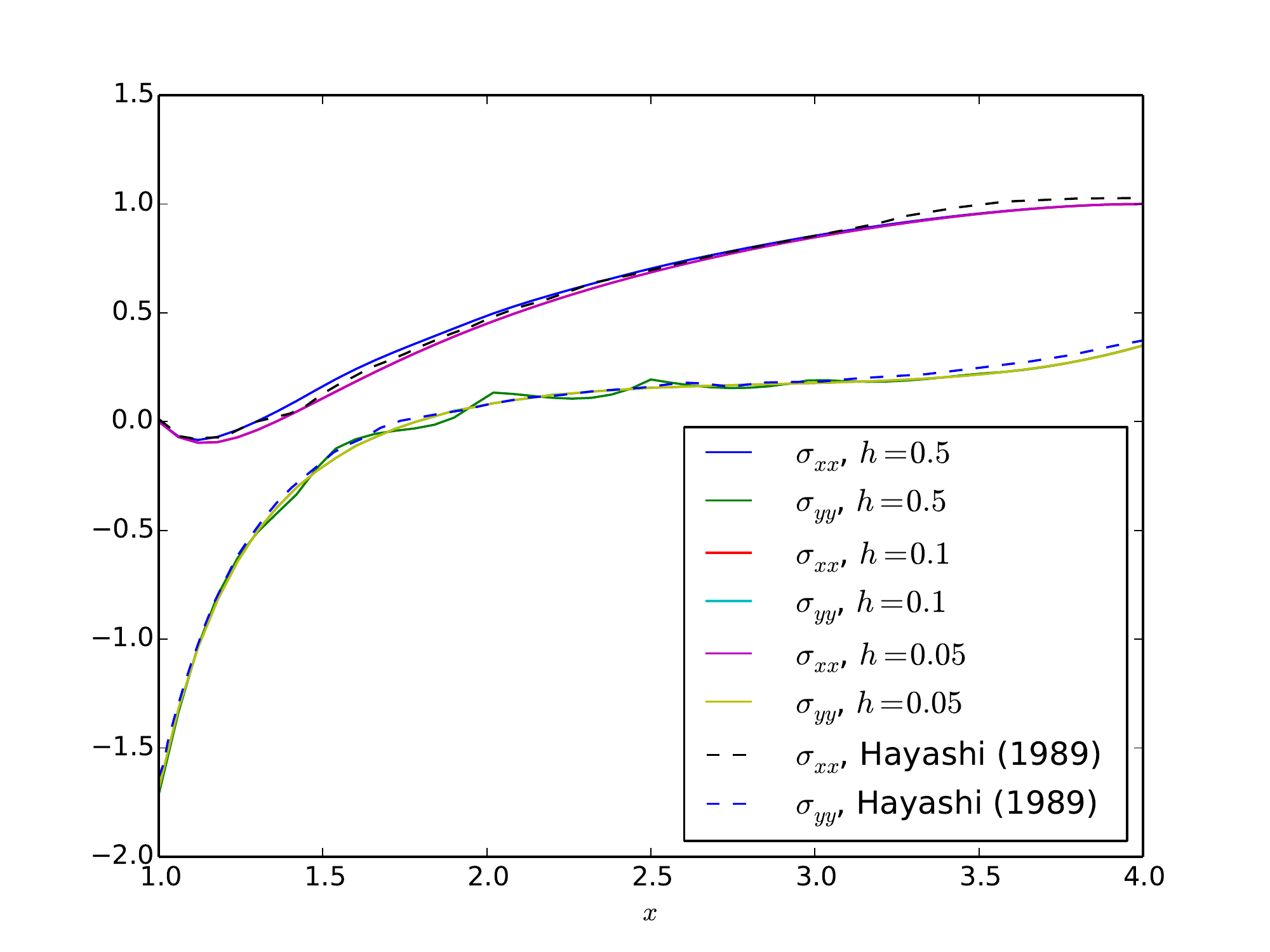}
\caption{Stress along the positive part of the $ x $ axis for the
problem sketched in figure \ref{fig:StressCentralHole}. The numerical
results for $ p = 3 $ and different choices of $ h $ are compared
to the results obtained by \cite{Hayashi1989}. The data 
for Hayashi has
been digitized form figure 1 in \cite{Hayashi1989}.
The lines for $ h = 0.1 $ and $ h = 0.05 $ are indistinguishable on
a plotting scale.}
\label{fig:ComparisonHayashi}
\end{figure}

For the case of a uniform tension force on the right and left boundary 
and a free boundary at the upper and lower edges and at the hole, cf.
figure \ref{fig:StressCentralHole}, a reference
solution is given in figure 1 in \cite{Hayashi1989}. 
Applying the present method to
the case treated in \cite{Hayashi1989}, 
we can compare the values of $ \sigma_{xx} $ and $ \sigma_{yy} $
along the $ x $-axis found by the present method with the 
values obtained in \cite{Hayashi1989}, cf. 
figure \ref{fig:ComparisonHayashi}. The present
numerical solution has been obtained by means of splines of degree 3. 
For the relatively coarse resolution of $ h = 0.5 $, i.e six cells,
along the $ x $-axis in figure \ref{fig:ComparisonHayashi}, we observe
kinks in the solution due to the fact that 
the third derivative of the basis functions for $ p = 3 $
is in general discontinuous at the cell boundaries. These kinks are
no longer visible for the finer resolutions $ h = 0.1, 0.05 $,
when the solution is indistinguishable on a plotting scale. 
When choosing $ p = 5 $, the solution is smooth even for coarse
resolutions (figure not shown). In general the lines follow
the solution by \cite{Hayashi1989} closely. The discrepancy
might be associated to the modest plotting quality of figure 1
in \cite{Hayashi1989}, which is passed further to
the digitized data. As the stress components
enter the buckling energy (\ref{eq:Energy3}), it is advantageous to
use B-splines of degree $ p + 2 $ for Airy's stress function, when
B-splines of degree $ p $ are used for the vertical displacement $ w $. 

\subsubsection{Buckling of a square plate with central hole} \label{sec:bucklingCentralHole}

\begin{figure}
\includegraphics[width=\linewidth]{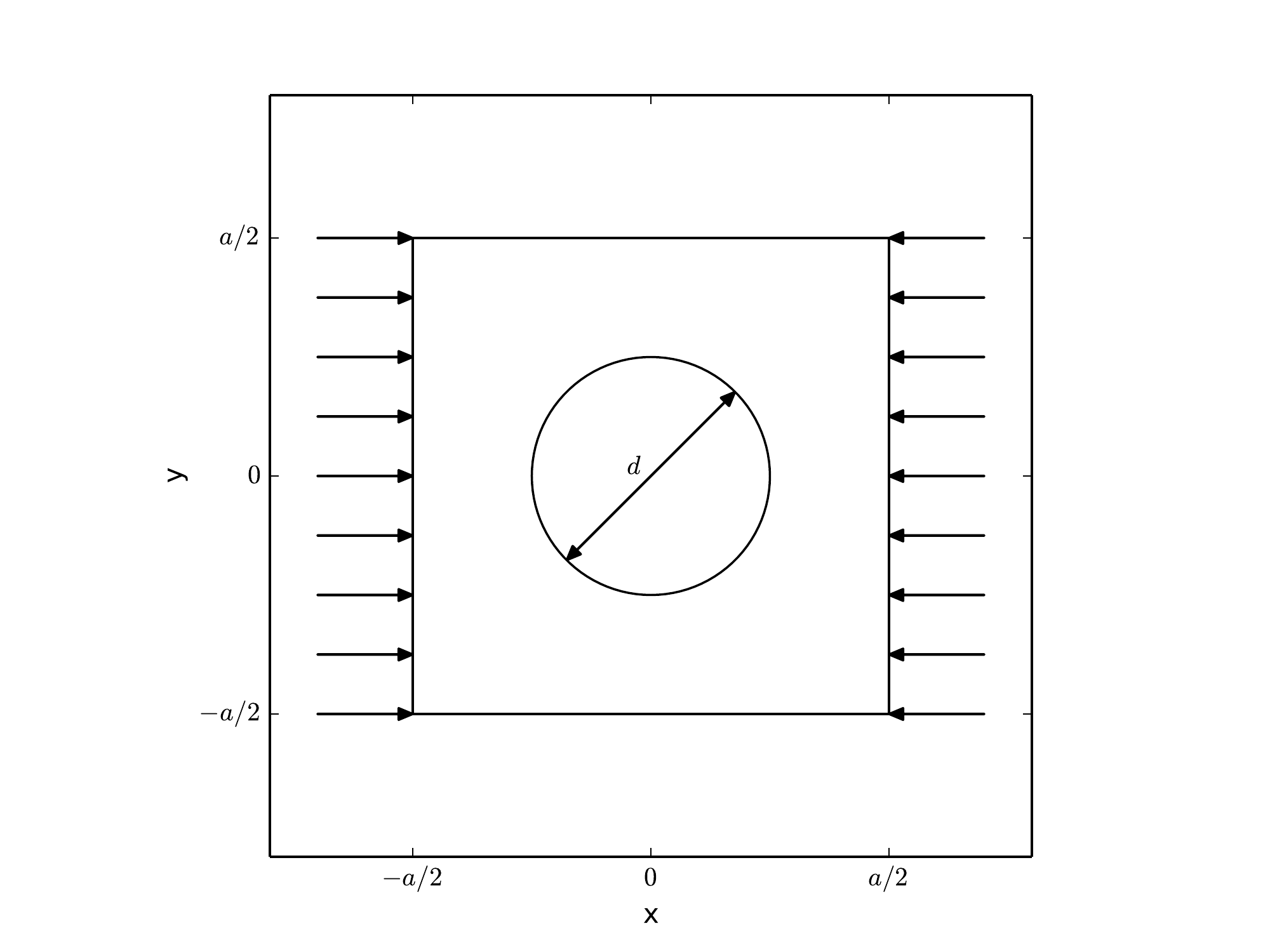}
\caption{Geometry of a square plate with central hole.
A uniform in-plane loading is applied onto the left and right
boundaries. }
\label{fig:sketchBucklingSquarePlateWithHole}
\end{figure}

\begin{figure}
\includegraphics[width=\linewidth]{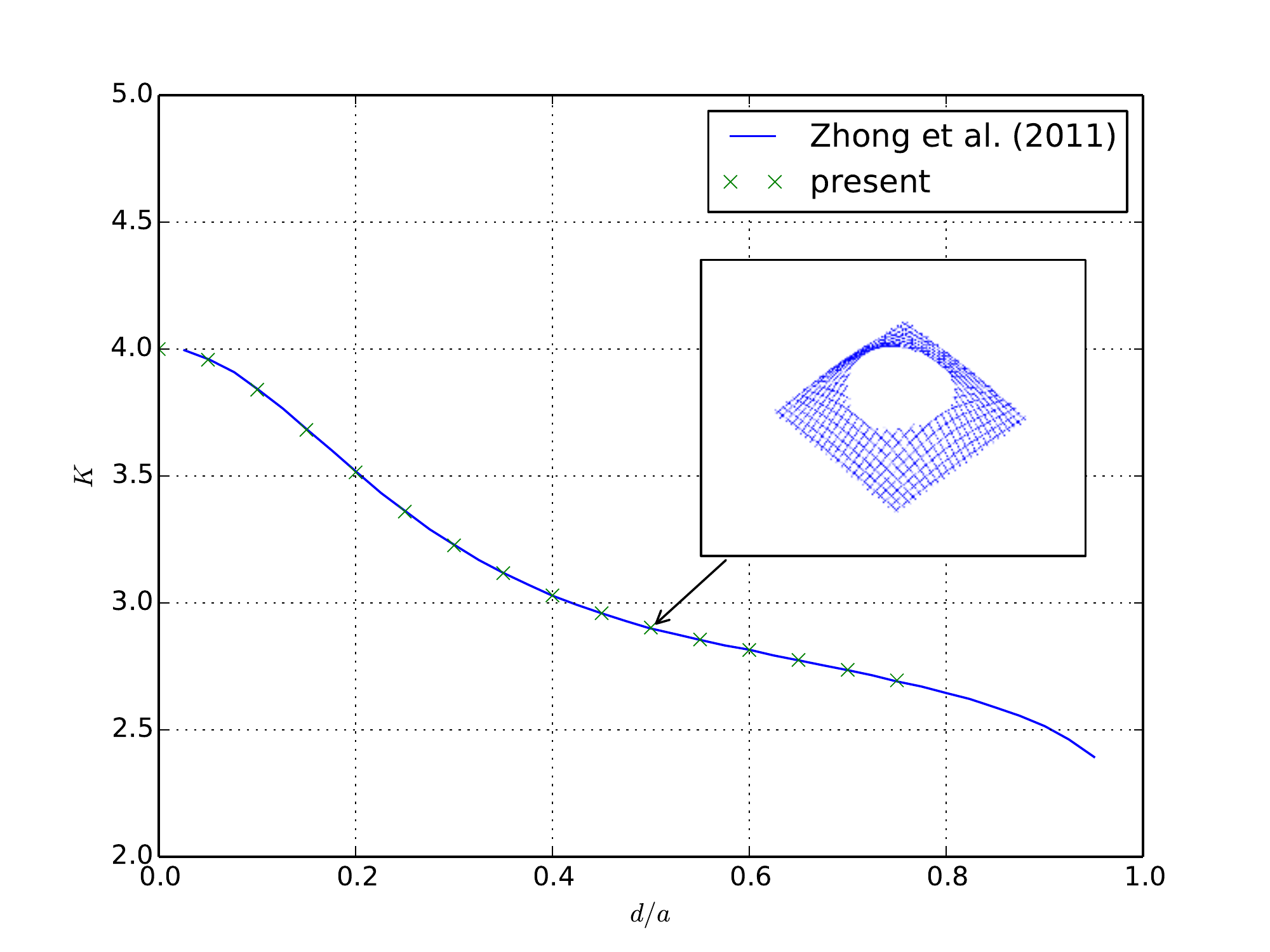}
\caption{Buckling stress for a square plate with central hole in function
of the ratio between hole diameter $ d $ and plate dimension $ a$. 
The data for Zhong {\em et al.} has been digitized from figure 13 in
\cite{ZhongPanYu2011}. The embedded figure shows the eigenfunction for the case $ d/a = 0.5$.}
\label{fig:bucklingPlateCentralHole}
\end{figure}

This case has been solved in several works,
for example \cite{KawaiOhtsubo1968,ZhongPanYu2011,DjelosevicTepicTanackovKostelac2013}.
A square plate with a free hole in its center is uniformly loaded in-plane 
at its left and right boundary,
cf. figure \ref{fig:sketchBucklingSquarePlateWithHole}. 
The pre-buckling stress field has been computed by 
means of the method in section \ref{sec:numericsStress}. 
This stress field is then used to solve the buckling problem
(\ref{eq:masterBuckling}) by the present method with splines of
degree $ p = 3 $, a resolution of $ h = 0.025,0.05,0.1 $, a side
length of $ a = 2$,
and for a Poisson module of $ \nu = 0.3 $. 
The $ K $ values of the buckling stress computed by the present method
are compared to the results obtained by 
Zhong {\em et al.}\cite{ZhongPanYu2011}
in figure \ref{fig:bucklingPlateCentralHole}.
The results by Zhong {\em et al.}
were digitized from figure 13 in \cite{ZhongPanYu2011}. The present results
lie neatly on top of the graph by the data obtained by \cite{ZhongPanYu2011}.
In addition, we checked that, apart from the smallest hole for $ h = 0.1$,
the present results are indistinguishable on a plotting scale
for $ h = 0.025,0.05,0.1 $.  

\subsection{Polygonal plate} \label{sec:polygonal}

In this section we turn to our last test, 
where we employ the present solver on a more general
domain. We consider two cases, a plate with holes bounded by a convex polygon
and a plate with holes bounded by a simple, non convex polygon, cf. figures
\ref{fig:sketchBendingConvexPolygon} and \ref{fig:sketchBendingNonConvexPolygon}.
We shall first consider bending, cf. section \ref{sec:bendingPolygon}, before
turning our attention to buckling, cf. section \ref{sec:bucklingPolygon}. 
As before, a definition for the weight function
in equation (\ref{eq:webSpline}) 
needs to be given for the present domains. 
For the geometry bounded by a convex polygon, a weight function
can be written as:
\be
\omega(x,y) = \left( \prod \limits_{i=0}^4 \omega_i(x,y) \right)^2,
\label{eq:weightConvexPolygon}
\ee
where $ \omega_i $ is defined in equation (\ref{eq:edgeweight}).
The right hand side in (\ref{eq:weightConvexPolygon}) is squared in
order to model clamped boundary conditions. For the external boundary 
of the non-convex polygon, on the other hand, a simple
product of the edge weights $ \omega_i $ as in (\ref{eq:weightConvexPolygon})
cannot be used. For the inward facing corner at $ \mbfv_0$, set operations
as for the R-function method \cite{RvachevSheiko1995} need to be used
in order to define a weight function. The resulting weight
function is in this case given by:
\be
\omega(x,y) = \left( \omega_{\corner}\prod \limits_{i=1}^3 \omega_i(x,y) \right)^2
\label{eq:weightNonConvexPolygon},
\ee
where
\be
\omega_{\corner} = \omega_0 + \omega_4 + \sqrt{ \omega_0^2 + \omega_4^2 } .
\ee
For the computation of the stress field for the buckling problem,
we need in addition 
a simple weight for the holes in the domain. The weight
\be
\omega_{\hole} = \left( x - x_0 \right)^2 + \left(y-y_0 \right)^2 - R^2
\ee 
is used for a circular hole of radius $ R $ centered at $ (x_0,y_0)$.
The coordinates of the vertices of the bounding polygons
and the position and radius of the holes
can be found in \ref{sec:coordinatesPolygons}.

\subsubsection{Bending of a polygonal plate with holes} \label{sec:bendingPolygon}

\begin{figure}
\includegraphics[width=\linewidth]{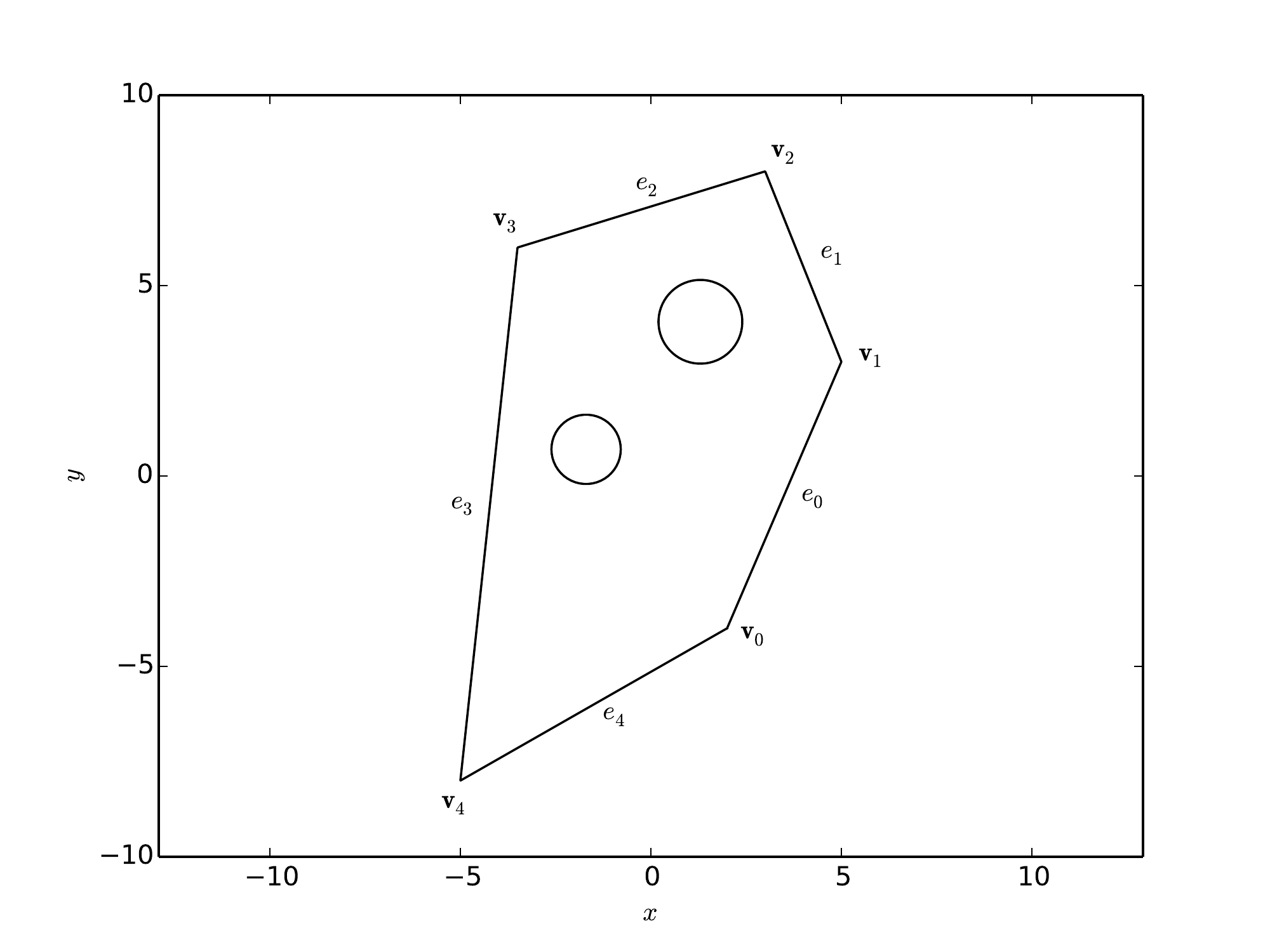}
\caption{Geometry of a plate with holes bounded by a convex polygon.
A uniform lateral loading is applied onto the plate.}
\label{fig:sketchBendingConvexPolygon}
\end{figure}

\begin{figure}
\includegraphics[width=\linewidth]{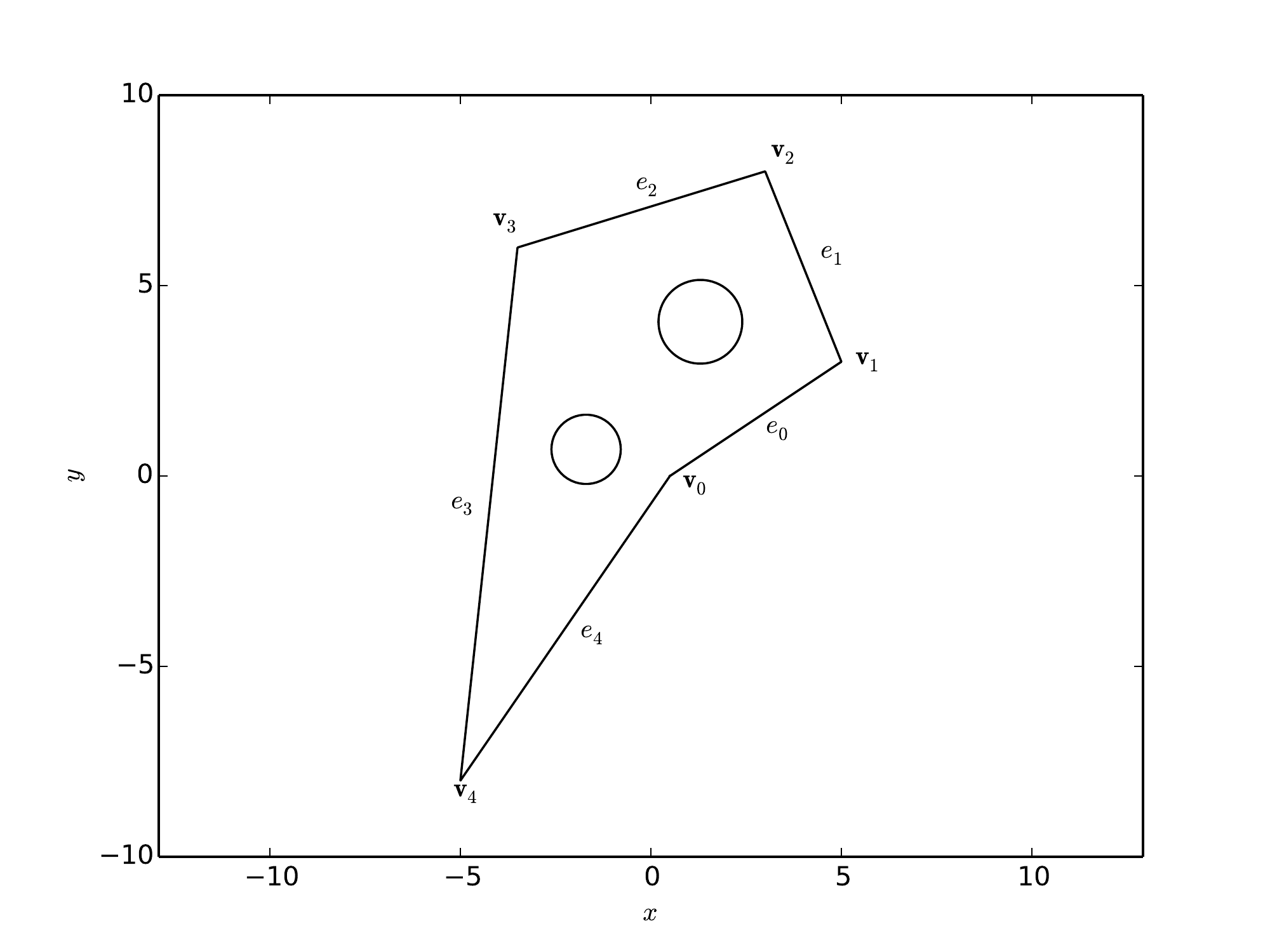}
\caption{Geometry of a plate with holes bounded by a non convex polygon.
A uniform lateral loading is applied onto the plate.}
\label{fig:sketchBendingNonConvexPolygon}
\end{figure}

\begin{figure}
\includegraphics[width=\linewidth]{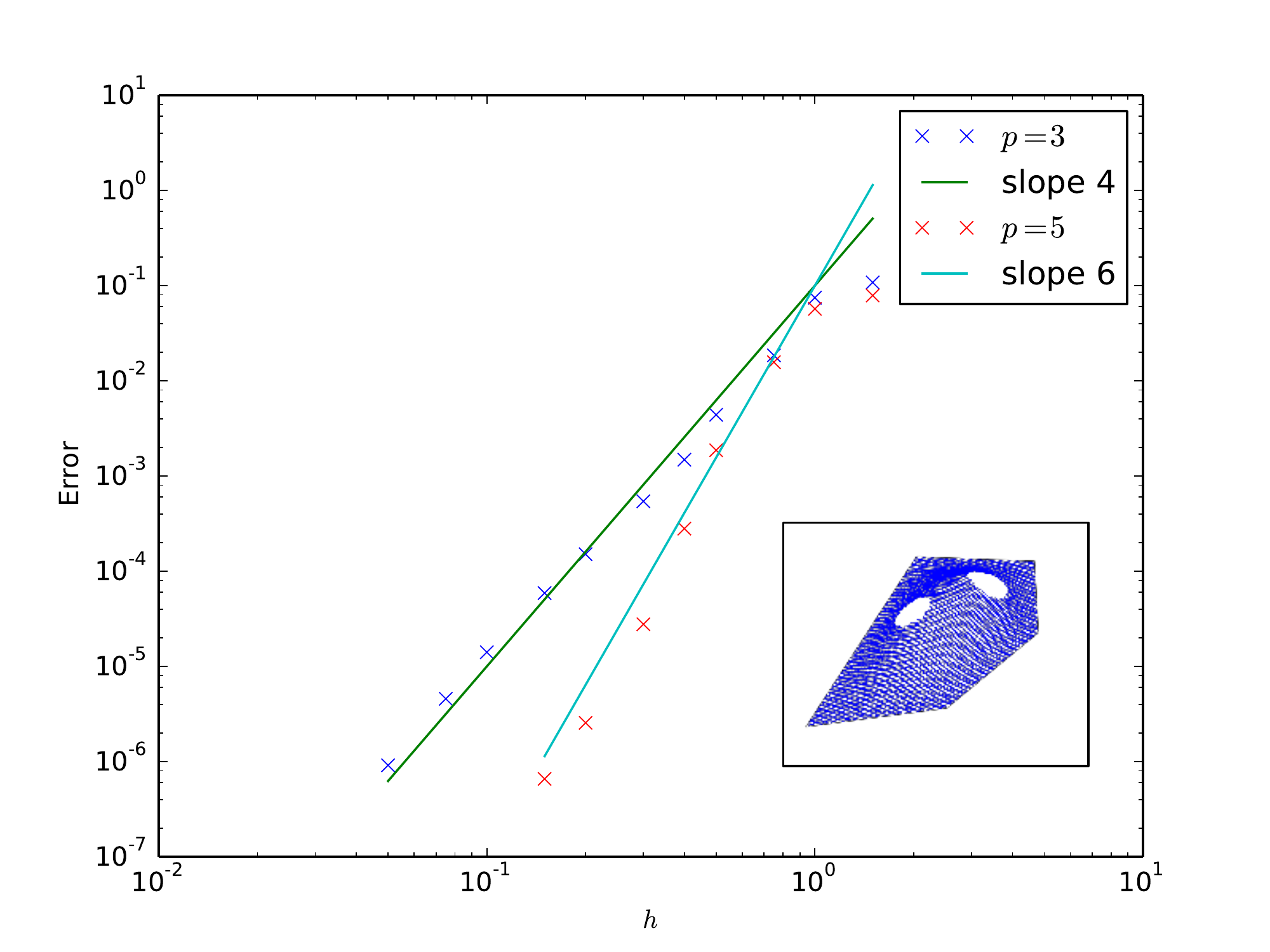}
\caption{Error convergence of the bending problem sketched in
figure \ref{fig:sketchBendingConvexPolygon}.}
\label{fig:convergenceBendingConvexPolygon}
\end{figure}

\begin{figure}
\includegraphics[width=\linewidth]{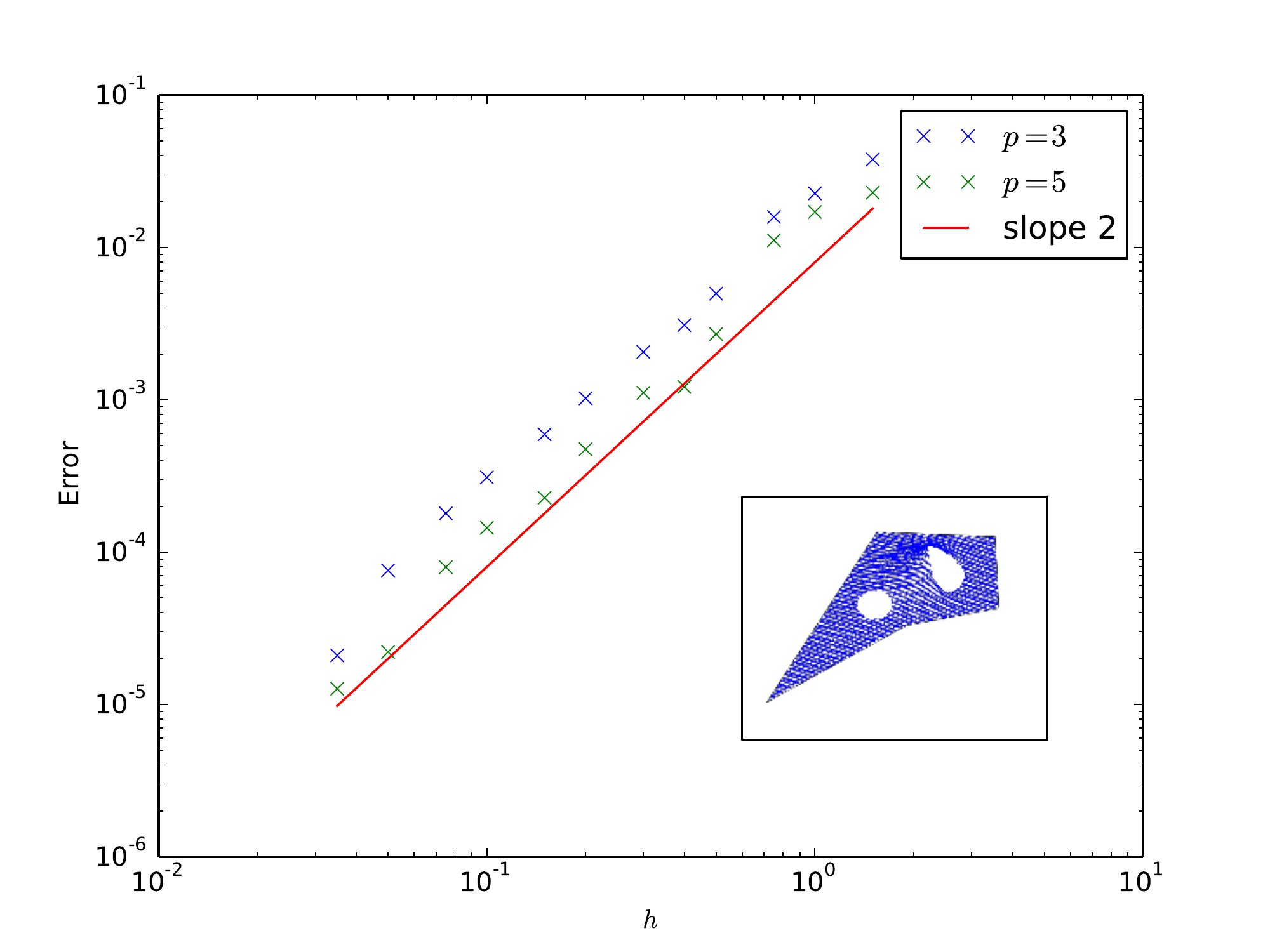}
\caption{Error convergence of the bending problem sketched in
figure \ref{fig:sketchBendingNonConvexPolygon}.}
\label{fig:convergenceBendingNonConvexPolygon}
\end{figure}

In the present section, a constant lateral loading $ p_0 $ is applied onto the plates
defined in figures \ref{fig:sketchBendingConvexPolygon} and 
\ref{fig:sketchBendingNonConvexPolygon}.  The parameters of the simulation
are given by:
\be
p_0 = 1.234, \quad D = 1, \quad \nu = 0.3. 
\ee
By using a numerical solution on a fine grid as reference solution, a
numerical error can be computed. 
As can be seen from figure \ref{fig:convergenceBendingConvexPolygon}, the
theoretical convergence rates for the solver are approximately recovered
for the domain with convex outer boundary. However,
when applying the method to the domain displaying an inward facing
corner at $ \mbfv_0$,
cf. figure \ref{fig:sketchBendingNonConvexPolygon}, the rate of convergence
drops to a value of approximately two, cf. figure \ref{fig:convergenceBendingNonConvexPolygon}.
The reason for this is the development of a singularity at the external corner 
reducing the convergence rate of the method, cf. \cite{BlumDobrowolski1982}. 
 
\subsubsection{Buckling of a polygonal plate with holes} \label{sec:bucklingPolygon}

\begin{figure}
\includegraphics[width=\linewidth]{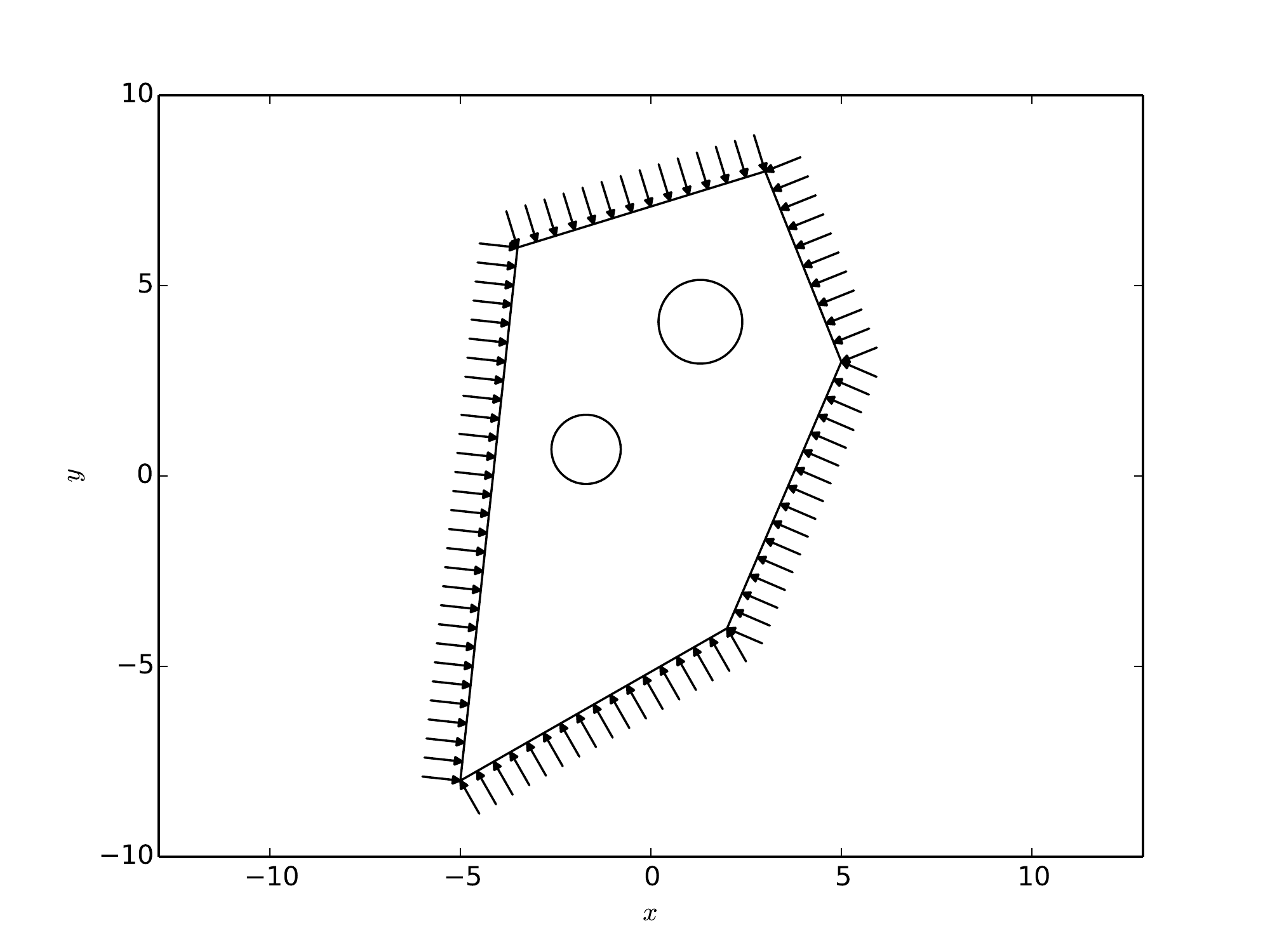}
\caption{Geometry of a plate with holes bounded by a convex polygon.
A uniform in-plane loading in-normal direction applied onto the edges of the outer boundary.
The holes are free.}
\label{fig:sketchBucklingConvexPolygon}
\end{figure}

\begin{figure}
\includegraphics[width=\linewidth]{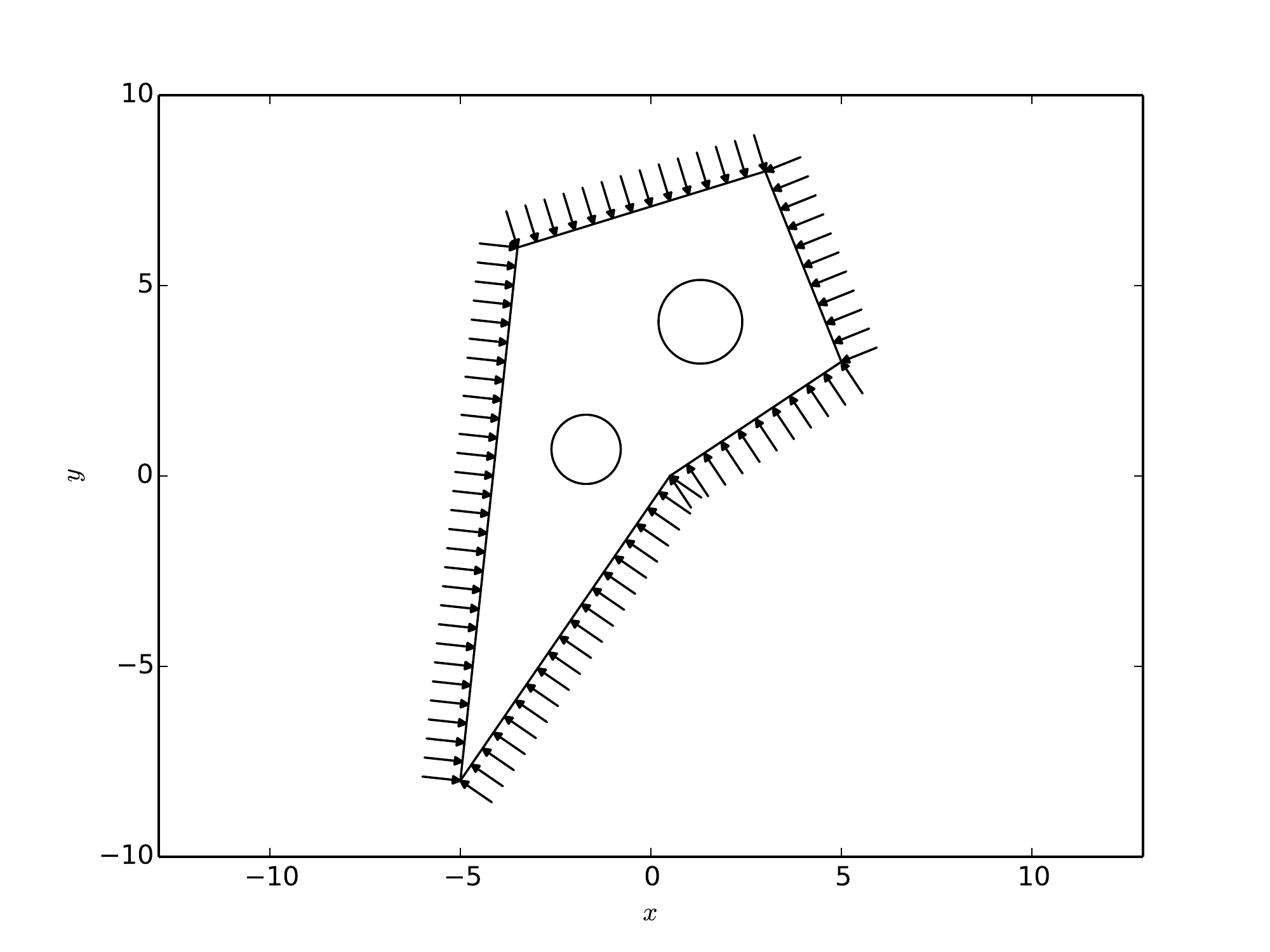}
\caption{Geometry of a plate with holes bounded by a non convex polygon.
A uniform in-plane loading in-normal direction applied onto the edges of the outer boundary.
The holes are free.}
\label{fig:sketchBucklingNonConvexPolygon}
\end{figure}

\begin{figure}
\includegraphics[width=\linewidth]{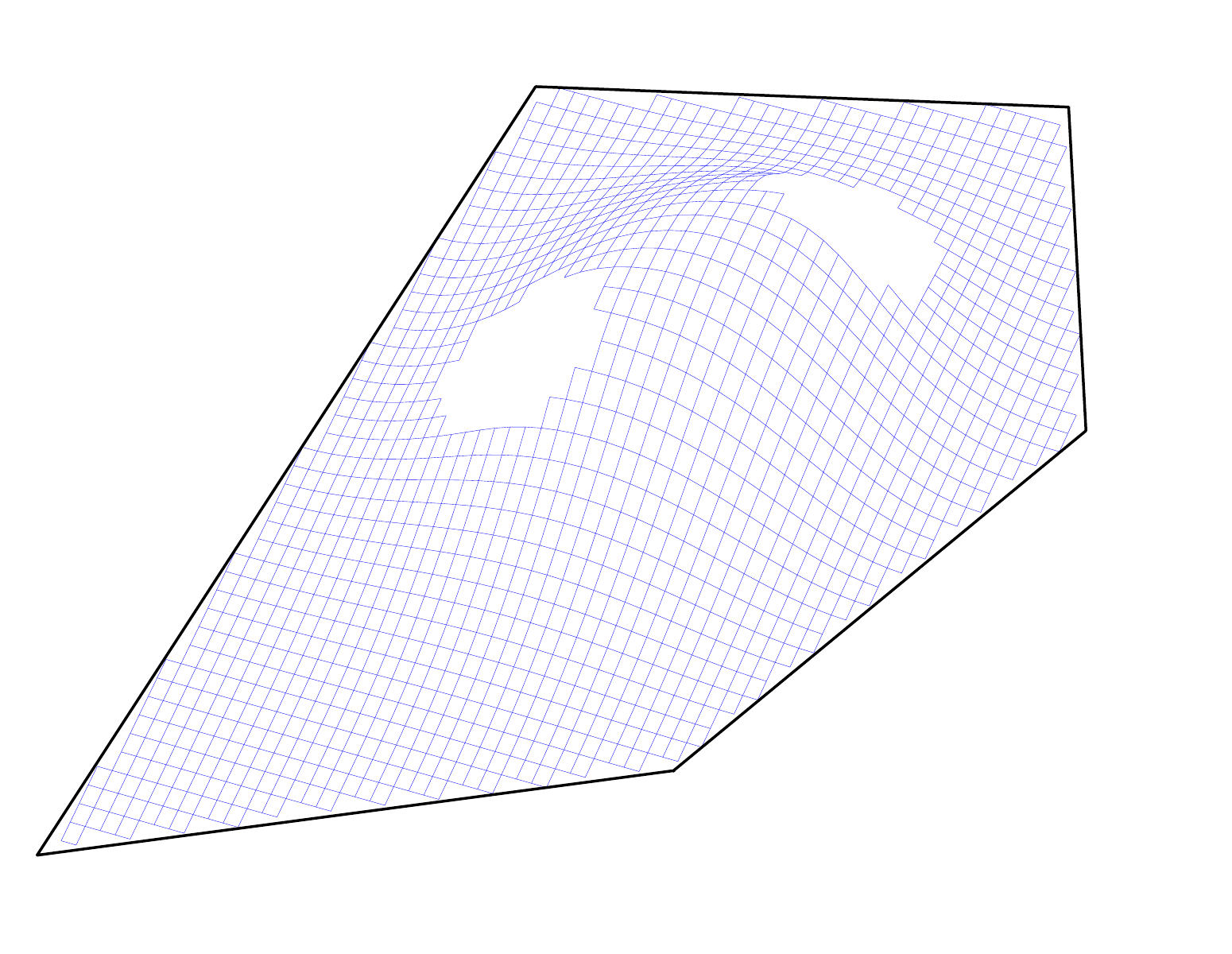}
\caption{Most unstable buckling eigenmode for the plate with holes bounded by
a convex polygon, cf. figure \ref{fig:sketchBucklingConvexPolygon}.}
\label{fig:surfaceBucklingConvexPolygon}
\end{figure}

\begin{figure}
\includegraphics[width=\linewidth]{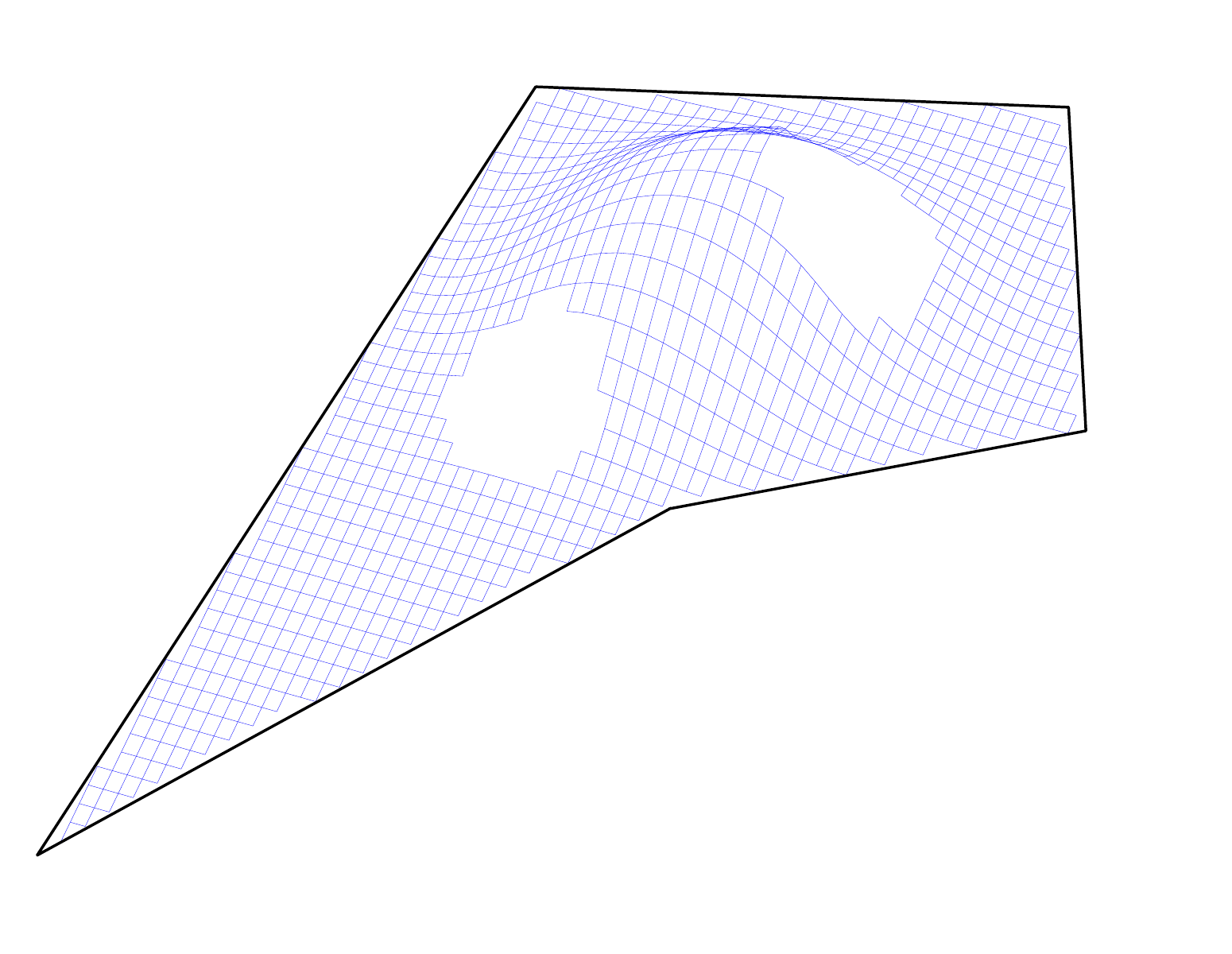}
\caption{Most unstable buckling eigenmode for the plate with holes bounded by
a non convex polygon, cf. figure \ref{fig:sketchBucklingNonConvexPolygon}.}
\label{fig:surfaceBucklingNonConvexPolygon}
\end{figure}

A simple buckling problem for the present plates 
with holes bounded by a polygon, 
can be defined by applying a normal compression onto the exterior boundary, cf.
figures \ref{fig:sketchBucklingConvexPolygon} and \ref{fig:sketchBucklingNonConvexPolygon}. 
In table \ref{tab:bucklingPolygon}, the buckling stresses $ \lambda_c $ 
are recorded for the two geometries
for two resolutions $ h$ using splines of degree $ 3$. 
As can be observed, 
for the convex polygon, four digits of
the buckling stress have been obtained, whereas for the
non-convex polygon, three digits have been reached. 
The eigenfunctions corresponding to the buckling stress are plotted
in figures \ref{fig:surfaceBucklingConvexPolygon} and
\ref{fig:surfaceBucklingNonConvexPolygon}. 

\begin{table}
\centering
\caption{Buckling of a polygonal plate with holes.}
\begin{tabular}{l|l}
$h$ & $\lambda_c/D $ \\ \hline 
\multicolumn{2}{c}{Convex polygon}\\
\hline
0.3 & 0.578828  \\
0.15 & 0.578869 \\ 
\hline 
\multicolumn{2}{c}{Non convex polygon}\\
\hline
0.3 & 0.851154 \\
0.15 & 0.851590
\end{tabular}
\label{tab:bucklingPolygon}
\end{table}

\section{Conclusions} \label{sec:conclusions}

In the present treatise, 
the weighted extended B-spline method by \cite{Hoellig2003} is applied
to bending and buckling problems of Kirchhoff plates of various shapes
for lateral and in-plane loading. The plate is allowed to
be supported by a stiffener. A range of benchmark tests is
applied to the present solver in order to document various aspects
affecting the accuracy of the present method. In particular,
we document how the jump conditions at the stiffener location
and singularities at the intersection between stiffener and boundary
of the domain or at inward facing corners reduce the accuracy of the
method. However, for smooth solutions, the present method displays
high order accuracy, making it a formidable choice for eigenvalue
problems, such as plate buckling problems, due to the relatively small
size of the stiffness matrix. 

As such, techniques have been developed for the treatment of
singularities and discontinuities 
in the framework of other methods \cite{SchultzLeeBoyd1989,BlumDobrowolski1982,LiLu2000,Boyd2001,BelytschkoGracieVentura2009}. 
However, it remains to future research how exactly these can be applied to
the weighted extended B-spline method
in order to recover the higher accuracy. As is pointed out in the
present work and also in \cite{LiApplegarthBullBettessBondThompson1997}, 
singularities have a bigger impact on the
accuracy than the reduced continuity at
the stiffener location and should therefore be addressed first. 

Concerning the computational efficiency of the present method, the
cost of the assembly of the matrices depends heavily on
the complexity of the domain. In general, the cost scales approximately
linearly with the number
of cells, i.e. $ \propto h^{-2} $. The sparse linear solver
by Tim Davis \cite{DavisWebsite} leads to a scaling of approximately
$ \propto h^{-4} $ for the solution of the bending problem (\ref{eq:discrete1}).
Concerning the eigenvalue solver, the sparsity is not taken
into account by using a Lapack solver \cite{Lapack} for full
matrices, i.e. the computational cost scales as $ \propto h^{-6} $. 
The application of an iterative eigenvalue solver 
would have reduced this cost significantly. 

A difficulty observed in the present work
for the weighted extended B-spline method is the treatment of small holes as in 
section \ref{sec:bucklingCentralHole}. As the smallest features of the
geometry dictate the resolution $ h$, small cells
are used also in regions of the domain where a coarser resolution would have
been enough. In section 4.5 in \cite{Hoellig2003}, ideas for the formulation
of a weighted extended B-spline method using hierarchical bases are presented.
Such a development allowing to use different resolutions for parts of the
domain would be a very welcome feature for the method. 

The embedded treatment of the boundary of the domain is an interesting
feature in order to facilitate the integration between CAD and structural
solver. However, it is also Achilles' heel of the method, as the cost
of the assembly of the system is dominating. We remark that the present
code might not have been fully optimized with respect to this point and
that a better balancing of the approximation error would give some improvement. 
On the other hand, a collocation formulation of the
weighted extended B-spline method as in \cite{ApprichHoelligHoernerReif2015}
would render this issue obsolete. 

The aim of the present treatise is to give an account
of the weighted extended B-spline method when applied to
plate bending and buckling problems relevant for the maritime industry.
We highlight a number of
issues which are important when
considering bending and buckling of plates. Although the accuracy of
the method is affected by these issues
(which holds true also for other methods),
it nevertheless is a remarkably accurate method applicable to a 
wide range of problems. 

\section{Acknowledgments}

The author would like to thank Dr Eivind Steen and Dr Lars Brubak
for the many interesting, helpful and motivating meetings at DNV-GL. 
Professor Jostein Hellesland
and Professor Brain Hayman are thanked for interesting discussions. 
Professor Michael Floater is heartily thanked for his patience
checking the extension algorithm. 

\begin{appendix}

\section{Extension algorithm for weighted extended B-Splines} \label{sec:extension}

\begin{figure}
\includegraphics[width=\linewidth]{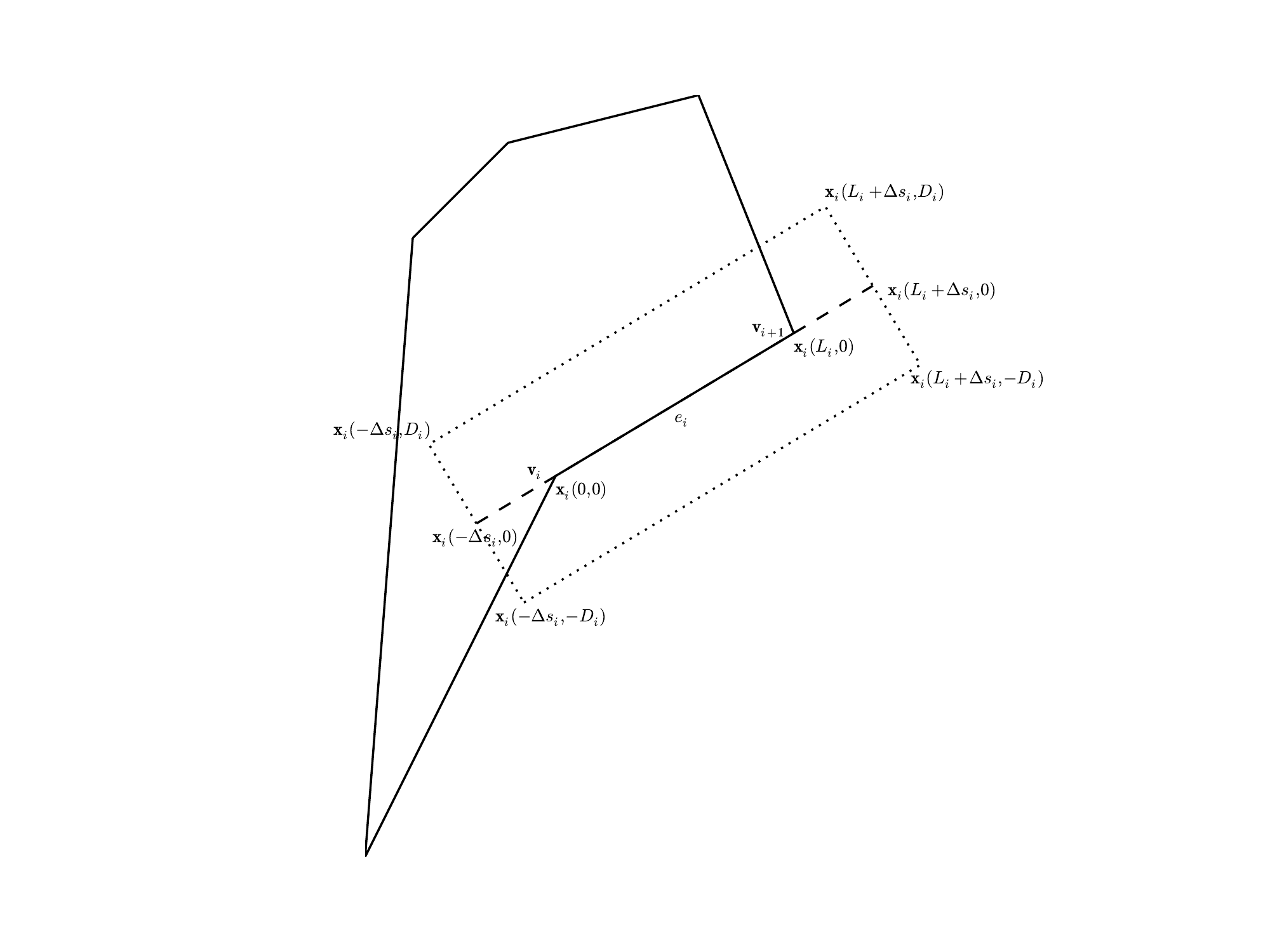}
\caption{Extension algorithm for a simple polygon. The mapping $ \mbfx_i $
defined by the edge $ e_i $ allows to mark the support of the
$ i^{th} $ term in the second sum in (\ref{eq:uPoisson}) which is
the dotted rectangle. The edge $ e_i $ is extended by an amount $ \Delta s_i $.}
\label{fig:extensionBox}
\end{figure}

\paragraph{Poisson problem} We start with the Poisson case,
meaning that boundary data is only given for the function values. 
The biharmonic case will be treated below. The aim is to find
a function $ \tilde{u} $ defined in $ \Omega $ sufficiently smooth,
such that $ \tilde{u} $ satisfies the boundary data:
\be
\tilde{u} = f \quad \mbox{on} \quad \partial \Omega  
\label{eq:PoissonExtension}\\
\ee
As sketched in figures \ref{fig:polygon} and \ref{fig:extensionBox}, 
the boundary of the domain $ \Omega $
consists of an $ n $ sided simple polygon. Each
edge $ e_i $ is connected by adjacent vertices $ \mathbf{v}_i $
and $ \mathbf{v}_{i+1} $ with components
\be
\mbfv_i= \left( { u_i \atop v_i } \right).
\ee
We assume that none of the vertices $ \mbfv_i $ is degenerated. 
Each edge $ e_i $ defines a coordinate transformation $ \mbfx_i(\alpha,\beta) $ by:
\bea
\mbfx_i(\alpha,\beta) &=& \left( \begin{array}{c} x_i(\alpha,\beta) \\ y_i(\alpha,\beta) \end{array} \right) \\
&=& \mbfv_i + \mbft_i \alpha + \mbfn_i \beta \\
&=& \left( { u_i \atop v_i } \right) +
\left( { t_x^i \atop t_y^i } \right) \alpha + 
\left( { n_x^i \atop n_y^i } \right) \beta,
\label{eq:map1}
\eea
where $ \mbft_i $ and $ \mbfn_i $ are the tangential and normal
vectors of the edge $ e_i $:
\bea
\mbft_i &=& \frac{1}{L_i}
\left( { { u_{i+1} - u_{i} } \atop {v_{i+1} - v_i } } \right), \\
\mbfn_i &=& \left( { t^i_y \atop -t_x^i } \right), \\
L_i &=& \sqrt{ (u_{i+1} - u_{i} )^2 + (v_{i+1} - v_i)^2 }.
\eea
For the components of the tangential and normal vector, the index denoting
the edge is written as a superscript in order to simplify the notation. 
The inverse transformation of $ \mbfx_i $ is given by
\bea
\pmb{\alpha}_i(x,y) &=& \left( {\alpha_i(x,y) \atop \beta_i(x,y)  } \right) \\
& =& \left( \begin{array}{cc} t_x^i & t_y^i \\ t_y^i & - t_x^i
\end{array} \right) \left( { { x - u_i } \atop { y-v_i} } \right).
\eea
Each
edge is parametrized by the arclength $ s $ going from $ 0 $ to $ L_i $,
the total arclength of the edge. 
Thus for any point $ (x,y) $ on the edge $ e_i $, we have:
\be
\alpha_i(x,y) = s \in [0,L_i] \quad \mbox{and} \quad \beta_i(x,y) = 0,
\quad (x,y) \in e_i.
\ee
We assume in addition, that 
each edge $ e_i $ can be extended by an amount $ \Delta s_i $,
such that the extended edge $ \tilde{e}_i $ controlled by the parameter $ 
s \in [ -\Delta s_i , L_i + \Delta s_i ] $ has, as $ e_i $, only
two intersections with the other edges, namely
$ \mathbf{v}_i $ and $ \mathbf{v}_{i+1} $, cf. the
dashed line in figure \ref{fig:extensionBox}. 

Since $ \mbfn_i $ is pointing outward of the domain, we can define
a weight $ \omega_i(x,y) $ for each edge by
\be
\omega_i(x,y) = - \beta_i(x,y). \label{eq:edgeweight}
\ee
For a convex polygon $ \omega_i(x,y) $ would be positive in the interior. 
We note
\be
\Omega_i(x,y) = \omega_{i-1}(x,y)\omega_{i+1}(x,y),
\ee
which is zero at the edges $ e_{i-1} $ and $ e_{i+1} $ and nonzero
at $ e_i $ for $ s \in ] 0 , L_i [ $, and has simple roots at
$ s = 0 $ and $ s = L_i $. 
For each vertex $ \mathbf{v}_i $ we define a radial hat function 
$ h_i(\mathbf{x}) $ given by:
\be
h_i(\mathbf{x}) = \exp - \frac{|\mathbf{x} - \mathbf{v}_i|^2}{R_i^2},
\ee
which is a Gaussian with variance $ R_i $. Finally, 
the boundary data $ f $ in equation (\ref{eq:PoissonExtension})
is given by a function $ f_i(s) $ for each edge $ e_i $. We assume
that the boundary data is compatible at the vertices, in the sense that
\be
f_i(L_i) = f_{i+1}(0).  \label{eq:compatibility1}
\ee
The extension function $ \tilde{u} $ in the Poisson case, equation
(\ref{eq:PoissonExtension}), is then written as:
\be
\tilde{u}(x,y) = \sum \limits_{i=1}^n a_i h_i(x,y)
+ \sum \limits_{i=1}^n \phi_i(\alpha_i(x,y)) \Omega_i(x,y) r^{p+1}_i( \beta_i(x,y)),
\label{eq:uPoisson}
\ee
where $ r_i $ is a shape function defined by
\be
r_i (x) = \left\{ \begin{array}{ll} \frac{(D_i - x)(x + D_i)}{D_i^2} & 
-D_i \le x \le D_i \\
0 & \mbox{otherwise} \end{array} \right.
\ee
The length $ D_i $ must be chosen such that the quadrilateral defined by
the vertices:
\be
\mbfx_i( - \Delta s_i , - D_i ) , \, 
\mbfx_i( - \Delta s_i , D_i ) , \, 
\mbfx_i( L_i + \Delta s_i ,  D_i ) , \, 
\mbfx_i( L_i + \Delta s_i ,  - D_i )  \label{eq:cornersBox}
\ee
does only have intersections with the edges $ e_{i-1} $ and $ e_{i+1} $.
We remark that the function $ r^{p+1}_i(\beta_i(x,y))$ is of class $ \mathcal{C}^p $ in $ \Omega $. 
The coefficients $ a_i $ and functions $ \phi_i(s) $ are unknown a priori, but
we require that $ \phi_i $ vanishes outside the interval $ [-\Delta s_i,L_i+\Delta s_i] $:
\be
\phi_i(s) = 0 \quad \mbox{for} \quad s < -\Delta s_i \quad \mbox{and} \quad s > L_i + \Delta s_i. 
\ee
In figure \ref{fig:extensionBox}, we plotted the extended edge $ \tilde{e}_i $
by a dashed line. On this line we marked the end points of $ \tilde{e}_i $,
i.e. $ \mbfx_i(-\Delta s_i,0) $ and $ \mbfx_i(L_i+\Delta s_i,0) $,
by means of the mapping $ \mbfx_i $, equation (\ref{eq:map1}),
defined by the edge $ e_i $. The points $ \mbfx_i(0,0) $ 
and $ \mbfx_i(L_i,0) $ return the vertices $ \mbfv_i $ and $ \mbfv_{i+1} $,
respectively. Due to the finite support of $ \phi_i $ and $ r_i $, the
$ i^{th} $ term in the second sum of equation (\ref{eq:uPoisson}) has
a finite support given by the dotted box in figure \ref{fig:extensionBox}
with corners given in equation (\ref{eq:cornersBox}). \\

In order to determine the coefficients $ a_i $, we evaluate $ \tilde{u} $
at the vertices $ \mbfv_j $:
\be
\tilde{u}(\mbfv_j) = \sum \limits_{i=1}^n a_i h_i(\mbfv_j),
\quad  j = 1,\ldots, n,
\ee
since the second sum in (\ref{eq:uPoisson}) vanishes at the vertices. 
The function value $ \tilde{u}(\mbfv_j) $ must equal $ f_j(0)$ by definition,
which gives us a system of $ n $ equations for the $ n $ unknown $ a_i$. 
In order to determine the unknown functions $ \phi_i $, 
we evaluate $ \tilde{u} $ at a point $ (x,y) $ on the edge $ e_j $:
\be
\tilde{u}(x,y) = \sum \limits_{i=1}^n a_i h_i(x,y) + \phi_j(\alpha_j(x,y))
\Omega_j(x,y).
\ee
A first choice $ \tilde{\phi}_j(s) $ for $ \phi_j(s) $ would be
\be
\tilde{\phi}_j(s) = \frac{ f_j(s) - \sum \limits_{i=1}^n a_i h_i( x_j(s,0),y_j(s,0)) }
{\Omega_j(x_j(s,0),y_j(s,0) ) },
\ee
which is well defined for $ s = 0 $ and $ s = L_i $, since the numerator
goes equally fast (or faster) to zero as the denominator. If the boundary data
 $ f_j $ is of class $ \mathcal{C}^p([0,L_i]) $, so is $ \tilde{\phi}_j $.
However, the behavior
of $ \tilde{\phi}_j $ outside of the interval $ [ 0, L_i ] $ is not controlled. 
We shall therefore choose a function $ \phi_j(s) $ which is sufficiently smooth
on the entire $ \mathbb{R} $ axis and vanishes outside of 
the interval $ [-\Delta s_i, L_i + \Delta s_i ] $. This function $ \phi_j(s) $
can,
however, only be an approximation to $ \tilde{\phi}_j $. 
The extension $ \tilde{u} $
will therefore not be an exact extension but only an approximate extension,
which
is, however, not a problem, as long as the approximation is sufficiently
close in order not to reduce the accuracy of the overall solution. We use 
a spline interpolation of degree $ p+1 $ for $ \phi_j$,
where we assume that $ p $ 
is even for simplicity. For each edge $ e_j $, we choose a number $ N_j $ and
define a sequence of knots $ \xi_k$  with:
\be
\xi_{-(p+1)} = -\Delta s_j , \quad \xi_0 = 0, \quad \xi_{N_j}  = L_j, \quad \xi_{N_j+(p+1)} = L_j + \Delta s_j.
\ee
The position of the other knots may be chosen arbitrarily
as long as the sequence $ \xi_k $ is
monotonically increasing. In the present treatise, we choose a uniform 
distribution of knots $ \xi_k $. The interpolant $ \phi_j $ is
thus a linear combination of $ N_j + p +1 $ B-splines:
\be
\phi_j = \sum \limits_{k=-{p+1}}^{N_j-1} c_k b_k^{p+1}(s).
\ee
The $ N_j + p +1 $ conditions allowing to determine the $ c_k $ are given by
\bea
\phi_j(\xi_k) & = & \tilde{\phi}_j(\xi_k), \quad k = 0,\dots,N_j\\
\frac{d^q}{ds^q} \phi_j( \xi_0 ) & = & \frac{d^q}{ds^q} \tilde{\phi}_j(\xi_0), \quad q = 1,\ldots, p/2\\
\frac{d^q}{ds^q} \phi_j( \xi_N ) & = & \frac{d^q}{ds^q} \tilde{\phi}_j(\xi_N), \quad q = 1,\ldots, p/2
\eea
Finally, the extension $ \tilde{u} $ defined by the above procedure
will be of class $ \mathcal{C}^p $ inside the polygon $ \Omega $
and on the boundary
$ \partial \Omega $ if the boundary data is sufficiently smooth. \\

The biharmonic case is more involved. In this case not only the 
function value is imposed on the boundary, but also its normal derivative:
\bea
\tilde{u} &=& f \quad \mbox{on} \quad \partial \Omega \\
\frac{\partial}{\partial n} \tilde{u} &=& g \quad \mbox{on} \quad \partial \Omega 
\eea
In addition to an $ f_i $ for each edge, we are given a function $ g_i $
matching the normal derivative of $ \tilde{u} $ at $ e_i $. Next
to the compatibility condition (\ref{eq:compatibility1}),
constraints on the first
and second derivative arise. For the first derivative we can write:
\bea
\left( { \frac{d}{ds} f_i(L_i) \atop \frac{d}{ds} f_{i+1} (0) } \right)& =&
\left( \begin{array}{cc} t_x^i & t_y^i \\ t^{i+1}_x & t_y^{i+1} \end{array}
\right) \left( { \frac{ \partial }{\partial x} \tilde{u}(\mbfv_{i+1} ) 
\atop  \frac{ \partial }{\partial y} \tilde{u}(\mbfv_{i+1} ) } \right) 
\label{eq:compatibilitySystem1}\\
\left( { g_i(L_i) \atop g_{i+1} (0) } \right)& =&
\left( \begin{array}{cc} n_x^i & n_y^i \\ n^{i+1}_x & n_y^{i+1} \end{array}
\right) \left( { \frac{ \partial }{\partial x} \tilde{u}(\mbfv_{i+1} ) 
\atop  \frac{ \partial }{\partial y} \tilde{u}(\mbfv_{i+1} ) } \right) 
\label{eq:compatibilitySystem2}
\eea
Inverting (\ref{eq:compatibilitySystem1}) and (\ref{eq:compatibilitySystem2})
leads to two additional compatibility conditions to be fulfilled by
$ f $ and $ g $:
\bea
t_y^{i+1} \frac{d}{ds} f_i(L_i) - t_y^i \frac{d}{ds} f_{i+1}(0) &=&
- t_x^{i+1} g_i(L_i) + t_x^i g_{i+1}(0) \\
-t_x^{i+1} \frac{d}{ds} f_i(L_i) + t_x^i \frac{d}{ds} f_{i+1}(0) &=&
- t_y^{i+1} g_i(L_i) + t_y^i g_{i+1}(0) 
\eea
A fourth compatibility condition arises from the fact
that four derivatives are imposed at $ \mbfv_{i+1} $,
namely 
\be
\frac{d^2}{ds^2} f_i(L_i), \quad \frac{d^2}{ds^2} f_{i+1}(0),
\quad \frac{d}{ds} g_i(L_i), \quad \frac{d}{ds} g_{i+1}(0),
\ee
but only three unknowns are available, which are
\be
\frac{\partial^2}{\partial x^2} \tilde{u}(\mbfv_{i+1}), \quad
\frac{\partial^2}{\partial y^2} \tilde{u}(\mbfv_{i+1}), \quad
\mbox{and} \quad
\frac{\partial^2}{\partial x \partial y} \tilde{u}(\mbfv_{i+1}). 
\ee
After some algebra the compatibility equation can be obtained as:
\be
\mbft_i \cdot \mbfn_{i+1}  \left( \frac{d^2}{ds^2} f_i(L_i) + \frac{d^2}{ds^2} f_{i+1}(0) \right) + \mbft_i \cdot \mbft_{i+1}
\left( - \frac{d}{ds} g_i(L_i) + \frac{d}{ds} g_{i+1}(0) \right) = 0
\ee
Given compatible $ f $ and $ g $, we shall, 
as for the Poisson case, consider the vertices first,
before turning to the edges. However, for 
the biharmonic case, we need to impose in addition to the function
value $ \tilde{u} $, the first and second order derivatives of 
$ \tilde{u} $ at the vertices. This is done by choosing the 
following expansion for $ \tilde{u} $:
\bea
\tilde{u} &=& \sum \limits_{i=1}^n H_i(x,y) \nonumber \\
&+& \sum \limits_{i=1}^n \phi_i(\alpha_i(x,y) ) \Omega_i^3(x,y) r_i^{p+1}
\left( \beta_i(x,y) \right) \nonumber \\
&+& \sum \limits_{i=1}^n \psi_i(\alpha_i(x,y) ) \Omega_i^2(x,y) \omega_i(x,y)r_i^{p+1} \left( \beta_i(x,y) \right) \label{eq:utildeHermite}
\eea
where the function $ H_i(x,y) $ collects all the coefficients for
the conditions on the vertices:
\be
H_i = a_i h_i + b_i \omega_i h_i + c_i \omega_{i+1} h_i 
+ d_i \omega_{i}^2 h_i + e_i \omega_{i+1}^2 h_i + k_i \omega_{i} \omega_{i+1} h_i.
\ee
The weighting of $ \phi_i $ and $ \psi_i $ by $ \Omega_i^3 $ and $ \Omega_i^2 \omega_i $, respectively, is such that for derivatives up to the third order
the second and third sum in equation (\ref{eq:utildeHermite}) vanish
at the vertices. For this reason the 6 coefficients, 
$ a_i $, $ b_i $, $ c_i $, $ d_i $, $ e_i $,
and $ k_i $ in $ H_i $, 
can be determined by 6 conditions at the vertices $ \mbfv_j $,
$ j = 1,\ldots,n$:
\bea
\sum \limits_{i=1}^n H_i(\mbfv_j) &=& f_j(0) \\
\sum \limits_{i=1}^n \mbft_j \cdot \nabla H_i (\mbfv_j) & = & \frac{d}{ds} f_j(0)\\
\sum \limits_{i=1}^n \mbfn_j \cdot \nabla H_i (\mbfv_j) & = & g_j(0)\\
\sum \limits_{i=1}^n \mbft_j \cdot \left( \nabla \nabla H_i (\mbfv_j) \right) 
\cdot \mbft_j & = & \frac{d^2}{ds^2} f_j(0) \\
\sum \limits_{i=1}^n \mbfn_j \cdot \left( \nabla \nabla H_i (\mbfv_j) \right) 
\cdot \mbft_j & = & \frac{d}{ds} g_j(0) \\
\sum \limits_{i=1}^n \mbft_{j-1} \cdot \left( \nabla \nabla H_i (\mbfv_j) \right) 
\cdot \mbft_{j-1} & = & \frac{d^2}{ds^2} f_{j-1}(L_{j-1}) 
\eea
For a point $ (x,y) $ on the edge $ e_j $ the expansion $ \tilde{u}$, equation
(\ref{eq:utildeHermite}), becomes:
\be
\tilde{u}(x,y) = \sum \limits_{i=1}^n H_i(x,y)
+ \phi_j(\alpha_j(x,y) ) \Omega_j^3(x,y) 
\ee
Similarly as before, we define a function $ \tilde{\phi}_j(s) $ by
\be
\tilde{\phi}_j(s) = 
\frac{ f_j(s) - \sum \limits_{i=1}^n H_i(x_j(s,0),y_j(s,0))}{\Omega_j^3(x_j(s,0),y_j(s,0))},
\ee
which is well posed in the interval $ [0,L_j] $, since the cubic roots in
the numerator at $ s = 0 $ and $ s = L_j $ are balanced by roots of equal
or higher degree in the numerator. As for the Poisson case, a B-spline interpolation
of order $ p + 1 $ is used to interpolate $ \tilde{\phi}_j $.
The normal derivative of the expansion $ \tilde{u} $ at a point $ (x,y) $ on
the edge $ e_j $ can then be written:
\bea
\mbfn_j \cdot \nabla \tilde{u} ( x,,y) &=& \sum \limits_{i=1}^n \mbfn_j \cdot \nabla H_i(x,y) \\ &+&
 \mbfn_j \cdot \nabla \left\{ \phi_j (\alpha_j(x,y) ) \Omega_j^3(x,y) \right\}\\
&-& \psi_j(\alpha_j(x,y)) \Omega_j^2(x,y) 
\eea
A function $ \tilde{\psi}_j(s) $ is defined as
\bea
\tilde{\psi}_j(s) &=& \frac{1}{ \Omega_j^2(x_j(s,0),y_j(s,0)) } \Big[ g_j(s) - \sum \limits_{i=1}^n \mbfn_j 
\cdot \nabla H_i( x_i(s,0) , y_i(s,0) )  \nonumber \\
& -& \mbfn_j \cdot \nabla \left\{ \phi_j (\alpha_j(x_j(s,0),y_j(s,0)) ) \Omega_j^3(x_j(s,0),y_j(s,0)) \right\} \Big], \nonumber
\eea
which is as before by construction well defined for $ s = 0 $ and $ s =  L_j $. 
The B-spline interpolation $ \psi_j(s) $ of $ \tilde{\psi}_j(s) $ 
completes the present boundary data extension method. 

\section{Vertices of polygons} \label{sec:coordinatesPolygons}

The vertices of the convex polygon in section \ref{sec:polygonal}
are given by:
\bea
\mbfv_0 &=& (2,-4) \\
\mbfv_1 &=& (5,3) \\
\mbfv_2 &=& (3,8) \\
\mbfv_3 &=& (-3.5,6.) \\
\mbfv_4 &=& (-5,-8).
\eea
For the non convex polygon the definition of the vertices is
exactly the same apart from $ \mbfv_0 $, which is given by:
\be
\mbfv_0 = (0.5,0).  
\ee
The holes are for both cases given by circles with centers
\be
\mbfx_a = (-1.7, 0.7 ) \quad \mbox{and} \quad \mbfx_b = (1.3,4.05), 
\ee
respectively. The radii are $ 0.91 $ and $ 1.1 $, respectively. 

\end{appendix}


\end{document}